\documentclass[12pt, reqno]{amsart}
\usepackage{xcolor}
\usepackage{blindtext}
\usepackage{latexsym}
\usepackage[pagewise]{lineno}
\usepackage[utf8]{inputenc}

\usepackage{exscale}
\usepackage{amsmath}
\usepackage{amssymb}
\usepackage{amsfonts}
\usepackage{mathrsfs}
\usepackage{amsbsy}
\usepackage{relsize}

\usepackage[colorlinks, linkcolor=venetianred, citecolor=darkpowderblue, urlcolor=blue, hypertexnames=false]{hyperref}

\usepackage{comment}
\PassOptionsToPackage{reqno}{amsmath}

\usepackage{esint} 
\usepackage{amssymb} 
\usepackage{stmaryrd}
\usepackage{cite}
\calclayout

\usepackage{calc}    

\definecolor{falured}{rgb}{0.5, 0.09, 0.09}
\definecolor{darkpowderblue}{rgb}{0.0, 0.2, 0.6}
\definecolor{darkviolet}{rgb}{0.58, 0.0, 0.83}
\definecolor{cornellred}{rgb}{0.858, 0.188, 0.478}
\definecolor{mypink2}{RGB}{219, 48, 122}
\definecolor{mypink3}{cmyk}{0, 0.7808, 0.4429, 0.1412}
\definecolor{mygray}{gray}{0.6}
\definecolor{venetianred}{rgb}{0.78, 0.03, 0.08}
\definecolor{sapphire}{rgb}{0.03, 0.15, 0.4}
\definecolor{utahcrimson}{rgb}{0.83, 0.0, 0.25}
\definecolor{trueblue}{rgb}{0.0, 0.45, 0.81}
\definecolor{carminered}{rgb}{1.0, 0.0, 0.22}
\definecolor{cobalt}{rgb}{0.0, 0.28, 0.67}
\definecolor{cornflowerblue}{rgb}{0.39, 0.58, 0.93}
\definecolor{cadmiumgreen}{rgb}{0.0, 0.42, 0.24}
\definecolor{cornellred}{rgb}{0.7, 0.11, 0.11}

\newtheorem{thm}{Theorem}[section]
\newtheorem{prop}{Proposition}[section]
\newtheorem{defi}{Definition}[section]
\newtheorem{lem}{Lemma}[section]
\newtheorem{cor}{Corollary}[section]

\newtheorem{rem}{Remark}[section]

\newcommand{\R}{\mathbb{R}}

\numberwithin{equation}{section}

\newcommand{\eps}{\epsilon}

\newcommand{\wto}{\rightharpoonup}

\makeatletter \@addtoreset{equation}{section} \makeatother

\newcounter{const}
\definecolor{green}{rgb}{0.25,0.75,0}

\author[ T. Gou, M. Majdoub \& T. Saanouni]{Tianxiang Gou, Mohamed Majdoub and Tarek Saanouni}
\address[T. Gou]{School of Mathematics and Statistics, Xi'an Jiaotong University, Xi'an 710049, Shaanxi, People's Republic of China.}
\email{\tt {tianxiang.gou@xjtu.edu.cn}}
\address[M. Majdoub]{Department of Mathematics, College of Science, Imam Abdulrahman Bin Faisal University, P. O. Box 1982, Dammam, Saudi Arabia.\newline 
Basic and Applied Scientific Research Center, Imam Abdulrahman Bin Faisal University, P.O. Box 1982, 31441, Dammam, Saudi Arabia.}
\email{\tt {mmajdoub@iau.edu.sa}}
\email{\tt {med.majdoub@gmail.com}}
\address[T. Saanouni]{Departement of Mathematics, College of Sciences and Arts of Uglat Asugour, Qassim University, Buraydah, Kingdom of Saudi Arabia.}
\email{\tt {T.saanouni@qu.edu.sa}}



\title[Competing INLS]{NLS equation with competing inhomogeneous
nonlinearities: ground states, blow-up, and scattering}

\begin{document}

\begin{abstract} 
We investigate a class of  nonlinear equations of Schr\"odinger type  with competing inhomogeneous nonlinearities in the non-radial inter-critical regime,
\begin{align*}
\textnormal{i} \partial_t u +\Delta u &=|x|^{-b_1} |u|^{p_1-2} u - |x|^{-b_2} |u|^{p_2-2}u \quad \mbox{in} \,\, \R \times \R^N,
\end{align*}
where $N \geq 1$, $b_1, b_2>0$ and $p_1,p_2>2$.
 First, we establish the existence/nonexistence, symmetry, decay, uniqueness, non-degeneracy and instability of ground states. Then, we prove the scattering versus blowup  below the ground state energy threshold. Our approach relies on Tao's scattering criterion and Dodson-Murphy's Virial/Morawetz inequalities. We also obtain an upper bound of the blow-up rate.
The novelty here is that the equation does not enjoy any scaling invariance due to the presence of competing { nonlinearities} and the singular weights prevent the invariance by translation in the space variable.

 To the best of authors knowledge, this is the first time when { inhomegeneous} NLS equation with a focusing leading order nonlinearity and a defocusing perturbation is investigated.
\end{abstract}
 
\maketitle

\thispagestyle{empty}

\section{Introduction}

In this paper, we are concerned with solutions to the following nonlinear Schr\"odinger equation with competing inhomogeneous nonlinearities,
\begin{subequations}\label{equ}
\begin{align}
\textnormal{i} \partial_t u +\Delta u &=|x|^{-b_1} |u|^{p_1-2} u - |x|^{-b_2} |u|^{p_2-2}u, \quad\mbox{for}\;\;\; (t,x)\in\R \times \R^N,\label{equ1-1}\\
u(0,x) &= u_0(x), \quad\quad\quad \mbox{for}\;\;\;\; x\in \R^N, \label{equ1-2}
\end{align}
\end{subequations}
where $N \geq 1$, $b_1, b_2>0$ and $p_1,p_2>2$. Equation \eqref{equ1-1} can be used to describe many physical phenomena, we refer the readers to \cite{BPVT, KA} and references therein for further interpretations. The well-known nonlinear Schr\"odinger equation
\begin{align}\label{sequ}
\textnormal{i} \partial_t u +\Delta u = |u|^{p-2} u \quad \mbox{in} \,\, \R \times \R^N
\end{align}
for $N \geq 1$ and $2<p<\frac{2N}{(N-2)^+}$ arises in various physical contexts, for example in nonlinear optics as well as in the description of nonlinear waves such as the propagation of laser beams, water waves at the free surface of an ideal fluid and plasma waves. In particular, equation \eqref{sequ} models the propagation of intense laser beams in a homogeneous bulk medium with Kerr nonlinearity. It was suggested that stable high power propagation can be achieved in plasma by sending a preliminary laser beam that creates a channel with a reduced electron density, and thus reduces the nonlinearity inside the channel, see for instance \cite{Gi, LT}. Under these conditions, the beam propagation can be modeled in the simplest case by the following inhomogeneous nonlinear Schr\"odinger equation 
\begin{align} \label{kequ}
\textnormal{i} \partial_t u +\Delta u = K(x) |u|^{p-2} u \quad \mbox{in} \,\, \R \times \R^N,
\end{align}
where $u$ is the electric field in laser and optics, $2<p<\frac{2N}{(N-2)^+}$ is the power of nonlinear interaction and the potential $K(x)$ is proportional to the electron density.

Equation \eqref{kequ} has attracted much interest from mathematical point-view. When $K$ is constant, then \eqref{kequ} is reduced to the classical nonlinear Schr"odinger equations such as \eqref{sequ}. Such problems have been extensively studied in the last decades. Local well-posedness for \eqref{sequ} in the energy space $H^1$ was first established by Ginibre and Velo in \cite{GV}. The existence of finite-time blow-up solutions was early proved by Glassey in \cite{Gl}. For the case $p=2+\frac 4N$, Weinstein investigated the structure and formation of the singularity of solutions in \cite{We}. Successively, the concentration phenomenon of blowup solutions was considered by Merle and Tsutsumi in \cite{MT}. The exact blowup solutions with the critical mass was constructed by Merle in \cite{M1}. We also refer the readers to \cite{MR1, MR2, MR3} for further consideration of blow-up solutions in the mass critical case.  For the case $p>2+\frac{4}{N}$, scattering versus blowup of solutions below the ground state energy level were investigated in \cite{ADM, DM1, DM, DHR, FXC} under the energy subcritical case and in \cite{HR, KM, KV} under the energy critical case. We also refer the readers to the monographs \cite{Ca, SS, Ta} for more relevant topics in these directions.

When $K$ is bounded, Merle in \cite{Me} first discussed the existence and nonexistence of minimal blowup solutions to  \eqref{kequ} with $p=2+\frac{4}{N}$. Furthermore, Rapha\"el and Szeftel in \cite{MRS} studied the existence, uniqueness and characterization of minimal blowup solutions. Later on, Fibich and Wang in \cite{FW} and Liu et al. in \cite{LWW} investigated stability and instability of standing waves to \eqref{kequ} with $p \geq 2+\frac{4}{N}$.

When $K$ is unbounded, the problem becomes more involved. For the case $K(x)=|x|^{b}$ with $b>0$, sharp conditions on the global existence and blowup of solutions to \eqref{kequ} were established by Chen in \cite{Ch} as well as Chen and Guo in \cite{CG}. Afterwards, Zhu in \cite{Zh} derived the existence and concentration phenomenon of blowup solutions. Recently, Dinh et al. in \cite{DMS} further investigated scattering and blowup of solutions. For the case $K(x)=|x|^{-b}$ with $b>0$, the problem has received much attentions in latest years. The well-posedness for \eqref{kequ} was initially established in \cite{AT, GS, Gu}. Later on, the existence and dynamical properties of blowup solutions was revealed in \cite{CG', D, Fa, G} for $p \geq 2+ \frac{2(2-b)}{N}$. Sharp thresholds for scattering versus blowup of solutions below the ground stat energy level were attained in \cite{Cam, CC, CFGM, DK, FG, MMZ} for the energy sub-critical regime and in \cite{CHL, CL, GM} for the energy critical case.

Following the pioneer article of Tao et al. \cite{TVZ}, there are many works devoted to the study of solutions to nonlinear Schr\"odinger equations with combined power nonlinearities. Indeed, the scattering and blowup of solutions for an intercritical focusing perturbed equation was investigated by Bellazzini et al. \cite{BDF} for a defocusing perturbation and Xie \cite{Xie} for a focusing perturbation. In three space dimensions, the scattering versus blow-up dichotomy under the ground state threshold was considered in \cite{KOPV, KRV} for a defocusing energy-critical perturbed source term. The same questions were investigated for a focusing energy-critical perturbed non-linearity in three space dimensions in \cite{MXZ}. This result was extended to four space dimensions in \cite{MTZ} and to lower dimensions in \cite{CMZ}. In the inhomogeneous case, some well-posedness issues {and dynamical behaviors of solutions to} a Schr\"odinger equation with a combined source term were investigated in \cite{hms, Gou}. Inspired by these literature, it would be interesting to investigate solutions to \eqref{equ} having competing inhomogeneous nonlinearities. It is worth mentioning that the problem under our consideration loses any scaling invariance due to the competing nonlinearities. In addition, singular weights prevent any space translation. This is different from the problems treated previously, which in fact leads to difficulties in our study.

{
Unless otherwise specified, along  the rest of this article, we suppose that
\begin{equation}
\left\{\begin{array}{l}
 N \geq 1, 0<b_1, b_2< \min \{2, N\}, \\  2<\frac{N(p_1-2)+2b_1}{2}<\frac{N(p_2-2)+2b_2}{2},\\
2<p_1<2_{b_1}^*, 2<p_2<2_{b_2}^*,
\label{H}
\end{array}\right.
\end{equation}
where $2_b^*:=\frac{2(N-b)}{(N-2)^+}$.
It is easy to see that
$$
p_j>2+\frac{2(2-b_j)}{N} \Longleftrightarrow \frac{N(p_j-2)+2b_j}{2}>2, \quad j=1,2.
$$
}
The first objective of the present paper is to investigate standing wave solutions to \eqref{equ1-1}. Recall that a standing wave solution has the form
$$
u(t,x)=e^{\textnormal{i} \omega t} \phi(x), \quad \omega \in \R,
$$
where $\phi\in H^1(\R^N)$ solves the elliptic equation
\begin{align} \label{equ1}
-\Delta \phi + \omega \phi =|x|^{-b_2} |\phi|^{p_2-2}\phi-|x|^{-b_1} |\phi|^{p_1-2} \phi \quad \mbox{in} \,\, \R^N.
\end{align}
\begin{defi} [Ground state]
		A non-trivial $H^1$ solution $Q$ to \eqref{equ1} is called a ground state related to \eqref{equ1} if it minimizes the action functional
		\[
		S_{\omega}(\phi):=\frac 12 \int_{\R^N} |\nabla \phi|^2 \,dx + \frac {\omega}{2} \int_{\R^N} |\phi|^2 \,dx +\frac{1}{p_1} \int_{\R^N} |x|^{-b_1} |\phi|^{p_1} \,dx - \frac{1}{p_2} \int_{\R^N} |x|^{-b_2} |\phi|^{p_2} \,dx.
		\]
		over all non-trivial solution of \eqref{equ1}, that is,
		$$
		S_\omega(Q) = \inf \Big\{ S_\omega(\phi) \ : \ \phi \in H^1 \backslash \{0\}\;\; \text{solves } \eqref{equ1} \Big\}.
		$$
	\end{defi}

Actually, we will show that
$$
S_{\omega}(Q) = \inf\Big\{ S_{\omega}(\phi);\quad  \phi \in H^1\backslash\{0\},\; K(\phi)=0 \Big\},
$$

  where
\begin{align} \label{K-phi}
\begin{split}
K(\phi):&=\int_{\R^N} |\nabla \phi|^2 \,dx+\frac{N(p_1-2)+2b_1}{2p_1}\int_{\R^N}|x|^{-b_1}|\phi|^{p_1}\, dx\\
& \quad -\frac{N(p_2-2)+2b_2}{2p_2}\int_{\R^N}|x|^{-b_2}|\phi|^{p_2}\, dx.
\end{split}
 \end{align}
To study solutions to \eqref{equ1}, we shall make use of variational methods. In this sense, solutions to \eqref{equ1} correspond to critical points of the underlying energy functional $S_{\omega}$ in $H^1(\R^N)$.

It is standard to check that $S_{\omega}$ is of class $C^1$ in $H^1(\R^N)$.  In the variational framework, we shall establish the existence/nonexistence, symmetry, decay, uniqueness, non-degeneracy and instability of solutions to \eqref{equ1}. The second aim of the present paper is to consider scattering versus blowup of solutions below the energy threshold to the Cauchy problem \eqref{equ}.

Although the equation \eqref{equ1-1} doesn't enjoy any scaling invariance, we introduce here a useful scaling in the space variable
\begin{equation}
\label{Scal}
v_{\lambda}(x):=\lambda^{\frac{N}{2}}v(\lambda x), \quad x \in \R^N, \lambda>0.
\end{equation}
In particular, one can easily verify that
$$
K(\phi)=\partial_\lambda\Big(S_\omega(\phi_\lambda)\Big)_{|\lambda=1}.
$$
We denote by $H_{rad}^1(\R^N)$ the subspace of $H^1(\R^N)$ consisting of all radially symmetric functions in $H^1(\R^N)$.

\begin{thm} \label{thm1}
Let \eqref{H} hold. Then
\begin{itemize}
\item[$(\textnormal{i})$]  There exist positive, radially symmetric and decreasing ground states in $H^1(\R^N)$ to \eqref{equ1} for any $\omega>0$.
\item[$(\textnormal{ii})$] Moreover, assume that $b_1<b_2$, $p_1<p_2$. If $3 \leq N \leq 4$, then there exist no positive, radially symmetric, and decreasing solutions in $H^1(\R^N)$ to \eqref{equ1} for $\omega=0$. If $N \geq 5$, then there exist positive, radially symmetric, and decreasing ground states in $H^1(\R^N)$ to \eqref{equ1} for $\omega=0$.
\item[$(\textnormal{iii})$]  If $N \geq 2$, then there are no solutions in $H^1_{rad}(\R^N)$ to \eqref{equ1} for any $\omega<0$.
\end{itemize}
\end{thm}

When $\omega>0$, by introducing the following minimization problem with the help of the underlying Pohozaev manifold 
$$
m_{\omega}:=\inf_{u \in P} S_{\omega}(u), \quad P:= \Big\{\phi \in H^1(\R^N) \backslash \{0\} : K(\phi) =0\Big\}
$$ 
and utilizing the fact that $H^1(\R^N)$ is compactly embedded in $L^p(\R^N, |x|^{-b}dx)$ for any $2<p<2^*_b$, see lemma \ref{cembedding}, one can establish the existence of ground states in \eqref{equ1}. It is worth remaking that the Pohozaev manifold $P$ is a natural constraint and minimizers of the minimization problem are ground states to \eqref{equ}. While $\omega=0$, the existence of solutions becomes delicate. In this situation, one clearly finds that $H^1(\R^N)$ is no longer a natural Sobolev space to seek for solutions. For this, we shall introduce the associated Sobolev space $X$ defined by the completion of $C_0^{\infty}(\R^N)$ under the norm
$$
\|u\|_X:=\left(\int_{\R^N} |\nabla u|^2 \,dx \right)^{\frac 12} + \left(\int_{\R^N} |x|^{-b_1}|u|^{p_1}\,dx\right)^{\frac{1}{p_1}}.
$$
Here the extra assumptions $b_1<b_2$ and $p_1<p_2$ are used to guarantee that $X$ is continuously embedded in $L^{p_2}(\R^N, |x|^{-b_2}dx)$, see Lemma \ref{em}. Similarly, by introducing the following minimization problem with the aid of the underlying Pohozaev manifold, 
$$
m_{0}:=\inf_{u \in \mathcal{P}} S_{0}(u), \quad \mathcal{P}:= \Big\{\phi \in X \backslash \{0\} : K(\phi) =0\Big\},
$$ 
one can also derive the existence of ground states to \eqref{equ1}, because $X$ is compactly embedded in $L^{p_2}(\R^N, |x|^{-b_2}dx)$, see Lemma \ref{xem}. Then, using the optimal decay of the solutions, see Theorem \ref{decay}, one concludes that the solutions in $X$ to \eqref{equ1} for $\omega=0$ are also ones in $H^1(\R^N)$ for $N \geq 5$, because the solutions in $X$ belong to $L^2(\R^N)$ for $N \geq 5$. Furthermore, one can obtain the nonexistence of the solutions for $3 \leq N \leq 4$ by Theorem \ref{decay}, because the solutions do not belong to $L^2(\R^N)$ in this case. Here, to demonstrate radial symmetry and decrease of ground states, we shall adapt the polarization arguments developed in \cite{BWW}. Due to the presence of competing inhomogeneous nonlinearities, one cannot adopt the symmetric-decreasing rearrangement arguments in \cite{Lieb} to show the symmetry of the solutions. To show the existence of positive ground states, we shall apply the maximum principle. When $\omega<0$, one can employ Kato's arguments in \cite{Ka} along with the radial Sobolev embedding results, see Lemma \ref{hrc}, to obtain the nonexistence of solutions in $H^1_{rad}(\R^N)$.

\begin{rem}
{\rm Observe that if $Q$ is a ground state to \eqref{equ1}, then $|Q|$ is also a ground state solution to \eqref{equ1} and we have $\|\nabla |Q|\|_2=\|\nabla Q\|_2$. Hence, by \cite[Theorem 4.1]{HS}, any ground state $Q$ has the form of $Q=e^{\textnormal{i} \theta} |Q|(\cdot-y)$ for $\theta \in \R$ and $y \in \R^N$.}
\end{rem}

Let $\phi \in H^1(\R^N)$ be a positive and radially symmetric solution to \eqref{equ1} for $\omega>0$. When $N \geq 2$, by Lemma \ref{hrc}, there holds that
$$
-\Delta \phi = \left(-\omega + |x|^{-b_2} |\phi|^{p_2-2} - |x|^{-b_1} |\phi|^{p_1-2} \right) \phi \geq -\frac{\omega}{2} \phi, \quad |x| \geq R.
$$
It then follows that $\phi$ enjoys exponential decay at infinity for $N \geq 2$. However, this is different from the case for $\omega=0$, where solutions only possess algebraic decay at infinity, as presented below.

\begin{thm} \label{decay}
Let \eqref{H} hold, $N \geq 3$, $b_1<b_2$ and $p_1<p_2$. Let $\phi \in X$ be a positive, radially symmetric, and decreasing solution to \eqref{equ1} for $\omega=0$. Then there holds that
\begin{align*} 
\phi(x) \underset{|x| \to \infty}{\sim}\left\{
\begin{aligned}
&|x|^{-\beta} \quad &\mbox{if }\,\, p_1\neq (2N-2-b_1)/(N-2),\\
&|x|^{2-N} \left(\ln |x|\right)^{\frac{2-N}{2-b_1}}\quad  &\mbox{if }\,\, p_1 = (2N-2-b_1)/(N-2),
\end{aligned}
\right.
\end{align*}
where $\beta:=\max\{(2-b_1)/(p_1-2), N-2\}>0$.
\end{thm}

Apparently, Theorem \ref{decay} reveals that the solution $\phi \not \in L^2(\R^N)$ for $3 \leq N \leq 4$ and  $\phi \in L^2(\R^N)$ for $N \geq 5$.
\begin{rem}
\rm Let $\phi_{\omega}$ be the solution obtained in Theorem \ref{thm1} for $\omega \geq 0$. Then, by Theorem \ref{decay}, it is not difficult to show that $\phi_{\omega} \to \phi_0$ in $X$ for $3 \leq N \leq 4$ and in $H^1(\R^N)$ for $N \geq 5$ as $\omega \to 0$. 
\end{rem}

\begin{thm} \label{thmunique}
Let \eqref{H} hold, $N \geq 3$, $b_1<b_2$ and $p_1<p_2$. Moreover, assume that
\begin{align} \label{uniquec}
\frac{2(N-1)-b_1}{p_1+2}-\frac{2b_1+2(N-1)p_1}{p_2(p_1+2)}+\frac{b_2}{p_2} \leq 0.
\end{align}
Then there exists at most one positive, radially symmetric and decreasing solution to \eqref{equ1}.
\end{thm}

The study of the uniqueness of solutions to nonlinear elliptic equations has a very long history, and there are two important schemes which turn out to be powerful. The first one is based on shooting arguments concerning the analysis of the set of zeros and global behaviors of solutions to the associated ordinary differential equations by Sturm's oscillation theory. This goes back to Coffman and Kwong, who initially proved the uniqueness of solutions to equations with power type nonlinearities in \cite{Co, Kw}. Let us also refer to the readers to \cite{Ge, Ya} for the study of the uniqueness of solutions to equations with inhomogeneous nonlinearities in the same spirit. The second, which is based on the Pohozaev identity, was originally proposed by Yanagida in \cite{Y}. Later, it was further generalized by Shioji and Watanabe in \cite{SW1, SW}. Since the problem under our consideration possesses competing inhomogeneous nonlinearities, the shooting arguments are not available to discuss the uniqueness of solutions. In particular, one can check that the associated local energy does not decrease along the trajectories in our case. For this reason, to prove Theorem \ref{thmunique}, we shall take advantage of the second scheme. Generally speaking, we first need to calculate the corresponding Pohozaev identity satisfied by the solutions to the ordinary differential equation \eqref{ode1}. Equation \eqref{ode1} is solved by any radially symmetric solution to \eqref{equ1}. The desired Pohozaev identity is indeed given by \eqref{drj}. Next, we need to discuss the sign of the associated Pohozaev quantity $J(r; u)$ defined by \eqref{defj} and also to discover some properties of solutions to \eqref{ode1}. Here, the condition \eqref{uniquec}, which is technical, is actually to ensure that $J(r; u) \geq 0$ for any $r \geq 0$, see Lemma \ref{unique11}. The presence of such a restriction is due to the fact that our problem has competing inhomogeneous nonlinearities, under which the Pohozaev identity \eqref{drj} has two terms with indefinite signs, and thus the verification of the assertion that $J(r; u) \geq 0$ for any $r \geq 0$ becomes difficult. It is not clear whether Theorem~\ref{thmunique} continues to hold when condition~\eqref{uniquec} is weakened.

\begin{rem}
{\rm Arguing as in the proof of \eqref{drj} below and replacing the role of $p_2$ by $p_1$, one can show that
\begin{align*}
\frac{d}{dr} J(r;\phi)=\widetilde{G}(r)\phi(r)^2+\widetilde{H}(r)\phi(r)^{p_1},
\end{align*}
where
\begin{align*}
\widetilde{G}(r)&:=\omega\left(\frac{(N-1)(2-p_2)-2b_2}{p_2+2}\right)r^{\frac{2(b_2-1)+(2N-3)p_2}{p_2+2}}\\
& \quad +\frac{\left(2(N-1)-b_2\right)\left(2(N-b_2)-p_2(N-2)\right)\left((b_2-2)+(N-2)p_2\right)}{(p_2+2)^3}r^{\frac{2(b_2-3)+(2N-5)p_2}{p_2+2}},
\end{align*}
and
\begin{align*}
\widetilde{H}(r)&:=\left(\frac{2(N-1)-b_2}{p_2+2}-\frac{2b_2+2(N-1)p_2}{p_1(p_2+2)}+\frac{b_1}{p_1}\right)r^{\frac{2(b_2-1)+(2N-3)p_2}{p_2+2}-b_1}.
\end{align*}
Assume in addition that
\begin{align*}
\frac{2(N-1)-b_2}{p_2+2}-\frac{2b_2+2(N-1)p_2}{p_1(p_2+2)}+\frac{b_1}{p_1} \geq 0,
\end{align*}
we see that the conclusion of Theorem \ref{thmunique} remains true.}
\end{rem}

\begin{defi}
We say that a solution $\phi \in H^1(\R^N)$ to \eqref{equ1} is non-degenerate if
{
$$
Ker[\mathcal{L}_+]=0,
$$
where
\begin{equation}
    \label{L+}
    \mathcal{L}_+:=-\Delta +\omega+(p_1-1)|x|^{-b_1}|\phi|^{p_1-2}-(p_2-1)|x|^{-b_2}|\phi|^{p_2-2}.
\end{equation}
}
Moreover, we shall use $n(\mathcal{L}_+)$ to denote the Morse index of $\phi$, which is defined by the maximal dimension of a subspace of $H^1(\R^N)$ on which $\mathcal{L}_+$ is negative definite.
\end{defi}

\begin{thm} \label{nond}
Let \eqref{H} hold, $N \geq 3$ and $\phi \in H^1(\R^N)$ be a positive radially symmetric solution to \eqref{equ1} that satisfies the Morse index of $\phi$ is 1, i.e., $n(\mathcal{L}_+)=1$. Then $\phi$ is non-degenerate.
\end{thm}

{ 
\begin{rem} It should be pointed out that if \eqref{H} holds and $p_1 \leq p_2$, then the ground states derived in Theorem \ref{thm1} do satisfy the assumption of Theorem \ref{nond}, i.e. the Morse index of any ground state to \eqref{equ1} is one. Indeed, let $\phi \in P$, then $K(\phi)=0$ and
\begin{align*}
\langle K'(\phi), \phi \rangle &=2\int_{\R^N} |\nabla \phi|^2 \,dx+\frac{N(p_1-2)+2b_1}{2}\int_{\R^N}|x|^{-b_1}|\phi|^{p_1}\, dx\\
& \quad -\frac{N(p_2-2)+2b_2}{2}\int_{\R^N}|x|^{-b_2}|\phi|^{p_2}\, dx \\
&=(2-p_2)\int_{\R^N} |\nabla \phi|^2 \,dx+\frac{N(p_1-2)+2b_1}{2} \left(1-\frac{p_2}{p_1}\right)\int_{\R^N}|x|^{-b_1}|\phi|^{p_1}\, dx.
\end{align*}
Since $2<p_1 \leq p_2$, then $\langle K'(\phi), \phi \rangle <0$. It then follows from the implicit function theorem that $P$ is a $C^1$-manifold of codimension 1 and $H^1(\R^N)=T_{\phi} P \oplus \R \phi$ for any $\phi \in P$, where $T_{\phi} P$ denotes the tangent space of $P$ at $\phi$. Therefore, we know that the Morse index of any ground state to \eqref{equ1} is at most one, because any ground state to \eqref{equ1} corresponds to a minimizer to $S_{\omega}$ restricted on $P$. On the other hand, if $\phi \in H^1(\R^N)$ is a ground state to \eqref{equ1}, then
$$
\int_{\R^N}|x|^{-b_1}|\phi|^{p_1}<\int_{\R^N} |\nabla \phi|^2 \,dx +\omega \int_{\R^N} |\phi|^2 \,dx+ \int_{\R^N}|x|^{-b_1}|\phi|^{p_1} \,dx=\int_{\R^N}|x|^{-b_2}|\phi|^{p_2} \,dx.
$$
Therefore, we obtain that
\begin{align*}
\langle \mathcal{L}_+[\phi], \phi \rangle&= \int_{\R^N} |\nabla \phi|^2 \,dx +\omega \int_{\R^N} |\phi|^2 \,dx+ (p_1-1) \int_{\R^N}|x|^{-b_1}|\phi|^{p_1} \,dx -(p_2-1)\int_{\R^N}|x|^{-b_2}|\phi|^{p_2} \,dx \\
&=\left(p_1-2\right) \int_{\R^N}|x|^{-b_1}|\phi|^{p_1} \,dx - \left(p_2-2\right) \int_{\R^N}|x|^{-b_2}|\phi|^{p_2} \,dx \\
& < -\left(p_2-p_1\right) \int_{\R^N}|x|^{-b_2}|\phi|^{p_2} \,dx \\
& \leq 0.
\end{align*}
This in turn infers that the Morse index of any ground state to \eqref{equ1} is at least one. Consequently, the desired conclusion follows.
\end{rem}
}
The discussion of non-degeneracy of solutions to nonlinear elliptic equations, which relies principally on the spherical harmonic decomposition arguments, plays an important role in revealing quantitative properties of solutions. We refer the readers to the early works \cite{Oh} and \cite{We1} for the study of non-degeneracy of solutions to \eqref{sequ}. For the study of non-degeneracy of solutions to \eqref{kequ} with $K(x)=|x|^{-b}$ for $b>0$, we refer the readers to \cite{GS}. To establish Theorem \ref{nond}, we shall make use of the spherical harmonic decomposition arguments, which is inspired by \cite{St}. Note that our problem has competing inhomogeneous nonlinearities, then the proof of the result is not a direct application of these arguments.

In what follows, we shall turn to discuss dynamical behaviors of solutions to the Cauchy problem \eqref{equ}. Let us first state the well-posedness result for  \eqref{equ} in $H^1(\R^N)$.
\begin{prop}
Let \eqref{H} holds. Then, for any $u_0 \in H^1(\R^N)$, there exists  $T_{max}>0$ such that \eqref{equ} has a unique maximal solution $u \in C([0, T_{max}), H^1(\R^N))$ to \eqref{equ} satisfying the conservation of mass and energy, i.e. for any $ t \in [0, T_{max})$,
$$
M(u(t))=M(u_0), \quad E(u(t))=E(u_0),
$$
where
$$
M(u(t)):=\int_{\R^N} |u(t,x)|^2 \,dx,
$$
$$
E(u):=\frac 12 \int_{\R^N} |\nabla u(t,x)|^2 \,dx + \frac{1}{p_1} \int_{\R^N} |x|^{-b_1} |u(t,x)|^{p_1} \,dx - \frac{1}{p_2} \int_{\R^N} |x|^{-b_2} |u(t,x)|^{p_2} \,dx.
$$
Moreover, the solution map $u_0 \mapsto u$ is continuous from $H^1(\R^N)$ to $C([0, T_{max}), H^1(\R^N))$. There also holds that either $T_{max}< +\infty$ or $\displaystyle\lim_{t \to T^-_{max}} \|\nabla u(t)\|_2=+\infty.$
\end{prop}
Although we deal here with competing inhomogeneous nonlinearities, the proof of the above proposition mimics the same steps performed in \cite{Gu} where only one inhomogeneous nonlinearity is considered. The details are hence omitted.

The first result in this direction consists in blowup of solutions to the Cauchy problem \eqref{equ} with initial data belonging to the set $\mathcal{A}^-_{\omega}$ defined by
\begin{align} \label{a-}
\mathcal{A}^-_{\omega}:=\Big\{\phi \in H^1(\R^N) \backslash \{0\} : S_{\omega}(\phi)<m_{\omega}, K(\phi) <0\Big\},
\end{align}

where $m_{\omega}>0$ is the ground state energy level defined by
\begin{align} \label{smin}
m_{\omega}:=\inf_{\phi \in P} S_{\omega}(\phi),
\end{align}
with
$$
P:= \Big\{\phi \in H^1(\R^N) \backslash \{0\} : K(\phi) =0\Big\}\,\, \mbox{for} \,\,\omega>0,  \quad P:= \Big\{\phi \in X \backslash \{0\} : K(\phi) =0\Big\}\,\, \mbox{for} \,\,\omega=0.
$$
Note that by Theorem \ref{thm1} the infimum $m_{\omega}$ is positive and is actually a minimum  for any $\omega \geq 0$. Here $P$ is the so-called Pohozaev manifold and $K(\phi)=0$ is the Pohozaev identity related to \eqref{equ1}, see Lemma \ref{ph} below.

\begin{thm} \label{thm2}
Let \eqref{H} holds and $u \in C([0, T_{max}), H^1(\R^N))$ be the maximal solution to the Cauchy problem \eqref{equ} with $u_0 \in \mathcal{A}^-_{\omega}$. Then $u$ blows up in finite time under one of the following assumptions,
\begin{itemize}
\item [$(\textnormal{i})$] $|x| u_0 \in L^2(\R^N)$;
\item [$(\textnormal{ii})$] $u_0 \in H^1_{rad}(\R^N)$  and ${\max\{p_1,p_2\}\leq6}$;
\item [$(\textnormal{iii})$] $u_0 \in H^1(\R^N)$ and ${\max\{p_1,p_2\}\leq 2+\frac{4}{N}}$.
\end{itemize}
\end{thm}

\begin{rem}
{\rm The blow-up in finite time for NLS is standard for a finite variance data or a radial one. However, for a non-radial data with possible infinite variance, the blow-up in finite time is still open except for the one space dimension \cite{ogawa}. Here we obtain the finite time blow-up of solutions to \eqref{equ} for a non-radial data with possible infinite variance for the range $2+\frac{2(2-b_j)}{N}<p_j\leq2+\frac4N,\,j=1,2$. Note that the above range  is unmeaningful if $b_1\,b_2=0$.  }
\end{rem}
\begin{rem}
    {\rm In a recent paper \cite{BFG-23}, the authors establish the existence of finite-time blowing-up solutions below the ground state energy threshold in the 3D homogeneous case, that is $b_1=b_2=0$.}
\end{rem}

To prove Theorem \ref{thm2}, we need to employ the variational characterization of the ground state energy level given by \eqref{smin} along with the evolution of the related localized virial quantity, see Lemma \ref{vbe}. Further, utilizing the arguments in \cite{MRS}, we are able to derive the upper bound on blowup rate of solutions.

\begin{thm} \label{thm3}
Let \eqref{H} holds and $u \in C([0, T_{max}), H^1(\R^N))$ be the maximal solution of the Cauchy problem \eqref{equ} with $u_0 \in \mathcal{A}^-_{\omega}$. Assume that $u$ blows up in finite time, i.e. $T_{max}<+\infty$.
\begin{itemize}
\item [$(\textnormal{i})$] If $u_0 \in H^1(\R^N)$ is radial, $p_1<6$ and $p_2 < 6$, then, for any $t>0$ close to $T_{max}$,
\begin{align} \label{r1}
\int_t^{T_{max}} (T_{max}-\tau) \|\nabla u(\tau)\|_2^2 \,d\tau \leq C (T_{max}-t)^{\min \left\{\frac{2(6-p_1)}{(N-2)(p_1-2)+2(b_1+2)}, \frac{2(6-p_2)}{(N-2)(p_2-2)+2(b_2+2)}\right\}}.
\end{align}
Moreover, there exists a sequence $\{t_n\} \subset [0, T_{max})$ with $t_n \nearrow T_{max}$ as $n \to \infty$ such that
\begin{align} \label{r2}
\|\nabla u(t_n)\|_2 \leq \frac{C}{(T_{\max}-t_n)^{\max \left\{\frac{(N-1)(p_1-2)+2b_1}{(N-2)(p_1-2)+2(b_1+2)}, \frac{(N-1)(p_2-2)+2b_2}{(N-2)(p_2-2)+2(b_2+2)}\right\}}}.
\end{align}
\item [$(\textnormal{ii})$] If $u_0 \in H^1(\R^N)$, $p_1<2+\frac{4}{N}$, $p_2 < 2+\frac{4}{N}$, then, for any $t>0$ close to $T_{max}$,
\begin{align} \label{r11}
\int_t^{T_{max}} (T_{max}-\tau) \|\nabla u(\tau)\|_2^2 \,d\tau \leq C (T_{max}-t)^{\min\left\{\frac{2(4-N(p_1-2))}{4-N(p_1-2)+2b_1}, \frac{2(4-N(p_2-2))}{4-N(p_2-2)+2b_2}\right\}}.
\end{align}
Moreover, there exists a sequence $\{t_n\} \subset [0, T_{max})$ with $t_n \nearrow T_{max}$ as $n \to \infty$ such that
\begin{align} \label{r22}
\|\nabla u(t_n)\|_2 \leq \frac{C}{(T_{\max}-t_n)^{\max \left\{\frac{2b_1}{4-N(p_1-2)+2b_1}, \frac {2b_2}{4-N(p_2-2)+2b_2}\right\}}}.
\end{align}
\end{itemize}
\end{thm}

As a straightforward consequence of Theorem \ref{thm2}, we  have the strong instability of standing waves.

\begin{cor} \label{instability}
Let \eqref{H} holds, $\omega \geq 0$, and $e^{\textnormal{i} \omega t} \phi$  be a standing wave solution to \eqref{equ1-1}. Then the following assertions hold.
\begin{itemize}
\item [$(\textnormal{i})$] If $\omega >0$ or $\omega=0$ and $N \geq 5$, then standing wave $e^{\textnormal{i} \omega t} \phi$ is strongly unstable in the sense that, for any $\eps>0$, there exists $u_0 \in H^1(\R^N)$ such that $\|u_0-\phi\| <\eps$ and the maximal solution $u$ of \eqref{equ}  blows up in finite time.
\item[$(\textnormal{ii})$] If $\omega=0$ and $3 \leq N \leq 4$, then the standing wave $e^{\textnormal{i} \omega t} \phi$ is strongly unstable in the sense that, for any $\eps>0$, there exists $u_0 \in H^1(\R^N)$ such that $\|u_0-\phi\|_{X} <\eps$ and the maximal solution $u$ of \eqref{equ}  blows up in finite time.
\end{itemize}
\end{cor}

Since $\phi \not \in L^2(\R^N)$ for $\omega=0$ and $3 \leq N \leq 4$ by Theorem \ref{decay}, then the strong instability in the statement $(\textnormal{ii})$ holds true in the framework of the weak topology.

Next we shall investigate scattering of solutions to the Cauchy problem \eqref{equ} with initial data belonging to the set $\mathcal{A}^+_{\omega}$ defined by
$$
\mathcal{A}^+_{\omega}:=\left\{\phi \in H^1(\R^N) \backslash \{0\} : S_{\omega}(\phi)<m_{\omega}, K(\phi)>0\right\}.
$$
Our next result is the following energy scattering for \eqref{equ} for initial data belonging to $\mathcal{A}^+_{\omega}$.

\begin{thm} \label{scattering}
Let \eqref{H} holds, $N \geq 3$,   $0<b_1, b_2<\min\{2, {N}/{2}\}$ and $u_0 \in \mathcal{A}^+_{\omega}$. Then the solution $u$ to the Cauchy problem  \eqref{equ}  exists globally in time  and scatters  in the sense that there exist $u_0^{\pm} \in H^1(\R^N)$ such that
$$
\lim_{t \to \pm\infty} \left\|u(t)-e^{\textnormal{i} t \Delta} u_0^{\pm}\right\|_{H^1} =0.
$$
\end{thm} 
\begin{rem}
{\rm In view of the  scattering result stated in the above theorem, some comments are in order.
~\begin{itemize}
  \item[(i)] The 3D homogeneous case ($b_1=b_2=0$) is studied in \cite{BDF}. The main ingredient used in \cite{BDF} was the interaction Morawetz estimate which fails to hold in the inhomogeneous case.
  \item[(ii)] Our approach relies on Tao's scattering criterion \cite{Ta1}, and Dodson-Murphy's Virial/Morawetz inequalities \cite{DM1} (see also \cite{ADM, DK}).
  \item[(iii)] We are face here to two difficulties coming from the singular weights: the scattering criterion and the coercivity.
  \item[(iv)] Unlike to the homogeneous case, where the Strauss inequality is used to handle the radial setting, here the decay of the singular weights enables us to remove the radial assumption.
    \end{itemize}
    }
\end{rem}
To prove the global existence of solutions, one can take into account the variational characterization of the ground state energy level given by \eqref{smin} and the conservation laws. While, to prove scattering, we shall adapt the arguments from \cite{Ta1}, which avoids the use of the concentration-compactness-rigidity techniques due to Kenig and Merle in \cite{KM}. First we need to establish coercivity result, see Lemma \ref{coerc}. Next we need to prove that solutions have sufficiently small mass near the origin at sufficiently late time, see Lemma \ref{crt}, which is achieved by small data scattering theory, see Lemmas \ref{crt1} and \ref{smlsct}. Thus, we have the desired conclusions.
\medskip

{\noindent {\bf\large Outline of the paper.} The rest of this paper is organized as follows. In Section \ref{pre}, we present some auxiliary results and useful tools used to establish our main Theorems. In Section \ref{existnon}, we discuss the existence / nonexistence, symmetry, and decay of solutions to \eqref{equ1} and give the proofs of Theorem \ref{thm1} and Theorem \ref{decay}. Section \ref{unno} is devoted to the uniqueness and non-degeneracy of solutions and contains the proofs of Theorem \ref{thmunique} and Theorem \ref{nond}. Section \ref{blowupr} is concerned with the study of blow-up of solutions and contains the proofs of Theorem \ref{thm2} and Theorem \ref{thm3} as well as Corollary \ref{instability}. Finally, in Section \ref{sca}, we investigate the energy scattering for \eqref{equ} and establish Theorem \ref{scattering}.}
\medskip

{\noindent {\bf\large Notations.} For $1 \leq p \leq \infty$, the norm in the usual  Lebesgue space $L^p(\R^N)$ will be denoted simply $\|\cdot\|_p$. For two quantities $A$ and $B$, $A \lesssim B$ means that there exists  a positive constant $C$ such that $A \leq C B$, and $A\gtrsim B$ means that $A\geq C B$ . We use $p'$ to denote the conjugate exponent of $p$ defined by $p'=\frac{p}{p-1}$. Along the rest of the article, the letter $C$ stands for a generic positive constant, whose value may change from line to line. For simplicity and clarity in the presentation, we will always denote by $u$ a solution of \eqref{equ1-1} or \eqref{equ1}.}
\medskip

\section{Preliminaries} \label{pre}

In this section, we collect some auxiliary results and useful tools needed in the proofs of our main results. Let us begin with the well-known Gagliardo-Nirenberg inequality \cite{G} and radial Sobolev embedding  \cite{Str}.

\begin{lem}
Let $N \geq 1$, $0 \leq b<\min\{2, N\}$ and $2<p\leq 2^*_b$. Then there exists $C>0$ such that, for any $u \in H^1(\R^N)$,
\begin{align} \label{GN}
\int_{\R^N} |x|^{-b} |u|^{p} \, dx \leq C\left(\int_{\R^N} |\nabla u|^2 \, dx\right)^{\frac{N(p-2)}{4} + \frac b 2} \left(\int_{\R^N} |u|^2 \,dx\right)^{\frac p2 -\frac{N(p-2)}{4} - \frac b 2}.
\end{align}
\end{lem}

\begin{lem} \label{hrc}
Let $N \geq 2$. Then there exists $C>0$ such that, for any $u \in H^1_{rad}(\R^N)$,
\begin{align} \label{re}
\sup_{x \neq 0} |x|^{\frac{N-1}{2}}|u(x)| \leq C \left(\int_{\R^N}|\nabla u|^2 \,dx \right)^{\frac 14}\left(\int_{\R^N}|u|^2 \,dx \right)^{\frac 14}.
\end{align}
\end{lem}

\begin{lem} \cite[Lemma 2.1]{AC} \label{cembedding}
Let $N \geq 1$ and $0<b<\min\{2, N\}$. Then $H^1(\R^N)$ is compactly embedded into $L^p(\R^N, |x|^{-b}dx)$ for any $2<p<2^*_b$.
\end{lem}

\begin{lem} \label{ph}
Let $u \in H^1(\R^N)$ be a solution to \eqref{equ1}. Then $K(u)=0$, where $K$ is as in \eqref{K-phi}.
\end{lem}
\begin{proof}
First multiplying \eqref{equ1} by $x \cdot \nabla u$ and integrating on $B_R(0)$ results in
\begin{align} \label{ph11}
\begin{split}
&-\int_{B_R(0)} \Delta u (x \cdot \nabla u ) \, dx + \omega \int_{B_R(0)} u (x \cdot \nabla u ) \, dx \\
&=\int_{B_R(0)}|x|^{-b_2}|u|^{p_2-2} u (x \cdot \nabla u) \,dx -\int_{B_R(0)}|x|^{-b_1}|u|^{p_1-2} u (x \cdot \nabla u) \,dx.
\end{split}
\end{align}
In what follows, we are going to calculate every term in \eqref{ph11} with the help of the divergence theorem. Observe that
\begin{align*}
-\int_{B_R(0)} \Delta u (x \cdot \nabla u ) \, dx&= \int_{B_R(0)} \nabla u \cdot \nabla (x \cdot \nabla u) \, dx-R\int_{\partial B_R(0)} \left| \nabla u \cdot \bf{n} \right|^2 \,dS\\
&=\frac{2-N}{2} \int_{B_R(0)}|\nabla u|^2\, dx-\frac R 2 \int_{\partial B_R(0)} \left| \nabla u \cdot \bf{n} \right|^2 \,dS,
\end{align*}
\begin{align*}
\int_{B_R(0)} u (x \cdot \nabla u ) \, dx=\frac 12 \int_{B_R(0)} x \cdot \nabla \left(|u|^2\right) \, dx=-\frac N 2 \int_{B_R(0)} |u|^2 \, dx + \frac  R 2 \int_{\partial B_R(0)} |u|^2 \, dS,
\end{align*}
where the vector $\bf{n}$ denotes the outward normal to $\partial B_R(0)$. In addition, we see that
\begin{align*}
\int_{B_R(0)}|x|^{-b_1}|u|^{p_1-2} u (x \cdot \nabla u) \,dx&=\frac {1}{p_1} \int_{B_R(0)} \left(|x|^{-b_1}x\right) \cdot \nabla\left(|u|^{p_1}\right) \, dx\\
&=\frac {b_1-N}{p_1} \int_{B_R(0)} |x|^{-b_1}|u|^{b_1} \, dx-\frac  {R} {p_1} \int_{\partial B_R(0)} |x|^{-b_1} |u|^{p_1} \, dS,
\end{align*}
\begin{align*}
\int_{B_R(0)}|x|^{-b_2}|u|^{p_2-2} u (x \cdot \nabla u) \,dx=\frac {b_2-N}{p_2} \int_{B_R(0)} |x|^{-b_2}|u|^{p_2} \, dx-\frac  {R} {p_2} \int_{\partial B_R(0)} |x|^{-b_2} |u|^{p_2} \, dS.
\end{align*}
Therefore, applying \eqref{ph11}, we are able to derive that
\begin{align} \label{ph1}
\begin{split}
&\frac{N-2}{2} \int_{B_R(0)}|\nabla u|^2\, dx +\frac {\omega N}{2} \int_{B_R(0)} |u|^2 \, dx \\
&=\frac {b_2-N}{p_2} \int_{B_R(0)} |x|^{-b_2}|u|^{p_2} \, dx-\frac {b_1-N}{p_1} \int_{B_R(0)} |x|^{-b_1}|u|^{p_1} \, dx +I_R,
\end{split}
\end{align}
where
$$
I_R:=-\frac R 2 \int_{\partial B_R(0)} \left| \nabla u \cdot \bf{n} \right|^2 \,dS+\frac  R 2 \int_{\partial B_R(0)} |u|^2 \, dS +\frac {R}{p_1} \int_{\partial B_R(0)} |x|^{-b_1} |u|^{p_1}\, dS-\frac {R}{p_2} \int_{\partial B_R(0)} |x|^{-b_2} |u|^{p_2}\, dS.
$$
Since $u \in H^1(\R^N)$, then
\begin{align*}
&\int_{\R^N} |\nabla u|^2 + |u|^2 + |x|^{-b_1} |u|^{p_1} + |x|^{-b_2} |u|^{p_2} \, dx \\
&= \int_0^{\infty} \int_{\partial B_R(0)}|\nabla u|^2 + |u|^2 + |x|^{-b_1} |u|^{p_1} + |x|^{-b_2} |u|^{p_2} \,dSdR<\infty.
\end{align*}
It then follows that there exists a sequence $\{R_n\} \to \infty$ as $n \to \infty$ such that $I_{R_n}=o_n(1)$.  Making use of \eqref{ph1} with $R=R_n$ and taking the limit as $n \to \infty$, we then get that
\begin{align*} 
\frac{N-2}{2} \int_{\R^N}|\nabla u|^2\, dx +\frac {\omega N}{2} \int_{\R^N} |u|^2 \, dx=\frac {b_2-N}{p_2} \int_{\R^N} |x|^{-b_2}|u|^{p_2} \, dx-\frac {b_1-N}{p_1} \int_{\R^N} |x|^{-b_1}|u|^{p_1} \, dx.
\end{align*}
On the other hand, multiplying \eqref{equ1} by $u$ and integrating on $\R^N$ leads to
\begin{align*}
\frac{N}{2} \int_{\R^N}|\nabla u|^2\, dx +\frac {\omega N}{2} \int_{\R^N} |u|^2 \, dx=\frac {N}{2} \int_{\R^N} |x|^{-b_2}|u|^{p_2} \, dx+\frac N 2 \int_{\R^N} |x|^{-b_1}|u|^{p_1} \, dx.
\end{align*}
Thereby, we conclude that
$$
\int_{\R^N} |\nabla u|^2 \,dx=\frac{N(p_2-2)+2b_2}{2p_2}\int_{\R^N}|x|^{-b_2}|u|^{p_2}\, dx-\frac{N(p_1-1)+2b_1}{2p_1}\int_{\R^N}|x|^{-b_1}|u|^{p_1}\, dx.
$$
Thus the proof is completed.
\end{proof}

\section{Existence/Nonexistence, Symmetry and Decay} \label{existnon}

In this section, we shall prove Theorem \ref{thm1}. The first purpose is to verify the existence of ground states to \eqref{equ1} for $\omega>0$. For this, we shall introduce the following minimization problem,
\begin{align} \label{ming1}
m_{\omega}:=\inf_{u \in P} S_{\omega}(u),
\end{align}
where
$$
P:=\left\{ u \in H^1(\R^N) \backslash\{0\} : K(u)=0\right\}.
$$

\begin{lem} \label{maximum}
Let \eqref{H} holds. Then, for any $u \in H^1(\R^N)$, there exists a unique $t_u>0$ such that $u_{t_u} \in P(c)$ and
\begin{align} \label{max1}
\max_{t>0} E(u_t)=E(u_{t_u}).
\end{align}
Moreover, if $K(u) < 0$, then $0<t_u < 1$. In adiition,the function $t \mapsto E(u_t)$ is concave on $[t_u, \infty)$.
\end{lem}
\begin{proof}
Observe that
\begin{align*}
\frac{d}{dt} E(u_t)&=t\int_{\R^N} |\nabla u|^2 \,dx+\frac{N(p_1-2)+2b_1}{2p_1}t^{\frac N 2 (p_1-2)+b_1-1}\int_{\R^N}|x|^{-b_1}|u|^{p_1} \, dx\\
&\quad-\frac{N(p_2-2)+2b_2}{2p_2}t^{\frac N 2 (p_2-2)+b_2-1}\int_{\R^N}|x|^{-b_2}|u|^{p_2} \, dx \\
&=\frac 1 t K(u_t).
\end{align*}
In addition, we have that
\begin{align*}
\frac{d^2}{dt^2} E(u_t)&=\int_{\R^N} |\nabla u|^2 \,dx+\frac{(N(p_1-2)+2b)((N(p_1-2)-2(1-b_1))}{4p_1}t^{\frac{N}{2}(p_1-2)+b_1-2}\int_{\R^N}|x|^{-b_1}|u|^{p_1} \, dx \\
& \quad -\frac{(N(p_2-2)+2b_2)((N(p_2-2)-2(1-b_2))}{4p_2}t^{\frac{N}{2}(p_2-2)+b_2-2}\int_{\R^N}|x|^{-b_2}|u|^{p_2} \, dx.
\end{align*}
Note that $p_1>2+\frac{2(2-b_1)}{N}$ and $p_2>2+\frac{2(2-b_2)}{N}$. Then it is not hard to see that there exists a unique $t_u>0$ such that $K(u_{t_u})=0$. In addition, there holds that $\frac{d}{dt} E(u_t)>0$ if $0<t<t_u$ and $\frac{d}{dt} E(u_t)<0$ if $t>t_u$. This then gives rise to \eqref{max1}. Meanwhile, we can see that $0<t_u<1$ if $K(u)<0$. Furthermore, it is simple to check that there exists a unique $0<\tilde{t}_{u} \leq t_u$ such that
$$
\frac{d^2}{dt^2} E(u_t)\mid_{t=\tilde{t}_u}=0.
$$
There also holds that
$$
\frac{d^2}{dt^2} E(u_t)>0 \,\, \mbox{for} \,\, 0<t<\tilde{t}_u, \quad \frac{d^2}{dt^2} E(u_t)<0 \,\, \mbox{for} \,\, t>\tilde{t}_u.
$$
This implies that the function $t \mapsto E(u_t)$ is concave on $[t_u, \infty)$ and the proof is completed.
\end{proof}

\begin{lem}\label{existence1}
Let \eqref{H} holds. Then there exist ground states to \eqref{equ1} for $\omega>0$.
\end{lem}
\begin{proof}
First it is standard to check that $P$ is a natural constraint. Hence, any minimizer to \eqref{ming1} is a ground state to \eqref{equ1}. Therefore, we are going to prove that there exist minimizers to \eqref{ming1}. Let $\{u_n\} \subset P$ be a minimizing sequence to \eqref{ming1}, i.e. $S_{\omega}(u_n)+o_n(1)=m_{\omega}$. Since $P$ is a natural constraint, without restriction, then we may assume that $S_{\omega}'(u_n)=o_n(1)$. Note that $K(u_n)=0$, by \eqref{GN}, then
\begin{align*}
\int_{\R^N} |\nabla u_n|^2 \,dx+\frac{N(p_1-1)+2b_1}{2p_1} \int_{\R^N}|x|^{-b_1}|u_n|^{p_1}\, dx
&=\frac{N(p_2-2)+2b_2}{2p_2}\int_{\R^N}|x|^{-b_2}|u_n|^{p_2}\, dx \\
& \leq C \|\nabla u_n\|_2^{\frac{N(p_2-2)}{2}+b_2} \|u_n\|_2^{p_2-\frac{N(p_2-2)}{2}-b_2}.
\end{align*}
This obviously leads to
\begin{align} \label{bddbe}
\|\nabla u_n\|_2^{\frac{N(p_2-2)}{2}+b_2-2} \|u_n\|_2^{p_2-\frac{N(p_2-2)}{2}-b_2} \geq \frac{1}{C}.
\end{align}
On the other hand, we observe that
\begin{align} \label{bbe1} \nonumber
m_{\omega}=S_{\omega}(u_n)+o_n(1)&=S_{\omega}(u_n)-\frac{2}{N(p_2-2)+2b_2}K(u_n)+o_n(1) \\
&=\left(\frac 12-\frac{2}{N(p_2-2)+2b_2}\right)\|\nabla u_n\|_2^2 +\frac 12 \|u_n\|_2^2 \\ \nonumber
& \quad +\frac{1}{p_1} \left(\frac{N(p_2-p_1)+2(b_2-b_1)}{N(p_2-2)+2b_2}\right)\int_{\R^N} |x|^{-b_1}|u_n|^{p_1}\,dx+o_n(1).
\end{align}
In view of \eqref{bddbe}, we then have that $m_{\omega}>0$. Next we assert that $m_{\omega}$ is achieved. It follows immediately from \eqref{bbe1} that $\{u_n\}$ is bounded in $H^1(\R^N)$. Then there exists a nontrivial $u \in H^1(\R^N)$ such that $u_n \wto u$ in $H^1(\R^N)$ as $n \to \infty$ and $u_n \to u$ in $L^{p_1}(\R^N, |x|^{-b_1}dx) \cap L^{p_2}(\R^N, |x|^{-b_2}dx)$ as $n \to \infty$ by Lemma \ref{cembedding}. This infers that $S_{\omega}(u) \leq m_{\omega}$. Furthermore, we know that $S_{\omega}'(u)=0$, i.e $u \in H^1(\R^N)$ is a solution to \eqref{equ1}. As a consequence, by Lemma \ref{ph}, we get that $K(u)=0$. It then follows that $S_{\omega}(u) = m_{\omega}$. 
Thus the proof is completed.
\end{proof}

In the following, we are going to prove the existence of ground states to \eqref{equ1} for $\omega=0$, which induces the following zero mass equation,
\begin{align} \label{equ12}
-\Delta u+|x|^{-b_1} |u|^{p_1-2}u=|x|^{-b_2}|u|^{p_2-2}u \quad \mbox{in} \,\, \R^N.
\end{align}
To investigate the existence of solutions to \eqref{equ12}, we need to work in the Sobolev space $X$ defined by the completion of $C_0^{\infty}(\R^N)$ under the norm
$$
\|u\|_X:=\left(\int_{\R^N} |\nabla u|^2 \,dx \right)^{\frac 12} + \left(\int_{\R^N} |x|^{-b_1}|u|^{p_1}\,dx\right)^{\frac{1}{p_1}}.
$$
\begin{lem} \label{em}
Let \eqref{H} holds, $N \geq 3$, $p_1<p_2$ and $b_1<b_2$. Then $X$ is continuously embedding into $L^{p_2}(\R^N, |x|^{-b_2}dx)$, i.e. there exists $C>0$ such that, for any $u \in X$,
\begin{align} \label{ineq1}
\int_{\R^N} |x|^{-b_2}|u|^{p_2} \,dx \leq C \left(\int_{\R^N} |x|^{-b_1}|u|^{p_1}\,dx\right)^{\theta_1 \theta_2}\left(\int_{\R^N} |\nabla u|^2 \,dx \right)^{\frac{\theta_1(1-\theta_2)p_1}{2}+\frac{(1-\theta_1)(N-b_2)}{N-2}},
\end{align}
where $0<\theta_1, \theta_2<1$ are two constants such that
$$
p_2=\theta_1 p_1+(1-\theta_1)\frac{2(N-b_2)}{N-2}, \quad b_2=\theta_2b_1+(1-\theta_2)\frac{2N-p_1(N-2)}{2}.
$$
\end{lem}
\begin{proof}
Since $p_1<p_2<\frac{2(N-b_2)}{N-2}$, then there exists $0<\theta_1<1$ such that
$$
p_2=\theta_1 p_1+(1-\theta_1)\frac{2(N-b_2)}{N-2}.
$$
Using H\"older's inequality and \eqref{GN}, we then have that
\begin{align}\label{c1111}
\begin{split}
\int_{\R^N} |x|^{-b_2}|u|^{p_2} \,dx &\leq \left(\int_{\R^N} |x|^{-b_2}|u|^{p_1}\,dx\right)^{\theta_1} \left(\int_{\R^N} |x|^{-b_2}|u|^{\frac{2(N-b_2)}{N-2}}\,dx\right)^{1-\theta_1} \\
& \leq C \left(\int_{\R^N} |x|^{-b_2}|u|^{p_1}\,dx\right)^{\theta_1} \left(\int_{\R^N} |\nabla u|^2\,dx\right)^{\frac{(1-\theta_1)(N-b_2)}{N-2}}.
\end{split}
\end{align}
Define
$$
b_0:=\frac{2N-p_1(N-2)}{2}.
$$
It is simple to see that
$$
b_0<2, \quad b_0<N, \quad \frac{2(N-b_0)}{N-2}=p_1.
$$
Since $p_1<p_2<\frac{2(N-b_2)}{N-2}$, then $b_0>b_2$. Thus there exists $0<\theta_2<1$ such that $b_2=\theta_2b_1+(1-\theta_2)b_0$. Using again H\"older's inequality and \eqref{GN}, we then get that
\begin{align}\label{c222}
\begin{split}
\int_{\R^N} |x|^{-b_2}|u|^{p_1}\,dx &\leq \left(\int_{\R^N} |x|^{-b_1}|u|^{p_1}\,dx\right)^{\theta_2} \left(\int_{\R^N} |x|^{-b_0}|u|^{p_1}\,dx\right)^{1-\theta_2} \\
& \leq C \left(\int_{\R^N} |x|^{-b_1}|u|^{p_1}\,dx\right)^{\theta_2} \left(\int_{\R^N} |\nabla u|^2 \,dx \right)^{\frac{(1-\theta_2)p_1}{2}}.
\end{split}
\end{align}
Consequently, by \eqref{c1111} and \eqref{c222}, there holds that
$$
\int_{\R^N} |x|^{-b_2}|u|^{p_1}\,dx \leq C \left(\int_{\R^N} |x|^{-b_2}|u|^{p_1}\,dx\right)^{\theta_1 \theta_2}\left(\int_{\R^N} |\nabla u|^2 \,dx \right)^{\frac{\theta_1(1-\theta_2)p_1}{2}+\frac{(1-\theta_1)(N-b_2)}{N-2}}.
$$
This completes the proof.
\end{proof}

\begin{lem} \label{xem}
Let \eqref{H} holds, $N \geq 3$, $p_1<p_2$ and $b_1<b_2$. Then $X$ is compactly embedding into $L^{p_2}(\R^N, |x|^{-b_2}dx)$.
\end{lem}
\begin{proof}
Observe that
\begin{align} \label{c111'}
\begin{split}
\int_{\R^N} |x|^{-b_2} |u|^{p_2} \,dx &\leq \int_{|x| < R} |x|^{-b_2} |u|^{p_2} \,dx +\frac{1}{R^{b_2}}\int_{|x| \geq R} |u|^{p_2}\,dx \\
& \leq \left(\int_{|x|<R} |u|^{-\frac{\theta b_2}{\theta-1}}\,dx \right)^{\frac{\theta-1}{\theta}} \left(\int_{|x|<R} |u|^{\theta p_2}\,dx\right)^{\frac{1}{\theta}}+\frac{1}{R^{b_2}}\int_{|x|>R} |u|^{p_2}\,dx,
\end{split}
\end{align}
where $\theta>1$ satisfies that
$$
\frac{\theta b_2}{\theta-1}<N, \quad \theta p_2 <\frac{2N}{N-2}.
$$
Let $\{u_n\} \subset X$ be a bounded sequence. Then there exists $u \in X$  such that $u_n \wto u$ in $X$ as $n \to \infty$. Note that $X \subset D^{1,2}(\R^N)$ and $D^{1,2}(\R^N)$ is locally compactly embedded into $L^p(\R^N)$ for any $2<p<\frac{2N}{N-2}$. Utilizing \eqref{c111'} with $R>0$ large enough, we then have that $u_n \to u$ in $L^{p_2}(\R^N, |x|^{-b_2}dx)$ as $n \to \infty$ and the proof is completed.
 \end{proof}

We are now able to prove the existence of ground states to \eqref{equ1} for $\omega=0$. For this, we shall introduce the following minimization problems,
\begin{align} \label{min2}
m_0:=\inf_{u \in P} E(u),
\end{align}
where
$$
P:=\left\{ u \in X \backslash \{0\} : K(u)=0\right\}.
$$

\begin{lem} \label{existence2}
Let \eqref{H} hold, $N \geq 3$, $p_1<p_2$ and $b_1<b_2$. Then there exist ground states in $X$ to \eqref{equ1} for $\omega=0$.
\end{lem}
\begin{proof}
We argue as in the proof of Lemma \ref{existence1}. It suffices to assert that there exist nonnegative minimizers to \eqref{min2}. Let $\{u_n\} \subset P$ be a sequence that minimizes \eqref{min2}, i.e. $E(u_n)+o_n(1)=m_0$.

Observe first that
$$
p_1 \theta_1 \theta_2 + 2\left(\frac{\theta_1(1-\theta_2)p_1}{2}+\frac{(1-\theta_1)(N-b_2)}{N-2}\right)=p_2,
$$
where $0<\theta_1, \theta_2<1$ are the constants decided in Lemma \ref{em}. It then follows that
$$
\theta_0:=\theta_1 \theta_2 + \frac{\theta_1(1-\theta_2)p_1}{2}+\frac{(1-\theta_1)(N-b_2)}{N-2}>1.
$$
Owing to \eqref{ineq1} and using $K(u_n)=0$, we infer that
\begin{align*}
&\int_{\R^N} |\nabla u_n|^2 \,dx+\frac{N(p_1-1)+2b_1}{2p_1} \int_{\R^N}|x|^{-b_1}|u_n|^{p_1}\, dx
=\frac{N(p_2-2)+2b_2}{2p_2}\int_{\R^N}|x|^{-b_2}|u_n|^{p_2}\, dx \\
& \leq C\left(\int_{\R^N} |\nabla u_n|^2 \,dx+\int_{\R^N}|x|^{-b_1}|u_n|^{p_1}\, dx\right)^{\theta_0}.
\end{align*}
It follows that
\begin{align} \label{bddz}
\int_{\R^N} |\nabla u_n|^2 \,dx +\int_{\R^N} |x|^{-b_1} |u_n|^{p_1} \,dx \geq \frac{1}{C^{\frac{1}{\theta_0-1}}}.
\end{align}
Observe that
\begin{align*} 
m_0=E(u_n)+o_n(1)&=E(u_n)-\frac{2}{N(p_2-2)+2b_2}K(u_n)+o_n(1) \\
&=\left(\frac 12-\frac{2}{N(p_2-2)+2b_2}\right)\int_{\R^N}|\nabla u_n|^2 \,dx \\
& \quad +\frac{1}{p_1} \left(\frac{N(p_2-p_1)+2(b_2-b_1)}{N(p_2-2)+2b_2}\right)\int_{\R^N} |x|^{-b_1}|u_n|^{p_1}\,dx+o_n(1).
\end{align*}
Then we have $m_0>0$ by \eqref{bddz}. At this stage, following closely the line of the proof of Lemma \ref{existence1}, we are able to get the desired conclusion. Thus, the proof is completed.
\end{proof}

\begin{lem} \label{symmetry}
Let $u$ be a ground state to \eqref{equ1} for $\omega \geq 0$. Then $u$ is positive, radially symmetric and decreasing.
\end{lem}
\begin{proof}
Let us first introduce the definition of polarization of measurable functions. Here we denote by $\mathcal{H}$ the family of all affine closed half spaces in $\R^N$ and denote by $\mathcal{H}_0$ the family of all closed half spaces in $\R^N$, i.e. $H \in \mathcal{H}_0$ if and only if $H \in \mathcal{H}$ and $0$ lies in the hyperplane $\partial H$. For $H \in \mathcal{H}$, we denote by $R_H: \R^N \to \R^N$ the reflection with respect to the boundary of $H$. We define the polarization of a measurable function $u : \R^N \to \R$ with respect to $H$ by
\begin{align*}
u_H(x):=\left\{
\begin{aligned}
\max\left\{u(x), u(R_H(x))\right\}, &\quad x \in H, \\
\min \left\{u(x), u(R_H(x))\right\}, &\quad x \in \R^N \backslash H.
\end{aligned}
\right.
\end{align*}
For simplicity, we shall only prove the result for the case $\omega>0$. Let $u \in H^1(\R^N)$ be a ground state to \eqref{equ1} for $\omega>0$, which is indeed a minimizer to \eqref{ming1}. In view of \cite[Lemmas 2.2-2.3]{BWW}, we see that, for any $H \in \mathcal{H}_0$,
$$
\int_{\R^N} |\nabla u_H|^2 \,dx =\int_{\R^N} |\nabla u|^2 \, dx, \quad \int_{\R^N} |u_H|^2 \,dx =\int_{\R^N} |u|^2 \, dx
$$
and
$$
 \int_{\R^N} |x|^{-b_1}|u_H|^{p_1}\,dx=\int_{\R^N} |x|^{-b_1}|u|^{p_1}\,dx, \quad  \int_{\R^N} |x|^{-b_1}|u_H|^{p_2}\,dx=\int_{\R^N} |x|^{-b_2}|u|^{p_2}\,dx.
$$
As a consequence, we find that $u_H \in H^1(\R^N)$ is also a minimizer to \eqref{ming1} for any $H \in \mathcal{H}_0$. It follows from \cite[Theorem 1]{Sch} that there exist a sequence $\{H_n\} \subset \mathcal{H}_0$ and a sequence $\{u_n\} \subset X$ such that $u_n \to u^*$ in $L^p(\R^N)$ as $n \to \infty$ for any $ 1 \leq p<\infty$, where $u^*$ denotes the symmetric-decreasing rearrangement of $u$ and the sequence $\{u_n\}$ is defined by
$$
u_1:=u, \quad u_{n+1}:=(u_n)_{H_1H_2 \cdots H_{n+1}}.
$$
Therefore, we conclude that
$$
\int_{\R^N} |\nabla u^*|^2 \,dx  \leq \int_{\R^N} |\nabla u|^2 \, dx, \quad \int_{\R^N} |u^*|^2 \,dx  = \int_{\R^N} |u|^2 \, dx,
$$
and
$$
 \int_{\R^N} |x|^{-b_1}|u^*|^{p_1}\,dx=\int_{\R^N} |x|^{-b_1}|u|^{p_1}\,dx, \quad  \int_{\R^N} |x|^{-b_2}|u^*|^{p_2}\,dx=\int_{\R^N} |x|^{-b_2}|u|^{p_2-2}\,dx.
$$
It then yields that $S_{\omega}(u^*) \leq S_{\omega}(u)=m_{\omega}$. Note that $K(u)=0$, then $K(u^*) \leq 0$. Hence we have that there exists a constant $0<t_{u^*}\leq 1$ such that $K(u^*_{t_{u^*}})=0$ by Lemma \ref{maximum}. Further, we are able to derive that $u^* \in H^1(\R^N)$ is a minimizer to \eqref{min2} and
$$
\int_{\R^N} |\nabla u^*|^2 \,dx=\int_{\R^N} |\nabla u|^2 \, dx.
$$
In view of \cite[Theorem 1.1]{BZ}, we deduce that $u$ is nonnegative, radially symmetric and decreasing up to translations. Using the maximum principle, we then have the existence of positive ground states. Thus the proof is completed.
\end{proof}

\begin{proof}[Proof of Theorem \ref{decay}]
By applying standard bootstrap arguments, we first get that $u \in C(\R^N) \cap C^2(\R^N\backslash\{0\})$ and $u(x) \to 0$ as $|x| \to \infty$. Let us first show that $u(x) \sim |x|^{-\beta}$ as $|x| \to \infty$ by adapting some ingredients from \cite{DS}. For $R>0$, there holds $-\Delta (|x|^{2-N})=0$ in $\R^N \backslash B_R(0)$. On the other hand, since $p_1<p_2$ and $u(x) \to 0$ as $|x| \to \infty$, then $-\Delta u \leq 0$ in $\R^N \backslash B_R(0)$ for any $R>0$ large enough. From the maximum principle, we then obtain that $u(x) \leq C |x|^{2-N}$ in $\R^N \backslash B_R(0)$. This means that $u(r) \leq C r^{2-N}$ for any $r>0$. We now prove that $u(r) \leq C r^{\frac{b_1-2}{p_1-2}}$ for any $r>0$. Since $u \in X$ is positive and radially symmetric, then \eqref{equ12} can be rewritten as
\begin{align}\label{zequ1}
-u_{rr}-\frac{N-1}{r}u_r+r^{-b_1}u^{p_1-1}=r^{-b_2}u^{p_2-1}.
\end{align}
It is simple to conclude that $u_{rr}>0$ for any $r>0$ large enough. This shows that $u_r$ is increasing and $\lim_{r \to \infty}|u_r(r)|=0$. Multiplying \eqref{zequ1} by $u_r$ and integrating on $[t, \infty)$, we have that
\begin{align} \label{d0}
-\int_t^{\infty} u_{rr}u_r \,dr +\int_t^{\infty}r^{-b_1}u^{p_1-1} u_r \, dr=\int_t^{\infty}r^{-b_2}u^{p_2-1} u_r \, dr + (N-1)\int_t^{\infty}\frac{u_r^2}{r}\,dr.
\end{align}
Note that $u(r) \to 0$ and $u_r(r) \to 0$ as $r \to \infty$. Therefore, from \eqref{d0}, we obtain that
\begin{align*}
f(t):=\frac 12 u_t^2-\frac {1}{p_1} t^{-b_1}u^{p_1} +\frac {1}{p_2} t^{-b_2}u^{p_2}&=\int_t^{\infty} \frac{N-1}{r} u_r^2-\frac {b_1} {p_1} r^{-b_1-1}u^{p_1} + \frac {b_2}{p_2} r^{-b_2-1} u^{p_2} \, dr,
\end{align*}
from which we get that
\begin{align} \label{nonnegative}
\frac{d}{dr} f(r)= -\frac{b_2}{r}f(r)-\frac{2(N-1)-b_2}{2r} u_r^2-\frac{b_2-b_1}{p_1}r^{-b_1}u^{p_1}.
\end{align}
This means that
$$
\frac{d}{dr} \left(r^{b_2} f(r) \right)=-\frac{(2(N-1)-b_2)r^{b_2-1}}{2} u_r^2-\frac{b_2-b_1}{p_1}r^{b_2-b_1}u^{p_1}<0.
$$
As a result, we derive that $r^{b_2} f(r)$ is decreasing for any $r>0$. It then follows that $f(r)$ is decreasing for any $r>0$. Note that $f(r) \to 0$ as $r \to \infty$. Hence there holds that $f(r)>0$ for any $r>0$ large enough. Thus we conclude that $u_r^2 \geq r^{-b_1} u^{p_1}/{p_1}$ for any $r>0$ large enough, because of $u(r) \to 0$ as $r \to \infty$ and $p_1<p_2$. Then we see that, for any $r>0$ large enough,
\begin{align} \label{d01}
\left(r^{\frac {b_1}{2}} u^{\frac{2-p_1}{2}}\right)'=\frac {b_1}{2}r^{\frac {b_1}{2} -1} u^{\frac {2-p_1}{2}} -\frac{p_1-2}{2} r^{\frac {b_1} {2}}u^{-\frac{p_1}{2}} u_r \geq \frac{p_1-2}{2} r^{\frac{b_1}{2}}u^{-\frac{p_1}{2}} |u_r| \geq \frac{p_1-2}{2 \sqrt {p_1}}.
\end{align}
On the other hand, there holds that, for any $r>0$,
\begin{align} \label{d02}
\left(r^{\frac{b_1}{2}} u^{\frac{2-p_1}{2}}\right)'=\frac {b_1}{2}r^{\frac {b_1}{2} -1} u^{\frac {2-p_1}{2}} -\frac{p_1-2}{2} r^{\frac {b_1}{2}}u^{-\frac{p_1}{2}} u_r  \geq \frac {b_1}{2}r^{\frac {b_1}{2} -1} u^{\frac {2-p_1}{2}}.
\end{align}
Combining \eqref{d01} and \eqref{d02}, we conclude that there exists a constant $C>0$ such that, for any $r>0$,
\begin{align} \label{d1}
\left(r^{\frac {b_1}{2}} u^{\frac{2-p_1}{2}}\right)' \geq C.
\end{align}
Therefore, from \eqref{d1}, there holds that $ u(r) \leq C r^{\frac{b_1-2}{p_1-2}}$ for any $r>0$. Consequently, we have that $u(r) \leq C r^{-\beta}$ for any $r>0$.

Let us define $v(r):=r^{\beta}u(r)$ for $r>0$. We are going to prove that there exists a constant $l>0$ such that $v(r) \to l$ as $r \to \infty$. From the discussions above, we see that $v$ is bounded. In addition, by \eqref{zequ1}, it is not hard to verify that $v$ solves the following equation,
\begin{align} \label{zequ2}
v_{rr}-\frac{2\beta +1-N}{r} v_r=\frac{\beta \left(N-\beta -2\right)}{r^2} v+r^{-b_1+\beta(2-p_1)}v^{p_1-1}-r^{-b_2+\beta(2-p_2)}v^{p_2-1}.
\end{align}
First we consider the case that $\beta=N-2$. In this case, by \eqref{zequ2}, then $v$ satisfies the equation
\begin{align*}
v_{rr}-\frac{N-3}{r} v_r=r^{-b_1+(N-2)(2-p_1)}v^{p_1-1}-r^{-b_2+(N-2)(2-p_2)}v^{p_2-1}.
\end{align*}
Define $w(t):=v(r)$ for $t=\beta_0 r^{N-2}/(N-2)$ and
$
\beta_0=\left(N-2\right)^{\frac{b_1-2+p_1(N-2)}{b_1-2+(N-2)(p_1-2)}}>0.
$
Therefore, we have that $w$ enjoys the equation
\begin{align}\label{iddd1}
w''=t^{-\frac{b_1-2+p_1(N-2)}{N-2}}w^{p_1-1}-\left(N-2\right)^{\frac{(N-2)\left((b_1-b_2)+(N-2)(p_2-p_1)\right)}{b_1-2+(N-2)(p_1-2)}}t^{-\frac{b_2-2+p_2(N-2)}{N-2}}w^{p_2-1}.
\end{align}
It then follows that $w''(t) >0$ for any $t>0$ large enough, because of $p_1<p_2$ and $0<b_1-2+p_1(N-2)<b_2-2+p_2(N-2)$.
Note that $w$ is bounded, then $w'(t) \to 0$ as $t \to \infty$. As a consequence, by \eqref{iddd1}, we get that
$$
-w'(t)=\int_t^{\infty} \left(s^{-\frac{b_1-2+p_1(N-2)}{N-2}}w^{p_1-1}(s)-\beta_1s^{-\frac{b_2-2+p_2(N-2)}{N-2}}w^{p_2-1}(s)\right) \,ds,
$$
where
$$
\beta_1:=\left(N-2\right)^{\frac{(N-2)\left((b_1-b_2)+(N-2)(p_2-p_1)\right)}{b_1-2+(N-2)(p_1-2)}}>0.
$$
This gives that $w'(t) <0$ for any $t>0$ large enough. Furthermore, there holds that
\begin{align*}
w(t)&=\int_t^{\infty}\int_{\tau}^{\infty} s^{-\frac{b_1-2+p_1(N-2)}{N-2}}w^{p_1-1}(s)-\beta_1s^{-\frac{b_2-2+p_2(N-2)}{N-2}}w^{p_2-1}(s) \,ds d\tau \\
&=\int_t^{\infty}(s-t)\left(s^{-\frac{b_1-2+p_1(N-2)}{N-2}}w^{p_1-1}(s)-\beta_1s^{-\frac{b_2-2+p_2(N-2)}{N-2}}w^{p_2-1}(s)\right) \,ds
\end{align*}
Therefore, we find that, for any $t>0$ large enough,
\begin{align} \label{below}
\begin{split}
w(t)&=\int_{t}^{\infty}(s-t)\left(s^{-\frac{b_1-2+p_1(N-2)}{N-2}}w^{p_1-1}(s)-\beta_1s^{-\frac{b_2-2+p_2(N-2)}{N-2}}w^{p_2-1}(s)\right) \,ds \\
& \leq w^{p_1-1}(t) \int_{t}^{\infty}(s-t)s^{-\frac{b_1-2+p_1(N-2)}{N-2}}\, ds \\
&=w^{p_1-1}(t) t^{2-\frac{b_1-2+p_1(N-2)}{N-2}} \mathcal{B}\left(\frac{b_1-2+p_1(N-2)}{N-2}-2, 2\right),
\end{split}
\end{align}
where $\mathcal{B}(a, b)$ denotes the Beta function for $a,b>0$. If $\beta=N-2$ and $(2-b_1)/(p_1-2) \neq N-2$, then
$$
\frac{b_1-2+p_1(N-2)}{N-2} > 2.
$$
It then follows from \eqref{below} that there exists a constant $l>0$ such that $w(t) \to l$ as $t \to \infty$. Otherwise, we can reach a contradiction. This in turn leads to  $v(r) \to l$ as $r \to \infty$.

Next, we consider the case $\beta>N-2$.

If $N \geq 4$, we define $w(t):=v(r)$ for $r=\exp\left(\frac{t}{|N-2\beta -1|}\right)$ and $t>0$. In view of \eqref{zequ2}, we then see that
\begin{align} \label{idd1}
w''=\frac{1}{\left(N-2\omega -1\right)^2} f(w),
\end{align}
where
$$
f(w):=\beta\left(N-\beta-2\right) w+\exp\left(\frac{\left(2-b_1+\beta(2-p_1)\right)t}{|N-2\beta -1|}\right) w^{p_1-1}-\exp\left(\frac{\left(2-b_2+\beta(2-p_2)\right)t}{|N-2\beta -1|}\right)w^{p_2-1}.
$$
Integrating \eqref{idd1} on $[t, \infty)$, we then get that
\begin{align} \label{idd2}
-w'(t) =\frac{1}{\left(N-2\beta -1\right)^2} \int_t^{\infty}f(w(s)) \,ds.
\end{align}
Note that $\beta=\max\{(2-b_1)/(p_1-2), N-2\}$. If $\beta > N-2$, then $N-2\beta-1<0$ and
$$
\frac{2-b_2+\beta(2-p_2)}{|N-2\beta-1|}<\frac{2-b_1+\beta(2-p_1)}{|N-2\beta-1|} \leq 0.
$$
Then there holds that, for any $t>0$ large enough,
$$
f(w(t)) \leq -\frac{\beta\left(\beta-N+2\right)}{2} w(t)<0.
$$
Taking into account \eqref{idd2}, we then obtain that $w'(t)>0$ for any $t>0$ large enough. Due to $w>0$, then $w(t) \to l$ as $t \to \infty$. This shows that $v(r) \to l$ as $r \to \infty$.
In this case, if $N=3$, then \eqref{zequ2} reduces to
\begin{align} \label{vequ}
v_{rr}-\frac{2(\beta-1)}{r} v_r=\frac{\beta \left(1-\beta\right)}{r^2} v+r^{-b_1+\beta(2-p_1)}v^{p_1-1}-r^{-b_1+\omega(2-p_1)}v^{p_2-1}.
\end{align}
Utilizing a similar way as before, we can also prove that  $v(r) \to l$ as $ r \to \infty$ for some $l>0$ when $\beta=1$ or $\beta>1$. Hence we have the desired result, i.e. $u(x) \sim |x|^{-\beta}$ as $|x| \to \infty$.

Let us now turn to demonstrate that $u(x) \sim |x|^{2-N} \left(\ln |x| \right)^{\frac{2-N}{2-b_1}}$ as $|x| \to \infty$ for $p_1 = (2N-2-b_1)/(N-2)$. To do this, we shall follow some ideas from \cite{DSW}. Let us first define $v(r):=r^{N-2}u(r)$ for any $r>0$. In virtue of \eqref{zequ2} and $\beta=N-2$, we immediately know that $v$ satisfies the following equation,
\begin{align*}
v_{rr}-\frac{N-3}{r} v_r=r^{-2}v^{p_1-1}-r^{-b_2+(N-2)(2-p_2)}v^{p_2-1}.
\end{align*}
Define $w(t):=v(r)$ and $r=e^t$ for any $t>0$. Then we find that
\begin{align*}
w''-(N-2)w'=w^{q-1}-e^{(2-b)t+(N-2)(2-p)t}w^{p-1}.
\end{align*}
It then follows that
\begin{align} \label{zequ11}
\left(w' e^{-(N-2)t}\right)'=e^{-(N-2)t}\left(w^{p_1-1}-e^{(2-b_2)t+(N-2)(2-p_2)t}w^{p_2-1}\right).
\end{align}
By integrating \eqref{zequ11} on $[t, \infty)$, we then have that
$$
-w'(t)=e^{(N-2)t}\int_{t}^{\infty}e^{-(N-2)s}\left(w^{p_1-1}(s)-e^{(2-b_2)s+(N-2)(2-p_2)s}w^{p_2-1}(s)\right) \, ds.
$$
This then yields that $w'(t)<0$ for any $t>0$ large enough. Therefore, we get that $-w'(t) \leq C w^{p_1-1}(t)$ for any $t>0$ large enough. As a consequence, we have that $\left(w^{2-p_1}(t)\right)' \leq C(p_1-2)$ for any $t>0$ large enough, from which we obtain that
$$
w(t) \geq \left(C(p_1-2) t +w^{2-p_1}(t_0)\right)^{\frac {1}{2-p_1}} \geq \widetilde{C} t^{\frac{2-N}{2-b_1}}, \quad t>t_0>0,
$$
because of $p_1 = (2N-2-b_1)/(N-2)$. This yields that, for any $r>0$ large enough,
$$
u(r) \geq \widetilde{C} r^{2-N} \left(\ln r\right)^{\frac{2-N}{2-b_1}}.
$$
In the following, we are going to show the upper bound of $u$. To this end, we first set
$$
S(r):=kr^{2-N} \left(\ln r\right)^{\frac{2-N}{2-b_1}}, \quad r>0,
$$
where $k>0$ is a constant defined by
$$
k:=\left(\frac{2-b_1}{(N-2)^2}\right)^{\frac{N-2}{2-b_1}}.
$$
From direct computations, we see that
$$
-\Delta S + |x|^{-b_1} S^{p_1-1}=-\frac{(N-b_1) |x|^{-b_1}}{(2-b_1)(N-2)} \frac{S^{p_1-1}}{\ln |x|}.
$$
Let $\eta>0$ be a constant to be determined later. It follows from the above equation that
$$
-\Delta (\eta S) +  |x|^{-b_1} (\eta S)^{p_1-1}= |x|^{-b_1}(\eta^{p_1-1}-\eta)S^{p_1-1} - \frac{\eta (N-b_1) |x|^{-b_1}}{(2-b_1)(N-2)}\frac{S^{p_1-1}}{\ln |x|}.
$$
Therefore 
\begin{align*}
&-\Delta(u-\eta S) + \frac{ |x|^{-b_1}\left(u^{p_1-1}-(\eta S)^{p_1-1}\right)}{u-\eta S} \left(u-\eta S\right)=-\Delta u+ |x|^{-b_1} u^{p_1-1} + \left(\Delta (\eta S)- |x|^{-b_1}(\eta S)^{p_1-1}\right) \\
& = |x|^{-b_2}u^{p_2-1} - |x|^{-b_1}(\eta^{p_1-1}-\eta)S^{p_1-1} + \frac{\eta (N-b_1) |x|^{-b_1}}{(2-b_1)(N-2)}\frac{S^{p_1-1}}{\ln |x|} \\
&= |x|^{-b_1}S^{p_1-1} \left(|x|^{-b_2+b_1}\frac{u^{p_1-1}}{S^{p_1-1}}-(\eta^{p_2-1}-\eta)+ \frac{\eta (N-b)}{(2-b)(N-2)\ln |x|}\right).
\end{align*}
Note that $u(x) \leq C |x|^{-\beta}$ for $x \in \R^N$ and $p_1<p_2$,  then
$$
\frac{u^{p_2-1}(x)}{S^{p_1-1}(x)} \to 0 \quad \mbox{as} \,\, |x| \to \infty.
$$
As a consequence, there exists a constant $R>0$ large enough such that
$$
-\Delta(u-\eta S) +c(x) \left(u-\eta S\right) <0,
$$
where $\eta>0$ is a constant such that $\eta<\eta^{p_1-1}$ and
$$
c(x):=\frac{|x|^{-b_1}\left(u^{p_1-1}(x)-(\eta S)^{p_1-1}(x)\right)}{u(x)-\eta S(x)}, \quad x \in \R^N
$$
Using the maximum principle, we then have that $u(x) \leq \widehat{C} S(x)$ for any $|x| \geq R$. This completes the proof.
\end{proof}

\begin{lem} \label{nonexistence}
Let \eqref{H} holds and $N \geq 2$. Then there exists no solutions in $H_{rad}^1(\R^N)$ to \eqref{equ1} for $\omega<0$.
\end{lem}
\begin{proof}
Let $u \in H_{rad}^1(\R^N)$ be a solution to \eqref{equ1} for $\omega<0$. Then
$$
-\Delta u + \left(|x|^{-b_1} |u|^{p_1-2} - |x|^{-b_2} |u|^{p_2-2} \right) u=-\omega u.
$$
As an application of Lemma \ref{hrc}, we have that
$$
|x|^{-b_1} |u|^{p_1-2} + |x|^{-b_2} |u|^{p_2-2} \leq C |x|^{-b_1-\frac{(N-1)p_1}{2}}, \quad |x| \geq 1.
$$
This readily implies that
$$
|x|^{-b_1} |u|^{p_1-2} + |x|^{-b_2} |u|^{p_2-2} =o(|x|^{-\frac{N-1}{2}}) \quad \mbox{as}\,\, |x| \to \infty.
$$
It then follows from \cite{Ka} that $u=0$. This completes the proof.
\end{proof}

\begin{proof} [Proof of Theorem \ref{thm1}]
It follows from Theorem \ref{decay} that if $u \in X$ is a positive, radially symmetric and decreasing solution to \eqref{equ12}, then $u \in H^1(\R^N)$ for $N \geq 5$. Therefore, from Lemmas \ref{existence1}, \ref{existence2},\ref{symmetry} and \ref{nonexistence}, we have the desired conclusions and the proof is completed.
\end{proof}

\section{Uniqueness and Non-degeneracy} \label{unno}

In this section, we are going to discuss uniqueness and non-degeneracy of solutions to \eqref{equ1} and present the proofs of Theorems \ref{thmunique} and \ref{nond}. In the following, we shall always assume that \eqref{H} holds.

\subsection{Uniqueness}

Let $u \in H^1(\R^N)$ be a positive, radially symmetric and decreasing solution to \eqref{equ1}. We shall first introduce the following ordinary differential equation satisfied by $u$,
\begin{align} \label{ode1}
\left\{
\begin{aligned}
&u_{rr}+\frac{N-1}{r} u_r - \omega u +r^{-b_2}u^{p_2-1} -r^{-b_1}u^{p_1-1}=0,\\
&u(0)=a>0, \quad \lim_{r \to \infty} u(r)=0.
\end{aligned}
\right.
\end{align}
Define the corresponding Pohozaev quantity by
\begin{align} \label{defj}
\hspace{-0.2cm}J(r; u):=\frac 12 A(r)u_r(r)^2+B(r)u_r(r)u(r) +\frac 12 C(r)u(r)^2+\frac{1}{p_2}A(r)r^{-b_2}u(r)^{p_2}-\frac{1}{p_1}A(r)r^{-b_1}u(r)^{p_1},
\end{align}
where $A, B$ and $C : \R^+ \to \R$ are functions determined later. Then we obviously see that
\begin{align*}
\frac{d}{dr} J(r;u)&=\left(\frac 12 A_r(r)-\frac{N-1}{r}A(r)+B(r)\right)u_r(r)^2 \\
& \quad + \left(\omega A(r)-\frac{N-1}{r}B(r) +B_r(r)+C(r)\right)u_r(r)u(r) +\left(\omega B(r)+ \frac 12 C_r(r)\right)u(r)^2\\
& \quad + \left(-B(r)r^{-b_2}+\frac{1}{p_2}\left(A(r)r^{-b_2}\right)_r\right) u(r)^{p_2}+\left(B(r)r^{-b_1}-\frac{1}{p_1}\left(A(r)r^{-b_1} \right)_r\right) u(r)^{p_1}.
\end{align*}
Let $A(r)$, $B(r)$ and $C(r)$ satisfy
$$
\frac 12 A_r(r)-\frac{N-1}{r}A(r)+B(r)=0, \quad \omega A(r)-\frac{N-1}{r}B(r) +B_r(r)+C(r)=0.
$$
$$
B(r)r^{-b_1}-\frac{1}{p_1}\left(A(r)r^{-b_1}\right)_r=0.
$$
It then follows that
$$
\left(\frac 12 A_r(r)-\frac{N-1}{r}A(r)\right)r^{-b_1} +\frac{1}{p_1}\left(A(r)r^{-b_1}\right)_r=0.
$$
Therefore, we get that
$$
A(r)=r^{\frac{2b_1+2(N-1)p_1}{p_1+2}},
\quad B(r)=\frac{2(N-1)-b_1}{p_1+2}r^{\frac{2(b_1-1)+(2N-3)p_1}{p_1+2}},
$$
$$
C(r)=-\omega r^{\frac{2b_1+2(N-1)p_1}{p_1+2}}+\frac{\left(2(N-1)-b_1\right)\left(2(N-b_1)-p_1(N-2)\right)}{(p_1+2)^2}r^{\frac{2(b_1-2)+2(N-2)p_1}{p_1+2}}.
$$
Define
\begin{align*}
G(r)&:=\omega B(r)+ \frac 12 C_r(r)\\
&=\omega\left(\frac{(N-1)(2-p_1)-2b_1}{p_1+2}\right)r^{\frac{2(b_1-1)+(2N-3)p_1}{p_1+2}}\\
& \quad +\frac{\left(2(N-1)-b_1\right)\left(2(N-b_1)-p_1(N-2)\right)\left((b_1-2)+(N-2)p_1\right)}{(p_1+2)^3}r^{\frac{2(b_1-3)+(2N-5)p_1}{p_1+2}},
\end{align*}
and
\begin{align*}
H(r)&:=-B(r)r^{-b_2}+\frac{1}{p_2}\left(A(r)r^{-b_2}\right)_r\\
&=\left(-\frac{2(N-1)-b_1}{p_1+2}+\frac{2b_1+2(N-1)p_1}{p_2(p_1+2)}-\frac{b_2}{p_2}\right)r^{\frac{2(b_1-1)+(2N-3)p_1}{p_1+2}-b_2}.
\end{align*}
As a consequence, there holds that
\begin{align} \label{drj}
\frac{d}{dr} J(r;u)=G(r)u(r)^2+H(r)u(r)^{p_2}.
\end{align}

\begin{lem} \label{unique11}
Let $u$ be a positive and decreasing solution to \eqref{ode1} such that $\displaystyle\lim_{r \to 0} ru_r(r)$ exists and $J(r; u) \to 0$ as $r \to \infty$. Then $J(\cdot; u) \not\equiv 0$. Moreover, suppose that \eqref{uniquec} holds,
then $J(r; u) \geq 0$ for any $r \geq 0$.
\end{lem}
\begin{proof}
Since $N \geq 3$, then
$$
\frac{2b_1+2(N-1)p_1}{p_1+2}>2.
$$
Using the assumption that $\displaystyle\lim_{r \to 0} ru_r(r)$ exists, and owing to \eqref{defj}, we then get $J(r; u) \to 0$ as $r \to 0$. In addition, from \eqref{drj}, we see that $\frac{d}{dr}J(r; u)>0$ for any $r>0$ small enough. It then follows that $J(r;u)>0$ for any $r>0$ small enough. Thus we know that $J(\cdot; u) \not\equiv 0$.  Next, we shall prove that $J(r; u) \geq 0$ for any $r \geq 0$. Since $u$ is decreasing on $(0, \infty)$, by \eqref{uniquec}, then we are able to derive that there exists a unique $r_0>0$ such that $\frac{d}{dr} J(r_0; u)=0$. This is indeed justified by the fact that the function $f: \R^+ \to \R$ defined by $f(r):=-\omega_1+\omega_2 r^{-2} +\omega_3 r^{-b_2} u(r)^{p_2-2}$ admits only one zero in $(0, \infty)$, where
$$
\omega_1=-\omega\left(\frac{(N-1)(2-p_1)-2b_1}{p_1+2}\right) \geq 0,
$$
$$
\omega_2=\frac{\left(2(N-1)-b_1\right)\left(2(N-b_1)-p_1(N-2)\right)\left((b_1-2)+(N-2)p_1\right)}{(p_1+2)^3}>0,
$$
$$
\omega_3=-\frac{2(N-1)-b_1}{p_1+2}+\frac{2b_1+2(N-1)p_1}{p_2(p_1+2)}-\frac{b_2}{p_2} \geq 0.
$$
Since $J(r; u) \to 0$ as $r \to 0$, $J(r; u)>0$ for any $r>0$ small enough and $J(r; u) \to 0$ as $r \to \infty$. Therefore, we know that $J(r; u) \geq 0$ for any $r \geq 0$. This completes the proof.
\end{proof}

\begin{lem} \label{unique12}
Let $u, v$ be two positive solutions to \eqref{ode1} such that $\displaystyle\lim_{r \to 0}ru_r(r)$ and $\displaystyle\lim_{r \to 0} rv_r(r)$ exist. Then there holds that
$$
\frac{d}{dr}\left(\frac{v(r)}{u(r)}\right)=\frac{1}{u(r)^2}\int_{0}^{r}\frac{s^{N-1}}{r^{N-1}}\left(s^{-b_2}\left(u(s)^{p_2-2}-v(s)^{p_2-2}\right)-s^{-b_1}\left(u(s)^{p_1-2}-v(s)^{p_1-2}\right)\right)u(s)v(s)\,ds.
$$
\end{lem}
\begin{proof}
Since $u, v$ are solutions to \eqref{ode1}, then
$$
\left(r^{N-1}u_r\right)_r-\omega r^{N-1} u+r^{N-1-b_2}u^{p_2-1}-r^{N-1-b_1}u^{p_1-1}=0,
$$
$$
\left(r^{N-1}v_r\right)_r-\omega r^{N-1} v+r^{N-1-b_2}v^{p_2-1}-r^{N-1-b_1}v^{p_1-1}=0.
$$
It then follows that
\begin{align*}
r^{N-1} \left(u(r)v_r(r)-v(r)u_r(r)\right)&=\int_{0}^{r}s^{N-1-b_2}\left(u(s)^{p_2-2}-v(s)^{p_2-2}\right)u(s)v(s)\,ds\\
&\quad -\int_{0}^{r}s^{N-1-b_1}\left(u(s)^{p_1-2}-v(s)^{p_1-2}\right)u(s)v(s)\,ds.
\end{align*}
This leads to the desired conclusion and the proof is completed.
\end{proof}

\begin{lem}\label{lem13}
Let $u, v$ be two positive solutions to \eqref{ode1} such that $\displaystyle\lim_{r \to 0} ru_r(r)$, $\displaystyle\lim_{r \to 0} rv_r(r)$ exist, $u(0)<v(0)$ and $J(r;u) \geq 0$ for any $r>0$. Suppose that \eqref{uniquec} holds, $b_1<b_2$ and $p_1<p_2$, then
$$
\frac{d}{dr}\left(\frac{v(r)}{u(r)}\right)<0, \quad \forall \,\, r>0.
$$
\end{lem}
\begin{proof}
Suppose by contradiction that the desired conclusion does not hold. Define $w(r):=\frac{v(r)}{u(r)}$ for any $r>0$. Then we may assume that there exists $r_0>0$ such that $w_r(r_0)>0$. Since $0<u(0)<v(0)$ and $0<b_1<b_2$, by Lemma \ref{unique12}, then there exists $r_1>0$ small enough such that $w_r(r_1)<0$. Therefore, we know that there exists $r_*>0$ such that $w_r(r_*)=0$ and $w_r(r)<0$ for any $0<r<r_*$. 
Define
\begin{align} \label{defxr}
X(r):=w(r)^2J(r; u)-J(r; v), \quad r>0.
\end{align}
By the definitions of $J(r; u)$ and $J(r; v)$, we then have that
\begin{align} \label{x}
\begin{split}
X(r)&=\frac 12 A(r)\left(\frac{v(r)^2u_r(r)^2}{u(r)^2}-v_r(r)^2\right)+B(r)\left(\frac{v(r)^2u_r(r)}{u(r)}-v_r(r)v(r)\right) \\
& \quad +\frac{1}{p_2} A(r)r^{-b_2} v(r)^2\left(u(r)^{p_2-2}-v(r)^{p_2-2}\right)-\frac{1}{p_1} A(r)r^{-b_1} v(r)^2\left(u(r)^{p_1-2}-v(r)^{p_1-2}\right).
\end{split}
\end{align}
This clearly indicates that $X(r) \to 0$ as $r \to 0$. In addition, using the fact that $w_r(r_*)=0$ and Lemma \ref{unique12}, we find that
\begin{align} \label{xr} \nonumber
X(r_*)&=\frac{1}{p_2} A(r_*)r_*^{-b_2} v(r_*)^2\left(u(r_*)^{p_2-2}-v(r_*)^{p_2-2}\right)-\frac{1}{p_1} A(r_*)r_*^{-b_1} v(r_*)^2\left(u(r_*)^{p_1-2}-v(r_*)^{p_1-2}\right) \\
& > \left(\frac{1}{p_2}-\frac{1}{p_1}\right) A(r_*) v(r_*)^2 \left(u(r_*)^{p_1-2}-v(r_*)^{p_1-2}\right).
\end{align}
It is simple to calculate that
\begin{align} \label{dxr}
\begin{split}
X_r(r)&=2w(r)w_r(r)J(r; u)+w(r)^2J_r(r; u)-J_r(r; v)\\
&=2w(r)w_r(r)J(r; u)+H(r)\left(u(r)^{p_2-2}-v(r)^{p_2-2}\right)v(r)^2 \\
& \leq H(r)\left(u(r)^{p_2-2}-v(r)^{p_2-2}\right)v(r)^2, \quad \forall\,\, 0<r<r^*.
\end{split}
\end{align}
If $w(r) \geq 1$ for any $0<r \leq r^*$, by \eqref{xr}, then there holds that $X(r_*)>0$. In view of \eqref{dxr}, we have that $X_r(r) \leq 0$ for any $0<r \leq r^*$, because $H(r) \geq 0$ for any $r>0$ under the assumption \eqref{uniquec}. Since $X(r) \to 0$ as $r \to 0$, then $X(r^*) \leq 0$. This is a contradiction. Otherwise, there exists $0<r_2<r^*$ such that $w(r_2)=1$, $w_r(r_2)<0$ and $w_r(r)<0$ for any $0<r<r_2$, by \eqref{x} and \eqref{dxr}, then $X(r_2)>0$ and $X_r(r)<0$ for any $0<r<r_2$. Similarly, we are able to reach a contradiction. Thus the proof is completed.
\end{proof}

\begin{proof}[Proof of Theorem \ref{thmunique}]
Let $u \in H^1(\R^N)$ be a positive, radially symmetric and decreasing solution to \eqref{equ1}. Then there holds that $\displaystyle\lim_{r \to 0} ru_r(r)$ exists and $J(r; u) \to 0$ as $r \to \infty$. Let us now suppose that there exist two distinct positive, decreasing and decreasing solutions $u, v$ to \eqref{ode1} such that $\displaystyle\lim_{r \to 0} ru_r(r)$ and $ \displaystyle\lim_{\to 0} rv_r(r)$ exist, $J(r; u) \to 0$ and $J(r; v) \to 0$ as $r \to \infty$. Suppose further that $u(0)<v(0)$. Using Lemma \ref{unique11}, we first know that $J(r; u) \geq 0$ and $J(r ;v) \geq 0$ for any $r>0$. In addition, there holds that $J(\cdot; u) \not\equiv 0$ and $J(\cdot; v) \not\equiv 0$. Define $w(r):=\frac{v(r)}{u(r)}$ for $r>0$ as previously. From Lemma \ref{lem13}, then $w_r(r)<0$ for any $r>0$. Let us also define $X(r)$ by \eqref{defxr} for $r>0$. Then we see that $X(r) \to 0$ as $r \to 0$ and $X(r) \to 0$ as $r \to \infty$. If $u(r) < v(r)$ for any $r>0$, then $X_r(r)<0$ for any $r>0$. This is impossible, because $X(r) \to 0$ as $r \to 0$ and $X(r) \to 0$ as $r \to \infty$. Otherwise, reasoning as the proof of Lemma \ref{lem13}, we can also obtain a contradiction. Thus the proof is completed.
\end{proof}

\subsection{Non-degeneracy}
{
To begin with, we shall present some basic results. It is clear that the Laplacian operator can be represented as the following form in term of radial and angular variables,
$$
\Delta=\partial_{rr} + \frac{N-1}{r} \partial_r + \frac{\Delta_{\mathbb{S}^{N-1}}}{r^2}.
$$
The eigenvalues of the operator $\Delta_{\mathbb{S}^{N-1}}$ are given by $-k(k+N-2)$, whose multiplicity are given by
$$
\left(
\begin{aligned}
N+&k-1\\
&k
\end{aligned}
\right)
-
\left(
\begin{aligned}
N&+k-3\\
&k-2
\end{aligned}
\right).
$$
The corresponding eigenfunctions are denoted by $\{Y_k\}$, which are spherical harmonic functions that satisfy the equation
$$
\Delta_{\mathbb{S}^{N-1}}Y_k=-k(k+N-2)Y_k, \quad k=1,2, \cdots.
$$
In particular, when $k=1$, we see that
\begin{align} \label{fes}
\Delta_{\mathbb{S}^{N-1}} \frac{x_j}{r}=-(N-1)\frac{x_j}{r}, \quad j=1, 2, \cdots, N.
\end{align}
Define
$$
\mathcal{H}_0:=L_{rad}^2(\R^N), \quad \mathcal{H}_k:=\mbox{span}\left\{ f_k(r) Y_k : f_k \in L^2_{rad}(\R^N)\right\}, \quad k=1,2, \cdots.
$$
This gives an orthogonal decomposition of $L^2(\R^N)$, i.e.
$$
L^2(\R^N)=\bigoplus_{k=0}^{\infty} \mathcal{H}_k.
$$
Let $\mathcal{L}_{+, k} := \mathcal{L}_+|_{\mathcal{H}_k}$ for $u \in L^2_{\mathrm{rad}}(\mathbb{R}^N)$, where $\mathcal{L}_+$ is given by \eqref{L+}. The operator $\mathcal{L}_{+, k}$ acts on $L^2_{\mathrm{rad}}(\mathbb{R}^N)$ via the expression
\[
\mathcal{L}_{+, k} = -\partial_{rr} - \frac{N-1}{r} \partial_r + \frac{k(k+N-2)}{r^2} + \omega + (p_1-1)|x|^{-b_1}|u|^{p_1-2} - (p_2-1)|x|^{-b_2}|u|^{p_2-2},
\]
where $k \in \mathbb{N}_0 = \{0, 1, 2, \dots\}$.

\begin{proof}[Proof of Theorem \ref{nond}]
First, by the definitions of $\mathcal{L}_{+, k}$, we observe that
$$
\mathcal{L}_{+, 0}<\mathcal{L}_{+, 1}< \cdots<\mathcal{L}_{+, k}< \cdots.
$$
Since $n(\mathcal{L}_+)=1$,  then $n(\mathcal{L}_{+,0})=1$. This then gives rise to $\mathcal{L}_{+, k} \geq 0$ for any $k=1,2, \cdots$.

To establish Theorem \ref{nond}, we first show that
\begin{align} \label{nond1}
Ker[\mathcal{L}_{+, 0}]=\{0\} \,\,\, \mbox{or} \,\,\, Ker[\mathcal{L}_{+, 0}]=\mbox{span}\{\Psi_0\},
\end{align}
Since $n(\mathcal{L}_{+, 0})=1$, then $\mathcal{L}_{+, 0}$ has a negative eigenvalue $-\sigma_0$ for some $\sigma_0>0$. Let $\Psi_0$ be an eigenfunction of $\mathcal{L}_{+, 0}$ corresponding to the next eigenvalue $\sigma_1>-\sigma_0$. Hence we have that $\sigma_1 \geq 0$, because of $n(\mathcal{L}_{+})=1$. If $\sigma_1>0$, then $Ker[\mathcal{L}_{+, 0}]=0$. If $\sigma_1=0$, then $\mathcal{L}_{+, 0}[\Psi_0]=0$. Utilizing the ideas presented in the proof of \cite[Theorem 3]{St}, we further know that $\Psi_0$ has exactly one zero in $(0, \infty)$.

To complete the proof, we need to rule out the existence of $\Psi_0$. In order to further discuss, we shall denote by $u_{\omega} \in H^1(\R^N)$ the positive radially symmetric solution to \eqref{equ1} for emphasis the dependence on $\omega$. Then we can write that
\begin{align} \label{equ1a}
-\Delta u_{\omega} + \omega u_{\omega} =|x|^{-b_2} u_{\omega}^{p_2-1}-|x|^{-b_1} u_{\omega}^{p_1-1} \quad \mbox{in} \,\, \R^N,
\end{align}
$$
-\Delta u_{\omega+ \delta} + (\omega+ \delta) u_{\omega+ \delta} =|x|^{-b_2} u_{\omega+ \delta}^{p_2-1}-|x|^{-b_1} u_{\omega+ \delta}^{p_1-1}  \quad \mbox{in} \,\, \R^N.
$$
This clearly leads to
\begin{align*}
&\left(-\Delta + \omega \right) \left(\frac{u_{\omega+ \delta}-u_{\omega}}{\delta}\right) + |x|^{-b_1} \left(\frac{u_{\omega+ \delta}^{p_1-1} -u_{\omega}^{p_1-1}}{\delta}\right) =|x|^{-b_2} \left(\frac{u_{\omega+ \delta}^{p_2-1}-u_{\omega}^{p_2-1}}{\delta}\right)-u_{\omega+\delta}.
\end{align*}
Let $\psi \in H^1(\R^N)$, then
\begin{align*}
&\left\langle \frac{u_{\omega+ \delta}-u_{\omega}}{\delta}, \left(-\Delta + \omega \right)\psi \right\rangle + \left\langle \frac{u_{\omega+ \delta}^{p_1-1} -u_{\omega}^{p_1-1}}{\delta}, |x|^{-b_1} \psi \right\rangle  \\
& =\left\langle\frac{u_{\omega+ \delta}^{p_2-1}-u_{\omega}^{p_2-1}}{\delta}, |x|^{-b_2} \psi \right\rangle-\left \langle u_{\omega+\delta}, \psi \right\rangle.
\end{align*}
It then follows that $\left\langle \partial_{\omega} u_{\omega}, \mathcal{L}_{+, 0} [\psi] \right \rangle=-\left \langle u_{\omega}, \psi \right\rangle$ by taking the limit as $\delta \to 0$. This immediately infers that $u_{\omega} \bot Ker [\mathcal{L}_{+,0}]$.
On the other hand, we observe that 
$$
\mathcal{L}_{+, 0} u_{\omega}=(p_1-2)|x|^{-b_1}u_{\omega}^{p_1-1} -(p_2-2)|x|^{-b_2}u_{\omega}^{p_2-1}u_{\omega}.
$$ 
This then results in
\begin{align} \label{nond2}
(p_1-2)|x|^{-b_1}u_{\omega}^{p_1-1} -(p_2-2)|x|^{-b_2}u_{\omega}^{p_2-1} \bot Ker[\mathcal{L}_{+, 0}].
\end{align}
In addition, we note that
\begin{align*}
\mathcal{L}_{+, 0}\left[\sum_{j=1}^N x_j\partial_{x_j} u_{\omega}\right] &=(-\Delta + \omega)  \left(\sum_{j=1}^N x_j \partial_{x_j} u_{\omega}\right)  \\
& \quad + \left( (p_1-1)|x|^{-b_1}u_{\omega}^{p_1-1}-(p_2-1)|x|^{-b_2}u_{\omega}^{p_2-1} \right)\left(\sum_{j=1}^N x_j \partial_{x_j} u_{\omega}\right) \\
&=-2 \Delta u_{\omega} +b_1 |x|^{-b_1} u_{\omega}^{p_1-1}-b_2 |x|^{-b_2} u_{\omega}^{p_2-1}.
\end{align*}
Since $u_{\omega} \in H^1(\R^N)$ is a solution to \eqref{equ1}, then 
$$
\mathcal{L}_{+, 0}\left[\sum_{j=1}^N x_j\partial_{x_j} u_{\omega}\right] =-2\omega u_{\omega}+(b_1-2) |x|^{-b_1} u_{\omega}^{p_1-1}-(b_2-2) |x|^{-b_2} u_{\omega}^{p_2-1}.
$$
It then follows that
$$
-2\omega u_{\omega}+(b_1-2) |x|^{-b_1} u_{\omega}^{p_1-1}-(b_2-2) |x|^{-b_2} u_{\omega}^{p_2-1} \bot Ker[\mathcal{L}_{+, 0}].
$$
Further, since $u_{\omega} \bot Ker [\mathcal{L}_{+,0}]$, we know that
\begin{align} \label{nond21}
(b_1-2) |x|^{-b_1} u_{\omega}^{p_1-1}-(b_2-2) |x|^{-b_2} u_{\omega}^{p_2-1} \bot Ker[\mathcal{L}_{+, 0}].
\end{align}
Combining \eqref{nond2} and \eqref{nond21}, we then have that
$$
|x|^{-b_2}u_{\omega}^{p_2-1} \bot Ker[\mathcal{L}_{+, 0}].
$$
Let $r_0>0$ the unique zero of $\Psi_0$ in $(0, \infty)$. Then we assume that $\Psi_0(r)<0$ for $0 <r<r_0$ and $\Psi_0(r)>0$ for $r>r_0$. Define
\begin{align*}
\varphi &:=c_0 u_{\omega}-|x|^{-b_2}u_{\omega}^{p_2-1}=u_{\omega} \left(c_0-|x|^{-b_2}u_{\omega}^{p_2-2}\right),  
\end{align*}
where $c_0:=r_0^{-b_2}u^{p_2-1}_{\omega}(r_0)$. Then we see that $\varphi \bot Ker[\mathcal{L}_{+, 0}]$, which shows that $\langle \varphi, \Psi_0 \rangle=0$. However, there holds that $\varphi(x) <0$ for $|x| <r_0$ and  $\varphi(x) >0$ for $|x| >r_0$. This means that $\langle \varphi, \Psi_0 \rangle>0$. We then reach a contradiction. This completes the proof.
\end{proof}
}
\section{Blowup} \label{blowupr}
In this section, we are going to discuss blowup of solutions to the Cauchy problem \eqref{equ} and present the proofs of Theorem \ref{thm2} and \ref{thm3} and Corollary \ref{instability}. In the following, we shall always assume that \eqref{H} holds. First, the space defined in \eqref{a-} is stable under the flow of \eqref{equ1-1}.

\begin{lem} \label{invariant}
Let $u \in C([0, T_{max}), H^1(\R^N))$ be the maximal solution to the Cauchy problem \eqref{equ} with $u_0 \in \mathcal{A}^-_{\omega}$. Then $u(t) \in \mathcal{A}^-_{\omega}$ for any $t \in [0, T_{\max})$. Moreover, there exists $\beta_0>0$ such that $K(u(t)) \leq -\beta_0$ for any $t \in [0, T_{\max})$.
\end{lem}
\begin{proof}
Let $u_0 \in \mathcal{A}^-_{\omega}$ and $u \in C([0, T_{max}), H^1(\R^N))$ be the maximal solution to the Cauchy problem \eqref{equ}. Since $S_{\omega}(u_0) <m_{\omega}$, by the conservation laws, then $S_{\omega}(u(t))<m_{\omega}$ for any $t\in[0, T_{max})$. Assume that there exists $0<t_0<T_{max}$ such that $K(u(t_0))=0$. Then $S_{\omega}(u(t_0)) \geq m_{\omega}$. This is impossible. Hence $\mathcal{A}^-_{\omega}$ is invariant under the flow of the Cauchy problem \eqref{equ}. In particular, $K(u(t)) <0$ for all $0\leq t<T_{max}$. From Lemma \ref{maximum}, it then follows that there exists $0<t_u<1$ such that $K(u_{t_u})=0$. In addition, we see that
$$
S_{\omega}(u)-S_{\omega}(u_{t_u})=(1-t_u) \frac{d}{dt}S_{\omega}(u_t)\mid_{t=\xi} \geq (1-t_u) \frac{d}{dt}S_{\omega}(u_t)\mid_{t=1}=(1-t_u)K(u)>K(u), \quad \xi \in (t_u, 1).
$$
Therefore, by the conservation laws, we conclude that $K(u(t))<S_{\omega}(u_0)-m_{\omega}:=-\beta_0$ for any $t \in [0, T_{max})$ and the proof is completed.
\end{proof}

\begin{lem} \label{tbelow}
Let $u \in C([0, T_{max}), H^1(\R^N))$ be the maximal solution to the Cauchy problem \eqref{equ} with $u_0 \in \mathcal{A}^-_{\omega}$. Then there exists $\eps_0>0$ such that $\|\nabla u(t)\|_2 \geq \eps_0$ for any $t \in [0, T_{\max})$. Moreover, there exists $\delta_0>0$ such that $K(u(t))\leq -\delta_0 \|\nabla u(t)\|_2^2$ for any $t \in [0, T_{\max})$.
\end{lem}
\begin{proof}
Observe that
\begin{align} \label{tbe111}
\begin{split}
K(u(t))&=\frac{N(p_2-2)+2b_2}{2} E(u(t))-\frac{N(p_2-p_1)+2(b_2-b_1)}{2p_1} \int_{\R^N} |x|^{-b_1}|u(t,x)|^{p_1}\,dx\\
& \quad -\left(\frac{N(p_2-2)+2b_2}{4}-1\right) \int_{\R^N} |\nabla u(t,x)|^2 \,dx.
\end{split}
\end{align}
This readily shows that there exists $\eps_0>0$ such that $\|\nabla u(t)\|_2 \geq \eps_0$ for any $t \in [0, T_{max})$. Otherwise, we may suppose that there exists $\{t_n\} \subset [0, T_{max})$ such that $\|\nabla u(t_n)\|=o_n(1)$. Using \eqref{GN} and \eqref{tbe111}, we then have that $K(u(t_n))=o_n(1)$. This is a contradiction by Lemma \ref{invariant}.

Next, we shall prove that the second assertion holds true. Suppose that there exist $\{\delta_n\} \subset \R^+$ with $\delta_n=o_n(1)$ and $\{t_n\} \subset [0, T_{\max})$ such that
\begin{align} \label{tb1}
-\delta_n \int_{\R^n}|\nabla u(t_n)|^2 \,dx<K(u(t_n))<0.
\end{align}
From \eqref{tbe111}, we see that
\begin{align} \label{tb2}
\begin{split}
\int_{\R^N} |\nabla u(t,x)|^2 \,dx&=\frac{2N(p_2-2)+4b_2}{N(p_2-2)-2(2-b_2)} E(u(t))-\frac{4}{N(p_2-2)-2(2-b_2)}K(u(t))
\\
& \quad -\frac{2N(p_2-p_1)+4(b_2-b_1)}{p_1\left(N(p_2-2)-2(2-b_2)\right)} \int_{\R^N} |x|^{-b_1}|u(t,x)|^{p_1}\,dx \\
& < \frac{2N(p_2-2)+4b_2}{N(p_2-2)-2(2-b_2)} E(u(t))-\frac{4}{N(p_2-2)-2(2-b_2)}K(u(t)).
\end{split}
\end{align}
Combining \eqref{tb1} and \eqref{tb2}, we then have that
\begin{align} \label{tb3}
0<\frac{\left(2N(p_2-2)+4b_2\right)\delta_n}{N(p_2-2)-2(2-b_2)} E(u(t_n))+\left(1-\frac{4\delta_n}{N(p_2-2)-2(2-b_2)}\right) K(u(t_n)).
\end{align}
Using the conservation of energy and the fact that $K(u(t_n))<-\beta_0$ by Lemma \ref{invariant}, we then reach a contradiction from \eqref{tb3} as $n \to \infty$. Thus the proof is completed.
\end{proof}

To establish the blowup for the Cauchy problem  \eqref{equ}, we introduce the localized virial identity as
$$
I_{\psi}(t):=\int_{\R^N} \psi(x) |u(t,x)|^2 \,dx,
$$
where $\psi : \R^N \to \R$ is a smooth cut-off function. The standard virial identity is stated in the following lemma.
\begin{lem} \label{virial}
Let $u \in C([0, T_{max}), H^1(\R^N))$ be the solution to the Cauchy problem \eqref{equ}. Then there holds that
$$
I'_{\psi}(t)=2 \mbox{Im} \int_{\R^N} \left(\nabla \psi \cdot \nabla u \right) \overline{u} \,dx,
$$
\begin{align*}
I''_{\psi}(t)&=4 \textnormal{Re} \sum_{j, k=1}^{N} \int_{\R^N} \partial_{j,k}^2\psi \partial_j u \partial_k \overline{u} \,dx -\int_{\R^N} \Delta^2 \psi |u|^2 \,dx \\
& \quad  +\frac{2(p_1-2)}{p_1} \int_{\R^N} |x|^{-b_1}|u|^{p_1} \Delta \psi \,dx-\frac{4}{p_1} \int_{\R^N} \nabla \left(|x|^{-b_1}\right) \cdot \nabla \psi |u|^{p_1} \,dx\\
& \quad -\frac{2(p_2-2)}{p_2} \int_{\R^N} |x|^{-b_2}|u|^{p_2} \Delta \psi \,dx +\frac{4}{p_2} \int_{\R^N} \nabla \left(|x|^{-b_2}\right) \cdot \nabla \psi |u|^{p_2} \,dx.
\end{align*}
\end{lem}
The proof of the above lemma is a straightforward application of  \cite[Lemma 5.3]{TVZ} for the nonlinearity $\mathcal{N}=|x|^{-b_1} |u|^{p_1-2} u - |x|^{-b_2} |u|^{p_2-2}u$. We shall apply Lemma \ref{virial} with a special choose of function $\psi$. For $R>0$, We introduce a smooth radial function  $\psi=\psi_R$ satisfying 
\begin{align} \label{defpsi}
\psi(r)&=\left\{
\begin{aligned}
&r^2, \quad & 0 \leq r \leq R,\\
& 0, \quad & r \geq 2R,
\end{aligned}
\right. \\
 0&\leq \psi(r) \leq r^2,\; \psi'(r) \leq 2r,\;\psi''(r) \leq 2,\; \psi^{(4)}(r) \leq \frac{4}{R^2}\; \,\, \mbox{for any}\,\,\; r \geq 0.
\end{align}

\begin{lem} \label{vbe}
Let $u \in C([0, T_{max}), H^1(\R^N))$ be the maximal solution to the Cauchy problem \eqref{equ} with $u_0 \in \mathcal{A}^-_{\omega}$. Then there exists $\delta_0>0$ such that $I''_{\psi}(t) \leq -2\delta_0 \|\nabla u(t)\|_2^2$ for any $t \in [0, T_{max})$ under one of the following assumptions,
\begin{itemize}
\item [$(\textnormal{i})$] $u_0 \in H^1(\R^N)$ is radial and $\max\{p_1,p_2\}< 6$;
\item [$(\textnormal{ii})$] $u_0 \in H^1(\R^N)$ and $\max\{p_1, p_2\}< 2+\frac{4}{N}$.
\end{itemize}
\end{lem}
\begin{proof}
Applying Lemma \ref{virial} together with the radial identity
\begin{equation}\label{''}
\partial_j\partial_k=\Big(\frac{\delta_{jk}}r-\frac{x_jx_k}{r^3}\Big)\partial_r+\frac{x_jx_k}{r^2}\partial_r^2,
\end{equation}
we get
\begin{align} \label{id2}
\begin{split}
I''_{\psi}(t)&=4 \int_{\R^N} \frac{\psi'}{r}|\nabla u|^2 \,dx+4 \int_{\R^N} \left(\frac{\psi''}{r^2}-\frac{\psi'}{r^3}\right) |x \cdot \nabla u|^2 \,dx -\int_{\R^N} \Delta^2 \psi |u|^2 \,dx \\
& \quad + \frac{2(p_1-2)}{p_1} \int_{\R^N} \left(\psi'' + \left(N-1\right) \frac{\psi'}{r}\right)|x|^{-b_1}|u|^{p_1} \, dx+\frac{4b_1}{p_1} \int_{\R^N} \frac{\psi'}{r} |x|^{-b_1}|u|^{p_1} \, dx\\
& \quad - \frac{2(p_2-2)}{p_2} \int_{\R^N} \left(\psi'' + \left(N-1\right)\frac{\psi'}{r}\right)|x|^{-b_2}|u|^{p_2} \, dx-\frac{4b_2}{p_2} \int_{\R^N} \frac{\psi'}{r} |x|^{-b_2}|u|^{p_2} \, dx \\
&=8 K(u(t)) +I_1+I_2+I_3+I_4,
\end{split}
\end{align}
where
\begin{align*}
I_1:= -\int_{\R^N} \Delta^2 \psi |u|^2 \,dx,
\end{align*}
\begin{align*} 
I_2:=4 \int_{\R^N} \left(\frac{\psi'}{r}-2 \right)|\nabla u|^2 \,dx+4 \int_{\R^N} \left(\frac{\psi''}{r^2}-\frac{\psi'}{r^3}\right) |x \cdot \nabla u|^2 \,dx,
\end{align*}
\begin{align*} 
I_3:= \frac{2(p_1-2)}{p_1} \int_{\R^N} \left(\psi'' + \left(N-1+\frac{2b_1}{p_1-2}\right) \frac{\psi'}{r}-\frac{2N(p_1-2)+4b_1}{p_1-2}\right)|x|^{-b_1}|u|^{p_1} \, dx,
\end{align*}
\begin{align*} 
I_4:=-\frac{2(p_2-2)}{p_2} \int_{\R^N} \left(\psi'' + \left(N-1+\frac{2b_2}{p_2-2}\right) \frac{\psi'}{r}-\frac{2N(p_2-2)+4b_2}{p_2-2}\right)|x|^{-b_2}|u|^{p_2} \, dx.
\end{align*}
Next we shall deal with the terms $I_1$, $I_2$, $I_3$ and $I_4$. First we see that
$$
|I_1| \leq \frac{C}{R^2}.
$$
Define
$$
\Omega_1:=\left\{x \in \R^N : \frac{\psi''(r)}{r^2}-\frac{\psi'(r)}{r^3}>0\right\}, \quad \Omega_2:=\R^N \backslash \Omega_1.
$$
Since $\psi'(r) \leq 2 r$ for any $r>0$, then we find that
\begin{align*}
I_2 &= 4 \int_{\R^N} \left(\frac{\psi'}{r}-2 \right)|\nabla u|^2 \,dx+4 \int_{\Omega_1} \left(\frac{\psi''}{r^2}-\frac{\psi'}{r^3}\right) |x \cdot \nabla u|^2 \,dx+4 \int_{\Omega_2} \left(\frac{\psi''}{r^2}-\frac{\psi'}{r^3}\right) |x \cdot \nabla u|^2\,dx \\
& \leq  4\int_{\R^N} \left(\frac{\psi'}{r}-2 \right)|\nabla u|^2 \,dx+4 \int_{\Omega_1} \left(\frac{\psi''}{r^2}-\frac{\psi'}{r^3}\right) |x \cdot \nabla u|^2 \,dx\\
& \leq 4 \int_{\Omega_1} \left(\frac{\psi'}{r}-2 \right)|\nabla u|^2 \,dx+4 \int_{\Omega_1} \left(\psi''-\frac{\psi'}{r}\right) |\nabla u|^2 \,dx =4\int_{\Omega_1} \left(\psi''-2 \right)|\nabla u|^2 \,dx \leq 0.
\end{align*}
Observe that
$$
 \int_{|x| \leq R} \left(\psi'' + \left(N-1+\frac{2b_1}{p_1-2}\right) \frac{\psi'}{r}-\frac{2N(p_1-2)+4b_1}{p_1-2}\right)|x|^{-b_1}|u|^{p_1} \, dx=0,
$$
where we used the fact that $\psi'(r)=2r$ and $\psi''(r)=2$ for any $0\leq r \leq R$. Therefore, we derive that
\begin{align*}
I_3&=\frac{2(p_1-2)}{p_1} \int_{|x|>R} \left(\psi'' + \left(N-1+\frac{2b_1}{p_1-2}\right) \frac{\psi'}{r}-\frac{2N(p_1-2)+4b_1}{p_1-2}\right)|x|^{-b_1}|u|^{p_1} \, dx\\
&\leq C \int_{|x|>R} |x|^{-b_1}|u|^{p_1} \,dx.
\end{align*}
Similarly, we can obtain that
$$
I_4 \leq C \int_{|x|>R} |x|^{-b_2}|u|^{p_2} \,dx.
$$
First suppose that $u_0 \in H^1(\R^N)$ is radial. It then follows from the radial Sobolev embedding \eqref{re} that
\begin{eqnarray*}
    I_3 &\leq& \frac{C}{R^{b_1}}\int_{|x|>R} |u(t,x)|^{p_1-2}|u|^2\,dx,\\
    &\leq& \frac{C}{R^{\frac{(N-1)(p_1-2)}{2}+b_1}} \|\nabla u(t)\|_2^{\frac{p_1-2}{2}} \|u(t)\|_2^{\frac{p_1+2}{2}},\\
    &\leq& \frac{C}{R^{\frac{(N-1)(p_1-2)}{2}+b_1}} \|\nabla u(t)\|_2^{\frac{p_1-2}{2}},\\
    I_4 &\leq& \frac{C}{R^{\frac{(N-1)(p_2-2)}{2}+b_2}} \|\nabla u(t)\|_2^{\frac{p_2-2}{2}}.
\end{eqnarray*}

Therefore, coming back to \eqref{id2} and using the elementary observation $a^\theta\leq \theta \,a+1-\theta,\,a\geq 0,\, 0\leq \theta\leq 1$, we obtain that, for any $R>0$ large enough, provided that $\max\{p_1,p_2\}\leq6$,
\begin{align} \label{v111}
\begin{split}
I_{\psi}''(t) & \leq 8K(u(t))+\frac{C}{R^2} + \left(\frac{C}{R^{\frac{(N-1)(p_1-2)}{2}+b_1}} \|\nabla u(t)\|_2^{\frac{p_1-2}{2}} +\frac{C}{R^{\frac{(N-1)(p_2-2)}{2}+b_2}} \|\nabla u(t)\|_2^{\frac{p_2-2}{2}}\right) \\
& \leq 8K(u(t))+\frac{C}{R^2} +  \frac{C}{R^{\frac{(N-1)(p_1-2)}{2}+b_1}}+ \frac{C}{R^{\frac{(N-1)(p_2-2)}{2}+b_2}} \\
& \quad + \frac{C}{R^{\frac{(N-1)(p_1-2)}{2}+b_1}}\|\nabla u(t)\|_2^2+\frac{C}{R^{\frac{(N-1)(p_2-2)}{2}+b_2}}\|\nabla u(t)\|_2^2.
\end{split}
\end{align}
Now suppose that $u_0 \in H^1(\R^N)$. In virtue of \eqref{GN} with $b=0$, we find that
\begin{eqnarray*}
I_3 &\leq& \frac{C}{R^{b_1}}\int_{|x|>R}|u(t,x)|^{p_1}\,dx,\\
&\leq& \frac{C}{R^{b_1}} \|\nabla u(t)\|_2^{\frac{N(p_1-2)}{2}} \|u(t)\|_2^{p_1-\frac{N(p_1-2)}{2}},\\& \leq& \frac{C}{R^{b_1}} \|\nabla u(t)\|_2^{\frac{N(p_1-2)}{2}},\\
I_4 &\leq& \frac{C}{R^{b_2}} \|\nabla u(t)\|_2^{\frac{N(p_2-2)}{2}}.
\end{eqnarray*}

In this case, we similarly get that, for any $R>0$ large enough, provided that $\max\{p_1,p_2\}\leq2+\frac4N$,
\begin{align} \label{v222}
\begin{split}
I_{\psi}''(t) & \leq 8K(u(t))+\frac{C}{R^2} + \left(\frac{C}{R^{b_1}} \|\nabla u(t)\|_2^{\frac{N(p_1-2)}{2}} +\frac{C}{R^{b_2}} \|\nabla u(t)\|_2^{\frac{N(p_2-2)}{2}}\right) \\
& \leq 8K(u(t))+\frac{C}{R^2} + \frac{C}{R^{b_1}}+\frac{C}{R^{b_2}} + \frac{C}{R^{b_1}}\|\nabla u(t)\|_2^2+ \frac{C}{R^{b_2}}\|\nabla u(t)\|_2^2.
\end{split}
\end{align}
At this point, taking into account Lemma \ref{tbelow}, \eqref{v111} and \eqref{v222}, we then have the desired conclusions and the proof is completed.
\end{proof}

\begin{proof}[Proof of Theorem \ref{thm2}]
Suppose that $u$ exists globally in time, i.e. $T_{max}=+\infty$. First we treat the case that $|x|u_0 \in L^2(\R^N)$. Hence, from Lemma \ref{tbelow}, we know that $I''_{|x|^2}(t)=8K(u(t)) \leq -2\eps_0\delta_0$ for any $t \geq 0$. Integrating twice the previous inequality, one gets a contradiction for large time. Next we handle other two cases. Let $\psi : \R^N \to \R$ be defined by \eqref{defpsi} in these two cases.

From Lemmas \ref{tbelow} and \ref{vbe}, we also have that $I''_{\psi}(t) \leq -2\eps_0\delta_0$ for any $t \geq 0$. Now integrating on $[0, t]$ for $t>0$, we then get that $I'_{\psi}(t) \leq -2\eps_0\delta_0 t + I'_{\psi}(0)$. Note that $|I'_{\psi}(0)| \leq C \|\nabla u_0\|_2 \|u_0\|_2 \leq C$ for some $C>0$. Thus there exists $T_0>0$ such that $I'_{\psi}(t)<0$ for any $t>T_0$. It then follows from Lemma \ref{vbe} that
$$
I_{\psi}'(t)=\int_{T_0}^t I''_{\psi}(s)\,ds +I_{\psi}'(T_0) \leq -2\delta_0 \int_{T_0}^t \|\nabla u(s)\|_2^2 \,ds, \quad t \geq T_0.
$$
Observe that
$$
|I_{\psi}'(t)| \leq 4 R \|\nabla u(t)\|_2 \|u(t)\|_2 \leq C \|\nabla u(t)\|_2.
$$
As a consequence, there holds that
\begin{align} \label{b1}
\int_{T_0}^t \|\nabla u(s)\|_2^2 \,ds \leq C |I'_{\psi}(t)| \leq C \|\nabla u(t)\|_2.
\end{align}
Define the real function
$$
f(t):=\int_{T_0}^t \|\nabla u(s)\|_2^2 \,ds.
$$
Then \eqref{b1} indicates that $f^2(t) \leq C f'(t)$ for any $t>T_0$. Taking $T_1>T_0$ and integrating on $[T_1, t]$, we then have that
\begin{align}\label{b2}
\frac{t-T_1}{C} \leq \int_{T_1}^t \frac{f'(s)}{f^2(s)} \,ds=\frac{1}{f(T_1)}-\frac{1}{f(t)} \leq \frac{1}{f(T_1)}.
\end{align}
On the other hand, using Lemma \ref{tbelow}, we know that $f(t) \geq \eps_0^2(t-T_0)$ for any $t>T_0$. In particular, there holds that $f(T_1) \geq \eps_0^2(T_1-T_0)$. Taking $t>0$ large enough in \eqref{b2}, we then reach a contradiction. Thus the proof is completed.
\end{proof}

\begin{proof}[Proof of Theorem \ref{thm3}]
Let $\psi : \R^N \to \R$ be defined by \eqref{defpsi}. First we consider the case that $u_0 \in H^1(\R^N)$ is radial. In this case, using \eqref{v111}, we have that, for any $R>0$ large enough,
$$
I_{\psi}''(t) \leq 8K(u(t))+\frac{C}{R^2} + \frac{C}{R^{\frac{(N-1)(p_1-2)}{2}+b_1}}\|\nabla u(t)\|_2^{\frac{p_1-2}{2}}+\frac{C}{R^{\frac{(N-1)(p_2-2)}{2}+b_2}}\|\nabla u(t)\|_2^{\frac{p_2-2}{2}}.
$$
By Young's inequality for $\eps>0$, we see that
$$
\frac{C}{R^{\frac{(N-1)(p_1-2)}{2}+b_1}}\|\nabla u(t)\|_2^{\frac{p_1-2}{2}} \leq \frac{\eps}{2} \|\nabla u(t)\|_2^2 + \frac{C_{\eps}}{R^{\frac{2(N-1)(p_1-2)+4b_1}{6-p_1}}},
$$
$$
\frac{C}{R^{\frac{(N-1)(p_2-2)}{2}+b_2}}\|\nabla u(t)\|_2^{\frac{p_2-2}{2}} \leq \frac{\eps}{2}\|\nabla u(t)\|_2^2 + \frac{C_{\eps}}{R^{\frac{2(N-1)(p_2-2)+4b_2}{6-p_2}}}.
$$
For simplicity, we shall assume that
$$
\frac{2(N-1)(p_1-2)+4b_1}{6-p_1} \leq \frac{2(N-1)(p_2-2)+4b_2}{6-p_2}.
$$
Since $p_2>2+\frac{2(2-b_2)}{N}$, then
$$
\frac{2(N-1)(p_2-2)+4b_2}{6-p_2}>2.
$$
As a consequence, we know that, for any $R>0$ large enough,
\begin{align*}
I_{\psi}''(t) & \leq 8K(u(t))+ \eps \|\nabla u(t)\|_2^2 +\frac{C}{R^2} + \frac{C}{R^{\frac{2(N-1)(p_2-2)+4b_2}{6-p_2}}} \\
&=\left(4N(p_2-2)+8b_2\right) E(u(t))-\left(2N(p_2-2)-4(2-b_2)\right)\|\nabla u(t)\|^2_2\\
& \quad -\frac{4N(p_2-p_1)+8(b_2-b_1)}{p_1} \int_{\R^N} |x|^{-b_1}|u(t)|^{p_1}\,dx + \eps \|\nabla u(t)\|_2^2 +\frac{C}{R^{\frac{2(N-1)(p_2-2)+4b_2}{6-p_2}}}\\
& \leq \left(4N(p_2-2)+8b_2\right) E(u_0)-\left(N(p_2-2)-2(2-b_2)\right)\|\nabla u(t)\|^2_2 +\frac{C}{R^{\frac{2(N-1)(p_2-2)+4b_2}{6-p_2}}} \\
& \leq -\left(N(p_2-2)-2(2-b_2)\right)\|\nabla u(t)\|^2_2+\frac{C}{R^{\frac{2(N-1)(p_2-2)+4b_2}{6-p_2}}},
\end{align*}
where $\eps>0$ is small enough. Thus we get that, for any $R>0$ large enough,
\begin{align} \label{b1111}
I_{\psi}''(t) +\left(N(p_2-2)-2(2-b_2)\right)\|\nabla u(t)\|^2_2 \leq \frac{C}{R^{\frac{2(N-1)(p_2-2)+4b_2}{6-p_2}}}.
\end{align}
In virtue of the definition of $I_{\psi}$ and Lemma \ref{virial}, we see that
$$
I_{\psi}(t) \leq C R^2 \|u(t)\|_2^2 \leq CR^2,
$$
$$
I'_{\psi}(t) \leq C R \|\nabla u(t)\|_2 \|u(t)\|_2 \leq C R \|\nabla u(t)\|_2.
$$
Then integrating \eqref{b1111} twice on $[t_0, t]$ for $0<t_0<t<T_{max}$ gives that
\begin{align}\label{blow1}
\begin{split}
I_{\psi}(t) + \int_{t_0}^t \int_{t_0}^s \|\nabla u(\tau)\|_2^2 \,d\tau ds &\leq +\frac{C(t-t_0)^2}{R^{\frac{2(N-1)(p_2-2)+4b_2}{p_2-6}}}+I'_{\psi}(t_0)(t-t_0)+I_{\psi}(t_0) \\
& \leq \frac{C(t-t_0)^2}{R^{\frac{2(N-1)(p_2-2)+4b_2}{p_2-6}}}+C R\|\nabla u(t_0)\|_2(t-t_0)+CR^2.
\end{split}
\end{align}
Observe that $I_{\psi} (t) \geq 0$ and
$$
\int_{t_0}^t \int_{t_0}^s \|\nabla u(\tau)\|_2^2 \,d\tau ds=\int_{t_0}^t (t-\tau) \|\nabla u(\tau)\|_2^2 \,d\tau.
$$
Therefore, from \eqref{blow1}, we conclude that
$$
\int_{t_0}^t (t-\tau) \|\nabla u(\tau)\|_2^2 \,d\tau \leq \frac{C(t-t_0)^2}{R^{\frac{2(N-1)(p_2-2)+4b_2}{p_2-6}}}+CR \|\nabla u(t_0)\|_2(t-t_0)+CR^2.
$$
Letting $t_n \nearrow T_{max}$ as $n \to \infty$, we then have that
\begin{align*}
\int_{t_0}^{T_{max}} (T_{max}-\tau) \|\nabla u(\tau)\|_2^2 \,d\tau &\leq \frac{C(T_{max}-t_0)^2}{R^{\frac{2(N-1)(p_2-2)+4b_2}{p_2-6}}} + CR \|\nabla u(t_0)\|_2(T_{max}-t_0)+CR^2.
\end{align*}
Now take $R>0$ such that
$$
\frac{(T_{max}-t_0)^2}{R^{\frac{2(N-1)(p_2-2)+4b_2}{6-p_2}}}=R^2, \quad \mbox{i.e.} \,\, R=(T_{max}-t_0)^{\frac{{6-p_2}}{(N-2)(p_2-2)+2(b_2+2)}}.
$$
As a consequence, there holds that
\begin{align} \label{blow2}
\begin{split}
& \int_{t_0}^{T_{max}} (T_{max}-\tau) \|\nabla u(\tau)\|_2^2 \,d\tau \\
& \leq C (T_{max}-t_0)^{\frac{2(6-p_2)}{(N-2)(p_2-2)+2(b_2+2)}}+C\|\nabla u(t_0)\|_2(T_{max}-t_0)^{\frac{{6-p_2}}{(N-2)(p_2-2)+2(b_2+2)}+1} \\
& \leq C (T_{max}-t_0)^{\frac{2(6-p_2)}{(N-2)(p_2-2)+2(b_2+2)}}+C\|\nabla u(t_0)\|_2^2(T_{max}-t_0)^2,
\end{split}
\end{align}
where we used the fact that
$$
\frac{{6-p_2}}{(N-2)(p_2-2)+2(b_2+2)}<1.
$$
Define
$$
g(t):=\int_t^{T_{max}} (T_{max}-\tau)\|\nabla u(\tau)\|_2^2 \,d\tau.
$$
Hence, by \eqref{blow2}, we have that
$$
g(t) \leq C (T_{max}-t_0)^{\frac{2(6-p_2)}{(N-2)(p_2-2)+2(b_2+2)}}-C(T_{max}-t) g'(t).
$$
It then follows that
$$
\left(\frac{g(t)}{T_{max}-t}\right)'=\frac{1}{(T_{max}-t)^2} \left((T_{max}-t)g'(t)+g(t)\right)\leq \frac{C}{(T_{max}-t)^{\frac{2(N-1)(p_2-2)+4b_2}{(N-2)(p_2-2)+2(b_2+2)}}}.
$$
Integrating on $[0, t]$, we then derive that
$$
\frac{g(t)}{T_{max}-t} \leq \frac{g(0)}{T_{max}} +\frac{C}{(T_{\max}-t)^{\frac{N(p_2-2)-2(2-b_2)}{(N-2)(p_2-2)+2(b_2+2)}}}-\frac{C}{T_{max}^{\frac{N(p_2-2)-2(2-b_2)}{(N-2)(p_2-2)+2(b_2+2)}}}.
$$
Therefore, there holds that, for any $t$ close to $T_{max}$,
$$
g(t) \leq C(T_{max}-t)^{\frac{2(6-p_2)}{(N-2)(p_2-2)+2(b_2+2)}}.
$$
This proves \eqref{r1}. Next we shall verify \eqref{r2}. Observe that, for any $t$ close to $T_{max}$,
\begin{align} \label{b112}
\frac{1}{T_{max}-t} \int_t^{T_{max}} (T_{max}-\tau)\|\nabla u(\tau)\|_2^2 \,d\tau \leq \frac{C}{(T_{\max}-t)^{\frac{N(p_2-2)-2(2-b_2)}{(N-2)(p_2-2)+2(b_2+2)}}}.
\end{align}
Taking $\{T_n\} \subset [0, T_{max})$ such that $T_n \nearrow T_{max}$ as $n \to \infty$ and using the mean value theorem, we know that there exists $t_n \in (T_n , T_{max})$ such that
$$
-(T_{max}-t_n) \|\nabla u(t_n)\|_2^2=g'(t_n)=\frac{g(T_{max})-g(T_n)}{T_{max}-T_n}=-\frac{\displaystyle \int_{T_n}^{T_{max}} (T_{max}-\tau)\|\nabla u(\tau)\|_2^2 \,d\tau}{T_{max}-T_n}.
$$
Using \eqref{b112}, we then see that
$$
(T_{max}-t_n) \|\nabla u(t_n)\|_2^2 \leq \frac{C}{(T_{\max}-T_n)^{\frac{N(p_2-2)-2(2-b_2)}{(N-2)(p_2-2)+2(b_2+2)}}}.
$$
Therefore, we conclude that
$$
\|\nabla u(t_n)\|_2^2 \leq \frac{C}{(T_{\max}-t_n)^{\frac{2(N-1)(p_2-2)+4b_2}{(N-2)(p_2-2)+2(b_2+2)}}}.
$$

Next we consider the case that $u_0 \in H^1(\R^N)$. In this case, using \eqref{v222}, we have that, for any $R>0$ large enough,
$$
I_{\psi}''(t) \leq 8K(u)+\frac{C}{R^2} + \frac{C}{R^{b_1}}\|\nabla u(t)\|_2^{\frac{N(p_1-2)}{2}}+\frac{C}{R^{b_2}}\|\nabla u(t)\|_2^{\frac{N(p_2-2)}{2}}.
$$
By Young's inequality for $\eps>0$ small enough, we see that
$$
\frac{C}{R^{b_1}}\|\nabla u(t)\|_2^{\frac{N(p_1-2)}{2}} \leq \frac{\eps}{2}\|\nabla u(t)\|_2^2+\frac{C_{\eps}}{R^{\frac{4b_1}{4-N(p_1-2)}}},
$$
$$
\frac{C}{R^{b_2}}\|\nabla u(t)\|_2^{\frac{N(p_2-2)}{2}} \leq \frac{\eps}{2}\|\nabla u(t)\|_2^2+\frac{C_{\eps}}{R^{\frac{4b_2}{4-N(p_2-2)}}}.
$$
For simplicity, we shall assume that
$$
\frac{4b_1}{4-N(p_1-2)} \leq \frac{4b_2}{4-N(p_2-2)}.
$$
Since $p_2>2+\frac{2(2-b_2)}{N}$, then
$$
\frac{4b_2}{4-N(p_2-2)}>2.
$$
As a result, we derive that, for {any $R>0$ large enough},
\begin{align*}
I_{\psi}''(t) &\leq 8 K(u(t)) +\frac{C}{R^2} + \eps \|\nabla u(t)\|^2_2+\frac{C_{\eps}}{R^{\frac{4b_2}{4-N(p_2-2)}}}\\
&=\left(4N(p_2-2)+8b_2\right) E(u(t))-\left(2N(p_2-2)-4(2-b_2)\right)\|\nabla u(t)\|^2_2\\
& \quad -\frac{4N(p_2-p_1)+8(b_2-b_1)}{p_1} \int_{\R^N} |x|^{-b_1}|u|^{p_1}\,dx+ \eps \|\nabla u(t)\|^2_2+\frac{C}{R^{\frac{4b_2}{4-N(p_2-2)}}} \\
& \leq -\left(N(p_2-2)-2(2-b_2)\right)\|\nabla u(t)\|^2_2 +\frac{C}{R^{\frac{4b_2}{4-N(p_2-2)}}},
\end{align*}
where $\eps>0$ is small enough. Hence we have that, for {any $R>0$ large enough},
\begin{align} \label{b111}
I_{\psi}''(t) +\left(N(p_2-2)-2(2-b_2)\right)\|\nabla u(t)\|^2_2 \leq \frac{C}{R^{\frac{4b_2}{4-N(p_2-2)}}}.
\end{align}
At this point, proceeding as before, we are able to get the desired conclusions \eqref{r11} and \eqref{r22}. This completes the proof.
\end{proof}

\begin{proof} [Proof of Corollary \ref{instability}]
Let us first define a cut-off function $\chi_R: \R^N \to \R$ by
$$
\chi_R(x):=\chi\left(\frac{|x|}{R}\right), \quad R>0,
$$
where $\chi \in C^{\infty}_0(\R^N, \R)$ is a function such that $\chi(r) \geq 0$ for $r \geq 0$, $\chi(r) =1$ for $0 \leq r \leq 1$ and $\chi(r)=0$ for $r \geq 2$. Note that $S_{\omega}(u)=m_{\omega}$ and $K(u)=0$. It follows from Lemma \ref{maximum} that there exists $\lambda>1$ such that $S_{\omega}(u_{\lambda})<m_{\omega}$ and $K(u_{\lambda})<0$. This infers that $u_{\lambda} \in \mathcal{A}^-_{\omega}$. It is simple to check that $\|u_{\lambda}-\chi_R u_{\lambda}\| \to 0$ as $R \to \infty$. Then we know that $\chi_R u_{\lambda} \in \mathcal{A}^-_{\omega}$ and $|x| \chi_R u_{\lambda} \in L^2(\R^N)$ for any $R>0$ large enough. As an application of Theorem \ref{thm2}, we have that the desired conclusion. By replacing the role of $H^1(\R^N)$ by $X$, we can also have the second assertion. Thus the proof is completed.
\end{proof}

\section{Scattering} \label{sca}

The aim of this section is to discuss scattering of solutions to the Cauchy problem \eqref{equ} with initial data belonging to $\mathcal{A}_{\omega}^+$ and give the proof of Theorem \ref{scattering}. For simplicity, we shall always assume that \eqref{H} holds. To begin with, we shall present some basic facts.

\begin{defi}
Let $N \geq 3$ and $0 \leq s <1$. A pair of real numbers $(q, r)$ is called $s$-admissible if
$$
2 \leq q, r \leq \infty, \quad \frac{2N}{N-2s} \leq r<\frac{2N}{N-2}, \quad \frac{2}{q}+\frac{N}{r}=\frac{N}{2}-s.
$$
\end{defi}

Denote the set of $s$-admissible pairs by $\Lambda_s$, i.e. $\Lambda_s:=\Big\{(q, r) : (q, r) \,\, \mbox{is} \,\, s\mbox{-admissible} \Big\}$. Define
$$
\|u\|_{\Lambda^s(I)}:=\sup_{(q ,r) \in \Lambda_s} \|u\|_{L^q(I, L^r)}.
$$
Likewise, we define
$$
\Lambda_{-s}:=\Big\{(q, r) : (q, r) \,\, \mbox{is} \,\, (-s)\mbox{-admissible} \Big\}, \quad \|u\|_{\Lambda^{-s'}(I)}:=\inf_{(q, r) \in \Lambda_{-s}}\|u\|_{L^{q'}(I, L^{r'})},
$$
where $(q', r')$ is the conjugate exponent pair of $(q ,r)$. When $s=0$, then we shall denote $\Lambda^s$ by $\Lambda$ and $\Lambda^{-s'}$ by $\Lambda'$.

\begin{lem} \label{str} (\cite{Fo, KT})
Let $N \geq 3$, $0 \leq s<1$ and a time slab $I \subset \R$. Then there holds that
\begin{itemize}
\item[$(\textnormal{i})$] $\displaystyle \|e^{\textnormal{i} t\Delta} f\|_{\Lambda^s(I)} \lesssim \|f\|_{\dot{H}^s}$;
\item[$(\textnormal{ii})$] $\displaystyle \left\|\int_0^t e^{\textnormal{i} (t-s)\Delta} g(\cdot, s) \,ds\right\|_{\Lambda^s(I)}+ \left\|\int_{\R} e^{\textnormal{i} (t-s)\Delta} g(\cdot, s) \,ds\right\|_{\Lambda^s(I)} \leq \|g\|_{\Lambda^{-s'}(I)}$;
\item[$(\textnormal{iii})$] $\displaystyle \left\|\int_I e^{\textnormal{i} (t-s)\Delta} g(\cdot, s) \,ds\right\|_{\Lambda^s(\R)} \leq \|g\|_{\Lambda^{-s'}(I)}$.
\end{itemize}
\end{lem}

Let $\chi \in C^{\infty}_0(\R^N, \R)$ be a smooth radial function such that $0 \leq \chi(r) \leq 1$ for $r \geq 0$, $\chi(r)=1$ for  $0 \leq r \leq \frac 12$ and $\chi(r)=0$ for $r \geq 1$. For $R>0$, we define
\begin{align} \label{chi-R}
\chi_R(x):=\chi \left(\frac{x}{R}\right), \quad x \in \R^N.
\end{align}
Clearly, there holds that $ \chi \leq \chi_R$ for $R>2$.

\begin{lem} Let $N\geq 3$, $0<s<1$, $0<b<\min\left\{2, {N}/{2}\right\}$, $2+\frac{4-2b}{N}<p<2_b^*$ and $I \subset \R$ be a bounded time slab. Then there exist $0<\theta_1, \theta_2<1$ and $\frac{pN}{N-b}< r< \frac{2N}{N-2}$  independent of $R$ such that
\label{NLEst}
\begin{align*} 
\left\||x|^{-b} |u|^{p-2}u\right\|_{\Lambda^{-s}(I)}
&\lesssim |I|^{\theta_1}\|u\|_{L^\infty(I, L^r)}^{p-2}\|\chi_R\,u\|_{L^\infty(I, L^r)}+|I|^{\theta_2}\|u\|_{L^\infty(I, L^p)}^{p-2}\|\chi_R\,u\|_{L^\infty(I, L^p)}\\
&\quad +R^{-b}|I|^{\theta_2}\|u\|_{L^\infty(I, L^p)}^{p-1},
\end{align*}
where $\chi_R$ is defined by \eqref{chi-R} for $R>2$.
\end{lem}
\begin{proof}
Denote by $\mathcal{N}=|x|^{-b} |u|^{p-2}u$ and decompose $\mathcal{N}$ as follows
\begin{align*}
\mathcal{N}&=\chi\,\mathcal{N}+(\chi_R-\chi)\mathcal{N}+(1-\chi_R)\mathcal{N}\\
&:= (I)+(II)+(III).
\end{align*}
To estimate the first term $(I)$, we choose $(q,r)\in \Lambda_{-s}$ and $\gamma> 1$ such that
\begin{equation*}
\frac{pN}{N-b}< r< \frac{2N}{N-2}, \quad \frac{1}{\gamma}=1-\frac{p}{r}> \frac{b}{N}.
\end{equation*}
This is possible because $2< \frac{pN}{N-b}< \frac{2N}{N-2}$. Applying H\"older's inequality yields that
\begin{align*}
\|(I)\|_{L^{q'}(I, L^{r'})}&\leq\left\|\left\||x|^{-b}\right\|_{L^\gamma(|x| \leq 1)}\|u\|_{L^r}^{p-2}\|\chi\,u\|_{L^r}\right\|_{L^{q'}(I)}\\
&\lesssim |I|^{\frac{1}{q'}}\|u\|_{L^\infty(I, L^r)}^{p-2}\|\chi_R\,u\|_{L^\infty(I, L^r)},
\end{align*}
where we used the fact that $\chi \leq \chi_R$ for $R>2$ in the second inequality. The term $(II)$ can be estimated in an easier way. Indeed, using the fact that $|x|^{-b} \lesssim 1$ on $\mbox{supp}\,(\chi_R-\chi)$ and choosing $d>1$ such that $(d, p)\in \Lambda_{-s}$, we then obtain that
\begin{align*}
\|(II)\|_{L^{d'}(I, L^{p'})}& \lesssim \left\|\|u\|_{L^p}^{p-2}\|(\chi_R-\chi)u\|_{L^p}\right\|_{L^{d'}(I)}\\
&\lesssim |I|^{\frac{1}{d'}}\|u\|_{L^\infty(I, L^p)}^{p-2}\|\chi_R\,u\|_{L^\infty(I, L^p)}.
\end{align*}
Arguing similarly as for $(II)$, we infer that
\begin{align*}
\|(III)\|_{L^{d'}(I, L^{p'})} &\lesssim R^{-b}\||u|^{p-2}u\|_{L^{d'}(I, L^{p'})} \\
&\leq R^{-b}|I|^{\frac{1}{d'}}\|u\|_{L^\infty(I, L^p)}^{p-1}.
\end{align*}
This finishes the proof of Lemma \ref{NLEst} with $\theta_1=\frac{1}{q'}$ and $\theta_2=\frac{1}{d'}$.
\end{proof}

\begin{lem} \label{global1}
Let $u \in C([0, T_{max}), H^1(\R^N))$ be the maximal solution to the Cauchy problem \eqref{equ} with $u_0 \in \mathcal{A}^+_{\omega}$. Then $u(t) \in \mathcal{A}^+_{\omega}$ for any $t \in [0, T_{\max})$. Moreover, $u(t)$ exists globally in time, i.e. $T_{max}=+\infty$.
\end{lem}
\begin{proof}
Arguing as the proof of Lemma \ref{invariant}, we can show that $\mathcal{A}_{\omega}^+$ is invariant under the flow of the Cauchy problem \eqref{equ}. Let us now prove the global existence of the solution. Observe that
\begin{align*} 
m_{\omega}>S_{\omega}(u(t))-\frac{2}{N(p_1-2)+2b_1}K(u(t)) &=\left(\frac 12-\frac{2}{N(p_1-2)+2b_1}\right)\|\nabla u(t)\|^2_2 + \frac {\omega}{2} \|u(t)\|_2^2\\
& \quad +\frac{1}{p_2} \left(\frac{N(p_2-p_1)+2(b_2-b_1)}{N(p_1-2)+2b_1}\right)\int_{\R^N} |x|^{-b_2}|u(t)|^{p_2}\,dx \\
& \gtrsim \|\nabla u(t)\|^2_2.
\end{align*}
Hence $T_{max}=+\infty$ and the proof is completed.
\end{proof}

\begin{lem} \label{vach}
There holds that
$$
m_{\omega}=\inf_{u \in \mathcal{P}}I(u),
$$
where
$$
I_{\omega}(u):=S_{\omega}(u)-\frac{2}{N(p_1-2)+2b_1}K(u), \quad \mathcal{P}:=\Big\{ u \in H^1(\R^N) \backslash \{0\}: K(u) \leq 0\Big\}.
$$
\end{lem}
\begin{proof}
It is clear to see that $\displaystyle\inf_{u \in \mathcal{P}}I(u) \leq m_{\omega}$. For any $u \in \mathcal{P}$, by Lemma \ref{maximum}, there exists a unique $0<t_u \leq 1$ such that $K(u_{t_u})=0$. Therefore, we have that
\begin{align*}
m_{\omega} \leq I_{\omega}(u_{t_u})&=S_{\omega}(u_{t_u})-\frac{2}{N(p_1-2)+2b_1}K(u_{t_u}) \\
&=t_u^2\left(\frac 12-\frac{2}{N(p_1-2)+2b_1}\right)\|\nabla u\|_2^2 \,dx+\frac 12 \|u\|_2^2 \\
& \quad +\frac{t_u^{\frac{N(p_2-2)+2b_2}{2}}}{p_2} \left(\frac{N(p_2-p_1)+2(b_2-b_1)}{N(p_1-2)+2b_1}\right)\int_{\R^N} |x|^{-b_2}|u|^{p_2}\,dx \\
& \leq S_{\omega}(u)-\frac{2}{N(p_1-2)+2b_1}K(u)=I_{\omega}(u).
\end{align*}
This implies that $m_{\omega} \leq \displaystyle\inf_{u \in \mathcal{P}}I(u) $. Then the proof is completed.
\end{proof}

\begin{lem} \label{coerc}
Let $u \in C([0, +\infty), H^1(\R^N))$ be the global solution to the Cauchy problem \eqref{equ} with $u_0 \in \mathcal{A}^+_{\omega}$. Then, for any $R>>1$ and $t\geq 0$, there holds that
\begin{equation*}
K(\chi_R u(t))\gtrsim\|\nabla(\chi_R u(t))\|^2_2,
\end{equation*}
where the function $\chi_R$ is given by \eqref{chi-R}.
\end{lem}
\begin{proof}
First observe that
\begin{align*}
I_{\omega}(u(t))&=\left(\frac 12-\frac{2}{N(p_1-2)+2b_1}\right)\|\nabla u(t)\|^2_2 + \frac {\omega}{2} \|u(t)\|_2^2 \\
& \quad + \frac{1}{p_2} \left(\frac{N(p_2-p_1)+2(b_2-b_1)}{N(p_1-2)+2b_1}\right)\int_{\R^N} |x|^{-b_2}|u(t)|^{p_2}\,dx \geq 0.
\end{align*}
It follows that
\begin{align} \label{c111}
\begin{split}
I_{\omega}(\chi_R u(t))&=\left(\frac 12-\frac{2}{N(p_1-2)+2b_1}\right)\|\nabla(\chi_R u(t))\|^2_2 + \frac {\omega}{2} \|\chi_R u(t)\|_2^2 \\
&\quad +\frac{1}{p_2} \left(\frac{N(p_2-p_1)+2(b_2-b_1)}{N(p_1-2)+2b_1}\right)\int_{\R^N} |x|^{-b_2}|\chi_R u(t)|^{p_2}\,dx.
\end{split}
\end{align}
It is clear to see from the definition of $\chi_R$ that
$$
\int_{\R^N}\chi_R^2|\nabla u(t)|^2\,dx=\int_{\R^N}\left(|\nabla(\chi_R u(t))|^2+\chi_R\Delta\chi_R|u(t)|^2\right)\,dx=\int_{\R^N} |\nabla(\chi_R u(t))|^2 \,dx + O(R^{-2}),
$$
$$
\int_{\R^N}\chi_R^2|\nabla u(t)|^2\,dx \leq \int_{\R^N} |\nabla u(t)|^2\,dx, \quad  \int_{\R^N} |\chi_R u(t)|^2 \,dx \leq \int_{\R^N} |u(t)|^2 \,dx,
$$
$$
\int_{\R^N} |x|^{-b_2}|\chi_R u(t)|^{p_2}\,dx \leq \int_{\R^N} |x|^{-b_2}|u(t)|^{p_2}\,dx.
$$
Using \eqref{c111}, we then conclude that
\begin{align} \label{c112}
I_{\omega}(\chi_R u(t)) \leq I_{\omega}(u(t)) +O(R^{-2}).
\end{align}
Note that $\mathcal{A}_{\omega}^+$ is invariant under the flow of the Cauchy problem \eqref{equ} by Lemma \ref{global1}. Thereby we know that $S_{\omega}(u(t))<m_{\omega}$ and $K(u(t))>0$ for any $t \geq 0$. In view of the conservation laws, then there exists $\eps>0$ such that $S_{\omega}(u(t)) \leq m_{\omega}-2 \eps$. This in turn leads to $I_{\omega}(u(t)) \leq m_{\omega}-2 \eps.$ Further, applying \eqref{c112}, we then have that
\begin{align} \label{c110}
I_{\omega}(\chi_R u(t)) \leq m_{\omega}- \eps, \quad R>>1.
\end{align}
It then follows that $K(\chi_R u(t)) >0$. Otherwise, there holds that $K(\chi_R u(t_0)) \leq 0$ for some $t_0>0$.
From Lemma \ref{vach}, then $I_{\omega}(\chi_R u(t_0)) \geq m_{\omega}$. This is impossible by \eqref{c110}.
For simplicity, we shall write $u=u(t)$. It is simple to compute that
\begin{align*}
\frac{dS_{\omega}((\chi_R u)_{\lambda})}{d \lambda}&=\lambda\int_{\R^N} |\nabla (\chi_R u)|^2 \,dx+\frac{N(p_1-2)+2b_1}{2p_1}\lambda^{\frac N 2 (p_1-2)+b_1-1}\int_{\R^N}|x|^{-b_1}|\chi_R u|^{p_1} \, dx\\
&\quad-\frac{N(p_2-2)+2b_2}{2p_2}\lambda^{\frac N 2 (p_2-2)+b_2-1}\int_{\R^N}|x|^{-b_2}|\chi_R u|^{p_2} \, dx \\
&=\frac {1} {\lambda} K(\chi_R u).
\end{align*}
In addition, we observe that
\begin{align} \label{c113} \nonumber
\frac{dK((\chi_R u)_{\lambda})}{d \lambda}&=\lambda\int_{\R^N} |\nabla (\chi_R u)|^2 \,dx+\frac{(N(p_1-2)+2b_1)^2}{4p_1}\lambda^{\frac{N}{2}(p_1-2)+b_1-1}\int_{\R^N}|x|^{-b_1}|\chi_R u|^{p_1} \, dx \\ \nonumber
& \quad -\frac{(N(p_2-2)+2b_2)^2}{4p_2}\lambda^{\frac{N}{2}(p_2-2)+b_2-1}\int_{\R^N}|x|^{-b_2}|\chi_R u|^{p_2} \,dx\\ \nonumber
& =-\frac{dS_{\omega}((\chi_R u)_{\lambda})}{d \lambda} + \lambda \left(2\int_{\R^N} |\nabla (\chi_R u)|^2 \,dx\right. \\ \nonumber
& \quad \left.+\frac{N(p_1-2)+2b_1}{2p_1}\left(1 +\frac{N(p_1-2)+2b_1}{2}\right)\lambda^{\frac{N}{2}(p_1-2)+b_1-2}\int_{\R^N}|x|^{-b_1}|\chi_R u|^{p_1} \, dx \right. \\ \nonumber
& \quad \left.-\frac{N(p_2-2)+2b_2}{2p_2}\left(1 +\frac{N(p_2-2)+2b_2}{2}\right)\lambda^{\frac{N}{2}(p_2-2)+b_2-2}\int_{\R^N}|x|^{-b_2}|\chi_R u|^{p_2} \, dx \right) \\
&:=-\frac{dS_{\omega}((\chi_R u)_{\lambda})}{d \lambda} +\lambda g(\lambda).
\end{align}
To proceed the proof, we shall consider two cases for $g(1)<0$ and $g(1)\geq 0$. Let us first assume that $g(1)<0$. This means that
\begin{align*}
&2\int_{\R^N} |\nabla (\chi_R u)|^2 \,dx+\frac{\left(N(p_1-2)+2b_1\right)\left(N(p_1-2)+2b_1+2\right)}{4p_1}\int_{\R^N}|x|^{-b_1}|\chi_R u|^{p_1} \, dx \\
& < \frac{\left(N(p_2-2)+2b_2\right)\left(N(p_2-2)+2b_2+2\right)}{4p_2}\int_{\R^N}|x|^{-b_2}|\chi_R u|^{p_2} \, dx.
\end{align*}
Direct computations give that, for any $\lambda>1$,
\begin{align*}
&\frac{g'(\lambda)}{\lambda^{\frac{N}{2}(p_1-2)+b_1-3}} \\
&< \left(\frac{\left(N(p_1-2)+2b_1\right)\left(\left(N(p_1-2)+2b_1+2\right)\left(N(p_1-2)+2b_1-4\right)\right)}{8p_1}\int_{\R^N}|x|^{-b_1}|\chi_R u|^{p_1} \, dx \right. \\
& \quad -\left.\frac{\left(N(p_2-2)+2b_2\right)\left(N(p_2-2)+2b_2+2\right)\left(N(p_2-2)+2b_2-4\right)}{8p_2} \int_{\R^N}|x|^{-b_2}|\chi_R u|^{p_2} \, dx \right) \\
&<\left(-\left(N(p_1-2)+2b_1-4\right)\int_{\R^N} |\nabla (\chi_R u)|^2 \,dx \right.\\
& \quad -\left.\frac{\left(N(p_2-2)+2b_2\right)\left(N(p_2-2)+2b_2+2\right)(N(p_2-p_2)+2(b_2-b_1))}{8p_2}\int_{\R^N}|x|^{-b_2}|\chi_R u|^{p_2} \, dx \right) \\
&<0.
\end{align*}
Therefore, we get that $g(\lambda)<g(1)<0$ for any $\lambda>1$. It then follows from \eqref{c113} that
\begin{align} \label{c114}
\frac{dK((\chi_R u)_{\lambda})}{d \lambda}<-\frac{dS_{\omega}((\chi_R u)_{\lambda})}{d \lambda}, \quad \forall \,\lambda>1.
\end{align}
Since $K(\chi_R u)>0$, by Lemma \ref{maximum}, then there exists $\widetilde{\lambda}>1$ such that $K((\chi_R u)_{\widetilde{\lambda}})=0$. Integrating \eqref{c114} on $[1, \widetilde{\lambda}]$, using Lemma \ref{vach} and \eqref{c110}, we then derive that
\begin{align*}
K(\chi_R u) &\geq K((\chi_R u)_{\widetilde{\lambda}}) + \left(S_{\omega}((\chi_R u)_{\widetilde{\lambda}})-S_{\omega}(\chi_R u)\right) \\
&=I_{\omega}((\chi_R u)_{\widetilde{\lambda}}) -I_{\omega}(\chi_R u) -\frac{2}{N(p_1-2)+2b_1} K(\chi_R u) \\
& \geq \eps-\frac{2}{N(p_1-2)+2b_1} K(\chi_R u).
\end{align*}
Therefore $K(\chi_R u) \gtrsim \eps$. On the other hand, since $K(\chi_R u)>0$, then $\|\nabla (\chi_R u)\|_2^2 \leq C_0$ for some $C_0>0$. Accordingly, there holds that $K(\chi_R u) \gtrsim \eps \|\nabla(\chi_R u)\|_2^2$. This proves \eqref{coerc} under the assumption that $g(1)<0$. Next we shall assume that $g(1) \geq 0$. This means that
\begin{align*}
&2\int_{\R^N} |\nabla (\chi_R u)|^2 \,dx+\frac{\left(N(p_1-2)+2b_1\right)\left(N(p_1-2)+2b_1+2\right)}{4p_1}\int_{\R^N}|x|^{-b_1}|\chi_R u|^{p_1} \, dx \\
& \geq \frac{\left(N(p_2-2)+2b_2\right)\left(N(p_2-2)+2b_2+2\right)}{4p_2}\int_{\R^N}|x|^{-b_2}|\chi_R u|^{p_2} \, dx.
\end{align*}
Therefore, we obtain that
\begin{align*}
K(\chi_R u) &\geq \left(1-\frac{4}{N(p_2-2)+2b_2+2}\right)\int_{\R^N} |\nabla (\chi_R u)|^2 \,dx  \\
& \quad + \frac{N(p_1-2)+2b_1}{2p_1} \left(1-\frac{N(p_1-2)+2b_1+2}{N(p_2-2)+2b_2+2}\right)\int_{\R^N}|x|^{-b_1}|\chi_R u|^{p_1} \, dx \\
& \gtrsim \|\nabla (\chi_R u)\|^2_2.
\end{align*}
Thus the proof is completed.
\end{proof}

\begin{lem}
Let $u \in C([0, +\infty), H^1(\R^N))$ be the global solution to the Cauchy problem \eqref{equ} with $u_0 \in \mathcal{A}^+_{\omega}$. Then there exists $R>>1$ such that, for any $T>0$,
\begin{equation}\label{mor}
\frac1T\int_0^T\int_{|x| \leq R}|u|^{\frac{2N}{N-2}}\,dx\,dt\lesssim \frac{R}T+R^{-\min\{b_1,b_2\}}.
\end{equation}
In particular, there exist $\{t_n\} \subset \R^+$ and $\{R_n\} \subset \R^+$ with $t_n \to \infty$ and $R_n \to \infty$ as $n \to \infty$ such that
\begin{align}\label{2}
\int_{|x| \leq R_n} |u(t_n)|^{\frac{2N}{N-2}} \,dx  \to 0 \quad\mbox{as} \,\, n\to\infty.
\end{align}
\end{lem}
\begin{proof}
Let us first present some notations in the spirit of \cite{DM1}. For $R>>1$, we define a smooth radial function $\zeta: \R^N \to \R$ by
$$
\zeta(x)=\left\{
\begin{aligned}
&\frac 12 |x|^2, &|x|\leq \frac R2, \\
&R|x|, &|x|>R.
\end{aligned}
\right.
$$
Moreover, we assume that in the centered annulus $\frac R 2 <|x| \leq R$,
$$
\partial_r\zeta > 0,\quad\partial^2_r\zeta \geq 0, \quad |\partial^\omega \zeta(x)| \leq C_\omega R|x|^{1-\omega},\quad\forall\,\; |\omega| \geq 1,
$$
where $\partial_r \zeta(x)= \nabla \zeta(x) \cdot \frac{x}{|x|}$ denotes the radial derivative. Note that on the centered ball $|x| \leq \frac{R}{2}$, there holds that
\begin{align} \label{z10}
\partial_{jk}^2 \zeta=\delta_{jk},\quad \Delta \zeta =N, \quad \Delta^2 \zeta = 0.
\end{align}
Moreover, for $|x| > R$,
\begin{align} \label{z11}
\partial_{jk}^2\zeta =\frac R{|x|}\Big(\delta_{jk}-\frac{x_jx_k}{|x|^2}\Big),\quad \Delta \zeta=\frac{(N -1)R}{|x|}, \quad \Delta^2 \zeta = 0.
\end{align}
Define
$$
I_{\zeta}(t):=\int_{\R^N} \zeta(x) |u(t)|^2 \,dx.
$$
It immediately follows from Lemma \ref{virial} and Cauchy-Schwarz's inequality that
\begin{align} \label{z1}
\left|I_{\zeta}'(t)\right|=2 \left|\textnormal{Im} \int_{\R^N} (\nabla\zeta\cdot\nabla u) \overline{u}\,dx\right| \lesssim R.
\end{align}
Furthermore, there holds that
\begin{align} \label{zv}
\begin{split}
I''_{\zeta}(t)&=4 \textnormal{Re} \sum_{j, k=1}^{N} \int_{\R^N} \partial_{j,k}^2\zeta \partial_j u \partial_k \overline{u} \,dx -\int_{\R^N} \Delta^2 \zeta |u|^2 \,dx \\
& \quad  +\frac{2(p_1-2)}{p_1} \int_{\R^N} |x|^{-b_1}|u|^{p_1} \Delta \zeta \,dx+\frac{4b_1}{p_1} \int_{\R^N} \frac{x \cdot \nabla \zeta}{|x|^2} |x|^{-b_1}|u|^{p_1} \,dx\\
& \quad -\frac{2(p_2-2)}{p_2} \int_{\R^N} |x|^{-b_2}|u|^{p_2} \Delta \zeta \,dx -\frac{4b_2}{p_2} \int_{\R^N} \frac{x \cdot \nabla \zeta}{|x|^2} |x|^{-b_2}|u|^{p_2} \,dx.
\end{split}
\end{align}
Observe from \eqref{z10} that
$$
\textnormal{Re} \sum_{j, k=1}^{N} \int_{|x| \leq \frac R 2} \partial_{j,k}^2\zeta \partial_j u \partial_k \overline{u} \,dx= \int_{|x| \leq \frac R 2} |\nabla u|^2 \,dx.
$$
Taking into account of the identity
\begin{equation*}
\partial_{jk}^2=\left(\frac{\delta_{jk}}r-\frac{x_jx_k}{r^3}\right)\partial_r+\frac{x_lx_k}{r^2}\partial_r^2,
\end{equation*}
we then obtain that
\begin{align*}
\textnormal{Re} \sum_{j, k=1}^{N} \int_{|x| > \frac R 2} \partial_{j,k}^2\zeta \partial_j u \partial_k \overline{u} \,dx
&=\textnormal{Re}\sum_{j, k=1}^{N} \int_{|x| > \frac R 2} \left(\left(\frac{\delta_{jk}}r-\frac{x_jx_k}{r^3}\right)\partial_r\zeta +\frac{x_jx_k}{r^2}\partial_r^2\zeta \right) \partial_j u \partial_k \overline{u}\,dx\\
&=\int_{|x| > \frac R 2}\left(|\nabla u|^2-\frac{|x\cdot\nabla u|^2}{|x|^2}\right)\frac{\partial_r\zeta }{|x|}\,dx
+\int_{|x| > \frac R 2}\frac{|x\cdot\nabla u|^2}{|x|^2}\partial_r^2\zeta \,dx.
\end{align*}
Therefore, from \eqref{z10}, \eqref{z11} and \eqref{zv}, we conclude that
\begin{align} \label{zv1} \nonumber
I''_{\zeta}(t)&=4\int_{|x| \leq \frac R 2} |\nabla u|^2 \,dx +\frac{2(p_1-2)N+4b_1}{p_1} \int_{|x| \leq \frac R 2} |x|^{-b_1}|u|^{p_1} \,dx \\ \nonumber
&\quad -\frac{2(p_2-2)N+4b_2}{p_2} \int_{|x| \leq \frac R 2} |x|^{-b_2}|u|^{p_2} \,dx-\int_{\frac R 2<|x| \leq R}\Delta^2 \zeta |u|^2 \,dx \\ \nonumber
& \quad + \frac{2(p_1-2)N+4b_1}{p_1} \int_{|x| > \frac R 2} |x|^{-b_1}|u|^{p_1} \Delta \zeta \,dx-\frac{2(p_2-2)N+4b_2}{p_2} \int_{|x| >\frac R 2} |x|^{-b_2}|u|^{p_2} \Delta \zeta \,dx \\
& \quad +4\int_{|x| > \frac R 2}\left(|\nabla u|^2-\frac{|x\cdot\nabla u|^2}{|x|^2}\right)\frac{\partial_r\zeta }{|x|}\,dx
+4\int_{|x| > \frac R 2}\frac{|x\cdot\nabla u|^2}{|x|^2}\partial_r^2\zeta \,dx.
\end{align}
Note that
$$
\int_{\frac R 2<|x| \leq R}\Delta^2 \zeta |u|^2 \,dx \lesssim R^{-2},
$$
$$
\int_{|x| >\frac R 2} |x|^{-b_1}|u|^{p_1} \Delta \zeta \,dx \lesssim R^{-b_1}, \quad \int_{|x| >\frac R 2} |x|^{-b_2}|u|^{p_2} \Delta \zeta \,dx \lesssim R^{-b_2},
$$
$$
\int_{|x| > \frac R 2}\left(|\nabla u|^2-\frac{|x\cdot\nabla u|^2}{|x|^2}\right)\frac{\partial_r\zeta }{|x|}\,dx
+\int_{|x| > \frac R 2}\frac{|x\cdot\nabla u|^2}{|x|^2}\partial_r^2\zeta \,dx \geq 0.
$$
Hence, by \eqref{zv1}, we get that
\begin{align*}
I''_{\zeta}(t)& \geq 4\int_{|x| \leq \frac R 2} |\nabla u|^2 \,dx +\frac{2(p_1-2)N+4b_1}{p_1} \int_{|x| \leq \frac R 2} |x|^{-b_1}|u|^{p_1} \,dx \\
&\quad -\frac{2(p_2-2)N+4b_2}{p_2} \int_{|x| \leq \frac R 2} |x|^{-b_2}|u|^{p_2} \,dx-R^{-\min \{b_1, b_2\}},
\end{align*}
where we also used the assumption that $b_1, b_2<2$. This along with Lemma \ref{coerc} then implies that
$$
I''_{\zeta}(t) \geq 4 K(\chi_R u)-R^{-\min \{b_1, b_2\}} +o_R(1) \gtrsim 4\|\nabla (\chi_R u)\|_2^2 -R^{-\min \{b_1, b_2\}}+o_R(1),
$$
where $o_R(1)$ denotes an infinitesimal quantity as $R \to \infty$. This clearly leads to
\begin{align} \label{i11}
I''_{\zeta}(t)  \gtrsim \|\nabla (\chi_R u)\|_2^2 -R^{-\min \{b_1, b_2\}}, \quad R>>1.
\end{align}
Integrating \eqref{i11} with respect to $t$ on $[0, T]$, using \eqref{z1} and Sobolve's embedding inequality gives that
$$
\frac 1 T \int_0^T\int_{\R^N}|\chi_Ru|^{\frac{2N}{N-2}}\,dx\,dt \lesssim \frac{R}{T} +R^{-\min \{b_1, b_2\}}.
$$
Since
$$
\int_{\R^N}|\chi_Ru|^{\frac{2N}{N-2}}\,dx \geq C \int_{|x|<R}|u|^{\frac{2N}{N-2}}\,dx, \quad R>>1,
$$
by taking $R=T^\frac1{1+\min\{b_1,b_2\}}$, then we get that
\begin{align} \label{1}
\frac 1 T \int_0^T\int_{|x|<T^\frac1{1+\min\{b_1,b_2\}}}|u|^{\frac{2N}{N-2}}\,dx\,dt\lesssim T^{-\frac{\min\{b_1,b_2\}}{1+\min\{b_1,b_2\}}}.
\end{align}
Now by the mean value theorem and  \eqref{1}, then \eqref{2} holds true. This completes the proof.
\end{proof}

The main ingredient for the proof of Theorem \ref{scattering} is the following scattering criterion in the spirit of arguments in \cite{Ta1}.

\begin{prop}\label{crt}
Let $u \in C([0, +\infty), H^1(\R^N))$ be the global solution to the Cauchy problem \eqref{equ} with $u_0 \in \mathcal{A}^+_{\omega}$. Assume that
$$
0<\sup_{t\geq0}\|u(t)\|_{H^1}:=E<\infty.
$$
If there exist $R, \epsilon>0$ depending on $E, p_i, b_i$ and $N$ such that
\begin{equation}\label{crtr}
\liminf_{t\to+\infty}\int_{|x| \leq R}|u(t,x)|^2\,dx<\epsilon^2,
\end{equation}
then $u$ scatters forward in time.
\end{prop}

To establish Proposition \ref{crt}, we need the following essential lemmas.

\begin{lem} \label{crt1}
Under the assumptions of Proposition \ref{crt}, there exist $T,\gamma>0$ such that
\begin{equation}\label{lm}
\left\|e^{\textnormal{i}(t-T)\Delta}u(T)\right\|_{\Lambda^{\max\{s_{c,1},s_{c,2}\}}(T,+\infty)}\lesssim \epsilon^\gamma,
\end{equation}
where
$$
s_{c, j}:=\frac{N}{2}-\frac{2-b_j}{p_j-2}, \,\,\, j=1,2.
$$
\end{lem}
\begin{proof}
For simplicity, we shall assume that $s_{1,c}<s_{2,c}$. Let  $\omega,\beta>0$ be two positive real numbers to be fixed later. Taking account of Lemma \ref{str}, we know that there exists $T_0>\epsilon^{-\beta}$ such that
\begin{equation}\label{fre}
\left\|e^{\textnormal{i}t\Delta}u_0\right\|_{\Lambda^{s_{c,2}}(T_0,+\infty)} \leq \epsilon^\omega.
\end{equation}
Let us now take the time slabs $J_1:=[0,T-\epsilon^{-\beta}]$ and $J_2:=[T-\epsilon^{-\beta},T]$, where $T>T_0$ is chosen later. The integral formula gives that
\begin{align} \label{int'}
\begin{split}
&e^{\textnormal{i}(t-T)\Delta}u(T)\\
&=e^{\textnormal{i}(t-T)\Delta}\left(e^{\textnormal{i}T\Delta}u_0-\textnormal{i}\int_0^Te^{\textnormal{i}(T-s)\Delta}\left(|x|^{-b_1}|u|^{p_1-2}u \right)\,dt+\textnormal{i}\int_0^Te^{\textnormal{i}(T-t)\Delta}\left(|x|^{-b_2}|u|^{p_2-2}u\right)\,ds\right)\\
&:=e^{\textnormal{i}t\Delta}u_0-\textnormal{i}\int_0^Te^{\textnormal{i}(t-s)\Delta}\mathcal N_1(u)\,dt+\textnormal{i}\int_0^Te^{\textnormal{i}(t-s)\Delta}\mathcal N_2(u)\,ds\\
&=e^{\textnormal{i}t\Delta}u_0-\textnormal{i} \int_{J_1} \left(e^{\textnormal{i}(t-s)\Delta}\left(\mathcal N_1(u)-\mathcal N_2(u)\right)\right)\,ds -\textnormal{i} \int_{J_2} \left(e^{\textnormal{i}(t-s)\Delta}\left(\mathcal N_1(u)-\mathcal N_2(u)\right)\right)\,ds\\
&:=e^{\textnormal{i}t\Delta}u_0-\textnormal{i}\mathcal F_1-\textnormal{i}\mathcal F_2.
\end{split}
\end{align}
In view of \eqref{crtr}, there exists $T>T_0$ such that
\begin{equation}\label{3.34}
\int_{\R^N}\chi_R|u(T)|^2\,dx<\epsilon^2, \quad R>>1,
\end{equation}
where $\chi_R$ is given by \eqref{chi-R}. On the other hand, we see that
\begin{equation}\label{3.35}
\frac\partial{\partial t}\int_{\R^N}\chi_R|u(t)|^2\,dx \leq2\left |\int_{\R^N}\left(\nabla\chi_R\cdot\nabla u\right)\bar u\,dx\right|\lesssim\frac 1 R.
\end{equation}
Integrating \eqref{3.35} in time on $[t, T]$ for $t \in J_2$ and using \eqref{3.34}, we then have that, for any $R>\epsilon^{-(2+\beta)}$,
\begin{align}\label{3.36}
\begin{split}
\int_{\R^N}\chi_R|u(t)|^2\,dx
&\leq \int_{\R^N}\chi_R|u(T)|^2\,dx+\frac{\epsilon^{-\beta}}R \\
&\lesssim \epsilon^2+\frac{\epsilon^{-\beta}}R \\
&\lesssim \epsilon^2.
\end{split}
\end{align}
By using Strichartz's estimates in Lemma \ref{str} and Sobolev's embedding inequality, we can write that
\begin{align}\label{3.38}
\begin{split}
\|\mathcal F_2\|_{\Lambda^{s_{c,2}}(T,+\infty)}
&\lesssim \left\|\int_{J_2}e^{\textnormal{i}(t-s)\Delta}\mathcal N_1\,ds\right\|_{\Lambda^{s_{c,2}}(T,+\infty)}+\left\|\int_{J_2}e^{\textnormal{i}(t-s)\Delta}\mathcal N_2\,ds\right\|_{\Lambda^{s_{c,2}}(T,+\infty)}\\
&\lesssim \sum_{j=1}^2\|\mathcal N_j\|_{\Lambda^{-s'_{c,2}}(J_2)}.
\end{split}
\end{align}
Making use of Lemma \ref{NLEst} and \eqref{3.36},  H\"older's inequality and Sobolev's embedding inequality, we then get that, for some $2<r<2^*$ and $0<\lambda_1, \lambda_2<1$,
\begin{align}\label{003.38}
\begin{split}
\|\mathcal N_1\|_{\Lambda^{-s'_{c,2}}(J_2)}
&\lesssim |J_2|^{\theta_1}\|u\|_{L^\infty(I, L^r)}^{p_1-2}\|\chi_R\,u\|_{L^\infty(J_2, L^r)}+|J_2|^{\theta_2}\|u\|_{L^\infty(J_2, L^{p_1})}^{p_1-2}\|\chi_R\,u\|_{L^\infty(J_2, L^{p_1})}\\
&+R^{-b_1}|J_2|^{\theta_2}\|u\|_{L^\infty(J_2, L^{p_1})}^{p_1-1}\\
&\lesssim \epsilon^{-\beta\theta_1}\|\chi_R\,u\|_{L^\infty(J_2, L^2)}^{\lambda_2}+\epsilon^{-\beta\theta_2}\|\chi_R\,u\|_{L^\infty(J_2, L^2)}^{\lambda_1}+R^{-b_1}\epsilon^{-\beta\theta_2}\\
&\lesssim \epsilon^{\lambda_2-\beta\theta_1}+\epsilon^{\lambda_1-\beta\theta_2}+R^{-b_1}\epsilon^{-\beta\theta_2}\\
&\lesssim \epsilon^\omega,
\end{split}
\end{align}
where in the last line it is sufficient to pick
$$
0<\omega<\min\{\lambda_1,\lambda_2\},\quad 0<\beta<\min\left\{\frac{\lambda_2-\omega}{\theta_1},\frac{\lambda_1-\omega}{\theta_2}\right\}, \quad R>\epsilon^{-\frac{\omega+\beta\theta_2}{b_1}}.
$$
Similarly, arguing as in \eqref{003.38}, we can also derive that $\|\mathcal N_2\|_{\Lambda^{-s'_{c,2}}(J_2)} \lesssim \epsilon^\omega$.
This together with \eqref{003.38} and \eqref{3.38} indicates that $\|\mathcal F_2\|_{\Lambda^{s_{c,2}}(T,+\infty)}\lesssim \epsilon^\omega$.

Now in the spirit of \cite{CC}, we shall take $(q_j, r_j)\in \Lambda^{s_{c,j}}$ for $j=1, 2$, a small real number $\nu>0$ and $(q^j,r^j)\in\Lambda$ such that
\begin{align*}
\frac1{q^j}=\frac{\frac1{q_j}-\nu s_{c,j}}{1-s_{c,j}}, \quad
\frac1{r^j}=\frac{1}{1-s_{c,j}}\left(\frac1{r_j}-\frac{s_{c,j}(N-2-4\nu)}{2N}\right).
\end{align*}
By the integral Duhamel's formula \eqref{int'}, we can write that
\begin{align*}
\mathcal F_1=e^{\textnormal{i}t\Delta}\left(e^{\textnormal{i}(\epsilon^{-\beta}-T)\Delta}u(T-\epsilon^{-\beta})-u_0\right).
\end{align*}
Using H\"older's inequality, Strichartz's estimates in Lemma \ref{str} and the dispersive estimate
\begin{align*}
\|e^{\textnormal{i}t\Delta}\varphi \|_r\lesssim |t|^{-N(\frac12-\frac 1r)}\|\varphi \|_{r'},\quad r\geq2,\; t\neq 0,
\end{align*}
we then have that
\begin{align*} 
\|\mathcal F_1\|_{L^{q_j}((T,+\infty),L^{r_j})}
&\leq \|\mathcal F_1\|_{L^{q^j}((T,+\infty),L^{r^j})}^{1-s_{c,j}}\|\mathcal F_1\|_{L^{\frac1\nu}((T,+\infty),L^{\frac{2N}{N-2-4\nu}})}^{s_{c,j}}\\
&=\left\|e^{\textnormal{i}t\Delta}\left(e^{\textnormal{i}(\epsilon^{-\beta}-T)\Delta}u(T-\epsilon^{-\beta})-u_0\right)\right\|_{L^{q^j}((T,+\infty),L^{r^j})}^{1-s_{c,j}}\|\mathcal F_1\|_{L^{\frac1\nu}((T,+\infty),L^{\frac{2N}{N-2-4\nu}})}^{s_{c,j}}\\
& \lesssim \|\mathcal F_1\|_{L^{\frac1\nu}((T, +\infty),L^{\frac{2N}{N-2-4\nu}})}^{s_{c,j}}\\
&\lesssim \sum_{j=1}^2\left\|\int_{J_1}|t-s|^{-(1+2\nu)}\|\mathcal N_i\|_{\frac{2N}{N+2+4\nu}}\,ds\right\|_{L^\frac1\nu(T,+\infty)}^{s_{c,j}}\\
&\lesssim \sum_{j=1}^2\|u\|^{(p-1)s_{c,j}}_{L^{\infty}H^1}\left\|(t-T-\epsilon^{-\beta})^{-2\nu}\right\|_{L^\frac1\nu(T,+\infty)}^{s_{c,j}}\\
&\lesssim \epsilon^{\beta\nu}.
\end{align*}
This apparently shows that $\|\mathcal F_1\|_{\Lambda^{s_{c,2}}(T,+\infty)}\lesssim \epsilon^{\beta\nu} $. Choose $\gamma=\min\{\omega, \epsilon^{\beta\nu}\}$, then the proof is completed.
\end{proof}

\begin{lem}\label{smlglb}
Let $u\in C([0, +\infty), H^1(\R^N))$ be the global solution to the Cauchy problem \eqref{equ} with $u_0 \in \mathcal{A}^+_{\omega}$. If there exists $\varepsilon>0$ such that
\begin{align} \label{3.58}
\sum_{j=1}^2\left\|e^{\textnormal{i}(t-T)\Delta}u(T)\right\|_{\Lambda^{s_{c,j}}(T, + \infty)}&\leq \varepsilon,
\end{align}
for some $T>0$, then there holds that
\begin{align*} 
\sum_{j=1}^2\|u\|_{\Lambda^{s_{c,j}}(T,+\infty)}&\leq 2\sum_{j=1}^2\left\|e^{\textnormal{i}(t-T)\Delta}u(T)\right\|_{\Lambda^{s_{c,j}}(T,+\infty)},
\end{align*}
\begin{align*} 
\|\left\langle \nabla\right\rangle u\|_{\Lambda(T, +\infty)}&\lesssim \|u(T)\|_{H^1},
\end{align*}
where $\left\langle \xi \right\rangle:=(1+|\xi|^2)^{\frac 12}$.
\end{lem}
\begin{proof}
For $a,b>0$,  we define the space
$$
Y_T:=\left\{v \in C([0, +\infty), H^1(\R^N)) : \sum_{j=1}^2\|v\|_{\Lambda^{s_{c,j}}(T,+\infty)}\leq a, \|\left\langle \nabla\right\rangle v\|_{\Lambda(T,+\infty)}\leq b\right\}
$$
equipped with the complete distance
$$
d(v, w):=\sum_{j=1}^2\|v-w\|_{\Lambda^{s_{c,j}}(T, +\infty)}.
$$
Let also define the integral functional by
\begin{align}\label{int}
\mathcal{F}(v)&:=e^{\textnormal{i}(t-T)\Delta}u(T)+\textnormal{i}\int_T^t e^{\textnormal{i}(t-s)\Delta}\left(\mathcal N_1(v)-\mathcal N_2(v)\right)\,ds,
\end{align}
where $\mathcal{N}_j(v)=|x|^{-b_j}|v|^{p_j-2}v$ for $j=1,2$. By means of Strichartz's estimates in Lemma \ref{str} and Sobolev's embedding inequality, we have that
\begin{align} \label{strf}
\|\mathcal{F}(v)\|_{\Lambda^{s_{c,2}}(T, +\infty)} \nonumber
& \lesssim \left\|e^{\textnormal{i}(t-T)\Delta}u(T)\right\|_{\Lambda^{s_{c,2}}(T, +\infty)}+\left\|\int_T^te^{\textnormal{i}(t-s)\Delta}\mathcal N_1(v)\,ds\right\|_{\Lambda^{s_{c,2}}(T, +\infty)} \\ \nonumber
& \quad +\|\mathcal N_2(v)\|_{\Lambda^{-s'_{c,2}}(T,+\infty)}\\
&\lesssim \left\|e^{\textnormal{i}(t-T)\Delta}u(T)\right\|_{\Lambda^{s_{c,2}}(T,+\infty)}+\left\|\left\langle \nabla\right\rangle \int_T^te^{\textnormal{i}(t-s)\Delta}\mathcal N_1(v)\,ds\right\|_{\Lambda(T, +\infty)}\\ \nonumber
& \quad +\|\mathcal N_2(v)\|_{\Lambda^{-s'_{c,2}}(T, +\infty)} \\ \nonumber
&\lesssim \left\|e^{\textnormal{i}(t-T)\Delta}u(T)\right\|_{\Lambda^{s_{c,2}}(T, +\infty)}+\left\|\left\langle \nabla\right\rangle \mathcal N_1(v)\right\|_{\Lambda'(T, +\infty)}+\|\mathcal N_2(v)\|_{\Lambda^{-s'_{c,2}}(T, +\infty)}.
\end{align}
Taking into account \cite[Lemma 2.7]{Cam}, we then get that, for certain $0<\omega_1, \omega_2<<1$,
\begin{align*} 
\|\mathcal{F}(v)\|_{\Lambda^{s_{c,2}}(T, +\infty)}
&\leq \left\|e^{\textnormal{i}(t-T)\Delta}u(T)\right\|_{\Lambda^{s_{c,2}}(T, +\infty)} +c\|v\|_{L^\infty H^1}^{\omega_1}\|v\|_{\Lambda^{s_{c,1}}(T, +\infty)}^{p_1-2-\omega_1}\|\left\langle \nabla\right\rangle v\|_{\Lambda(T, +\infty)}\\
& \quad +c \|v\|_{L^\infty H^1}^{\omega_2}\|v\|_{\Lambda^{s_{c,2}}(T, +\infty)}^{p_2-1-\omega_2}.
\end{align*}
Similarly, we are able to derive that
\begin{align*}
\|\mathcal{F}(v)\|_{\Lambda^{s_{c,1}}(T, +\infty)}
& \leq \left\|e^{\textnormal{i}(t-T)\Delta}u(T)\right\|_{\Lambda^{s_{c,1}}(T,+\infty)}+c\|v\|_{L^\infty H^1}^{\omega_2}\|v\|_{\Lambda^{s_{c,2}}(T,+\infty)}^{p_2-2-\omega_2}\|\left\langle \nabla\right\rangle v\|_{\Lambda(T,+\infty)}\\
& \quad + c\|v\|_{L^\infty H^1}^{\omega_1}\|v\|_{\Lambda^{s_{c,1}}(T,+\infty)}^{p_1-1-\omega_1}.
\end{align*}
Moreover, we can obtain that
\begin{align} \label{363}
\|\left\langle \nabla\right\rangle \mathcal{F}(v)\|_{\Lambda(T, +\infty)}
& \leq c\|u(T)\|_{H^1}+ c\left(\sum_{j=1}^2\|v\|_{L^\infty H^1}^{\omega_j}\|v\|_{\Lambda^{s_{c,j}}(T,+\infty)}^{p_j-2-\omega_j}\right)\|\left\langle \nabla\right\rangle v\|_{\Lambda(T,+\infty)}.
\end{align}
In light of Strichartz's estimates in Lemma \ref{str}, then
\begin{align*}
\|\mathcal{F}(v)-\mathcal{F}(w)\|_{\Lambda(T, +\infty)}
& \leq c \sum_{j=1}^2\left\||x|^{-b_j}(|v|^{p_j-2}+|w|^{p_j-2})|v-w|\right\|_{\Lambda'(T, +\infty)}\\
& \leq c \sum_{j=1}^2\left(\|v\|_{L^\infty H^1}^{\omega_j}\|v\|_{\Lambda^{s_{c,j}}(T, +\infty)}^{p_j-2-\omega_j}+\|w\|_{L^\infty H^1}^{\omega_j}\|w\|_{\Lambda^{s_{c,j}}(T, +\infty)}^{p_j-2-\omega_j}\right)\|v-w\|_{\Lambda(T,+\infty)}.
\end{align*}
Consequently, we have that
\begin{align*}
\sum_{j=1}^2\|\mathcal{F}(v)\|_{\Lambda^{s_{c,j}}(T,+\infty)}
&\leq \sum_{j=1}^2\left\|e^{\textnormal{i}(t-T)\Delta}u(T)\right\|_{\Lambda^{s_{c,j}}(T, +\infty)}+c\left(\sum_{j=1}^2a^{p_j-2-\omega_j}\,b^{\omega_j}\right)(a+b) ,
\end{align*}
\begin{align*}
\|\left\langle \nabla\right\rangle \mathcal{F}(v)\|_{\Lambda(T, +\infty)}
&\leq c\|u(T)\|_{H^1}+c\left(\sum_{j=1}^2a^{p_j-2-\omega_j}\,b^{\omega_j}\right)b,
\end{align*}
\begin{align*}
d(\mathcal{F}(v), \mathcal{F}(w))&\leq c\left(\sum_{j=1}^2a^{p_j-2-\omega_j}\,b^{\omega_j}\right)d(v,w).
\end{align*}
From \eqref{3.58}, we shall take $T>0$ such that
\begin{align*}
a:=2\sum_{j=1}^2\|e^{\textnormal{i}(t-T)\Delta}u(T)\|_{\Lambda^{s_{c,j}}(T,+\infty)}\label{368}\leq 2\varepsilon.
\end{align*}
Then we define $ b:=2c\|u(T)\|_{H^1}$. Thereby, for $\varepsilon>0$ small enough,  $F$ is a contraction mapping on $Y_T$ and the proof is achieved by the classical Picard arguments. This completes the proof.
\end{proof}

\begin{lem}\label{smlsct}
Let $u \in C([0, +\infty), H^1(\R^N))$ be the global solution to the Cauchy problem \eqref{equ} with $u_0 \in \mathcal{A}^+_{\omega}$. If there exists $\varepsilon>0$ such that
\begin{align*}
\left\|e^{\textnormal{i}(t-T)\Delta}u(T)\right\|_{\Lambda^{\max\{s_{c,1},s_{c,2}\}}(T,+\infty)}\lesssim \varepsilon
\end{align*}
for some $T>0$, then $u$ scatters forward in time.
\end{lem}
\begin{proof}
We may assume without loss of generality that $s_{c,1}<s_{c,2}$.  Let $(q_j, r_j)\in \Lambda^{s_{c,j}}$ and
$$
\lambda:=\frac{s_{c,1}}{s_{c,2}}=\frac{(N(p_2-2)-2(2-b_2))(p_1-2)}{(N(p_1-2)-2(2-b_1)(p_2-2)} \in (0,1).
$$
Let $ \rho ,\gamma>0$ be such that
\begin{align*}
\frac{1}{q_1}= \frac{\lambda}{q_2}+\frac{1-\lambda}{\rho}, \quad \frac{1}{r_1}= \frac{\lambda}{r_2}+\frac{1-\lambda}{\gamma}.
\end{align*}
It is clear that $(\rho, \gamma) \in \Lambda$. From H\"older's inequality, we then get that
 \begin{align*}
\left\|e^{\textnormal{i}(t-T)\Delta}u(T)\right\|_{\Lambda^{s_{c,1}}(T, +\infty)}
&\leq \left\|e^{\textnormal{i}(t-T)\Delta}u(T)\right\|_{\Lambda^{s_{c,2}}(T, +\infty)}^\lambda\left\|e^{\textnormal{i}(t-T)\Delta}u(T)\right\|_{\Lambda(T,+\infty)}^{1-\lambda}\\
&\lesssim \varepsilon^{\lambda}.
\end{align*}
This infers that
$$
\sum_{j=1}^2\left\|e^{\textnormal{i}(t-T)\Delta}u(T)\right\|_{\Lambda^{s_{c,j}}(T,+\infty)} \lesssim \varepsilon^{\lambda}+\varepsilon.
$$
Observe that, by \eqref{363} and Lemma \ref{smlglb},
\begin{align*}
\left\|e^{-\textnormal{i}t\Delta}u(t)-e^{-\textnormal{i}t'\Delta}u(t')\right\|_{H^1}
&\lesssim \left(\sum_{j=1}^2\|u\|_{L^\infty H^1}^{\omega_j}\|u\|_{\Lambda^{s_{c,j}}(t, t')}^{p_j-2-\omega_j}\right)\left\|\left\langle \nabla\right\rangle u\right\|_{\Lambda(t, t')}\\
&\to 0\quad\mbox{as} \,\, t, t'\to +\infty.
\end{align*}
This completes the proof.
\end{proof}

\begin{proof} [Proof of Proposition \ref{crt}]
Invoking Lemmas \ref{crt1} and \ref{smlsct}, we then complete the proof.
\end{proof}

\begin{proof} [Proof of Theorem \ref{scattering}]
The proof consists in a standard application of Lemma \ref{global1} and Proposition \ref{crt}.
\end{proof}

\hrule

\vspace{0.3cm}

{\noindent {\bf\large Acknowledgements.} { The authors express their sincere gratitude to the referee for the insightful comments and suggestions that have enhanced the manuscript.}\newline  T. Gou was supported by the Postdoctoral Science Foundation of China (No.2021M702620), the National Natural Science Foundation of China (No.12101483 $\&$ 12471113) and the Shaanxi Fundamental Science Research Project for Mathematics and Physics (No. 22JSZ003).}

\vspace{0.3cm}

{\noindent{\bf\large Declarations.}}
On behalf of all authors, the corresponding author states that there is no conflict of interest. No data-sets were generated or analyzed during the current study.

\vspace{0.3cm}

\hrule

\end{document}